\documentclass[11pt, a4paper, twoside]{book} 
\usepackage[hidelinks]{hyperref} 

\hypersetup{
    pdftitle={Supply Chain Analytics: A Data-Driven Approach},
    pdfauthor={Elioth Sanabria},
    pdfsubject={Graduate Textbook - Columbia IEOR},
    pdfkeywords={Supply Chain, Analytics, Optimization, Forecasting},
  pdfstartview={FitH}} 
\usepackage{fancyhdr}
\usepackage[centertags]{amsmath}
\usepackage{amsfonts}
\usepackage{amsthm}
\usepackage{bm}
\usepackage{latexsym,amssymb}
\usepackage{graphics}
\usepackage{subcaption} 
\usepackage[demo]{graphicx} 
\usepackage{hyperref}
\usepackage{tabularx}

\usepackage{algorithm} 
\usepackage{algorithmicx}
\usepackage{algpseudocode}

\usepackage[
backend=bibtex,
maxbibnames=99  
]{biblatex}

\usepackage{caption}
\usepackage{tikz}
\usepackage{subcaption}
\usetikzlibrary{graphs, positioning, calc, arrows.meta,3d,positioning, fit, backgrounds,matrix,decorations.pathreplacing}

\usepackage{pgfplots}
\pgfplotsset{compat=1.17}
\usepgfplotslibrary{fillbetween,groupplots}
\usepackage{pgfplotstable}

\usepackage[
a4paper,
twoside,             
textwidth=6.25in,
textheight=8.75in,
top=1.0in,           
bottom=1.0in,        
inner=1.0in,         
outer=1.0in,         
]{geometry}

\newcommand{\pr}{\mbox{\sf P}}
\newcommand{\ex}{{\bf\sf E}}               
\newcommand{\var}{\mbox{\sf Var}}
\newcommand{\cov}{\mbox{\sf Cov}}

\newcommand{\calm}{{\cal M}}

\newcommand{\cald}{{\cal D}}
\newcommand{\calp}{{\cal P}}

\newcommand{\cals}{{\cal S}}

\newcommand{\calh}{{\cal H}}

\newcommand{\calx}{{\cal X}}

\theoremstyle{definition} 
\newtheorem{exm}{Example}[chapter]

\def\eqd{{\buildrel {\rm d} \over =}}

\def\ind{{\bf 1}}


\title{Supply Chain Analytics: A Data-Driven Approach}
\author{Elioth Sanabria}
\date{Spring 2026} 

\begin{document}
\maketitle

\cleardoublepage

\tableofcontents

\cleardoublepage
\pagenumbering{arabic}

\linespread{1.5}	
\setlength{\parskip}{0.5em}
\baselineskip 18pt

\newcommand{\mbe}{\mathbb{E}}

\newcommand{\Int}{\mathbb{Z}}





\chapter{Demand Estimation and Forecasting} 
\label{sec:demand}
\section{Introduction}

Demand forecasting is an important tool to make decisions such as financial projections, capacity planning and allocation of resources for a productive activity. The demand for any product, service or good in the future is a random variable. That is, it is unknowable before it is realized. The demand that already occurred is called \emph{realized demand} or can also be called historical demand.

In fact, only the past is really known, everything happening in the future has to be guessed. In these lectures notes we deal with the most common practices in \emph{guessing} the future quantities of the demand. Mathematically, it is useful to model the demand as a random variable. Random variables are a mathematical construct to capture events and their randomness. We will cover the basics building up the necessary knowledge to manipulate and eventually model and estimate demand in the future.

\section{The Demand as a Random Variable}
\label{sec:demand_rv}
Suppose we want to model the demand quantity once in the future (say, next week or next month or next year). We denote this random quantity by the  capital letter $D$. This random variable $D$ is in fact a function that takes as input a random mechanism $\omega$ and outputs an event, in our case a number for the demand\footnote{The space of random mechanism is normally called sample space, and is denoted as the set $\Omega$. Then, formally the random variable is a function mapping $D:\Omega\rightarrow E$ to a set of possible events $E$.}. For example, the demand can be a number like $D=14$ with probability $p_{14}$, or it could be $D=5,000$ with a probability $p_{5,000}$.

As normally the demand is a \emph{count} (of products, goods or services) we will assume the demand cannot take negative values. Then, we can describe the probability of any outcome of the demand by a \emph{probability mass function} or pmf as $p_i=\pr(D=i)$ for all possible values $i={0,1,\dots}$ For the pmf to be correctly defined two things must happen: First, each outcome must occur with a probability between 0 and 1, that is $0\le p_i\le 1$ and adding over all outcomes should add up to 1, that is, $\sum_{i=0}^\infty p_i=1$.

With this framework, we can already solve interesting questions: For example, if $D$ is the quantity of semiconductor chips sold in the next quarter, we can ask \emph{what is the probability that the demand is at least 10 million chips?}. To answer this, the event $\{D>10^7\}$ is the same as the events $\{D=10^7+1,D=10^7+2,...\}$, then we can simply do:
\begin{eqnarray}
	\pr(D>10^7)=\pr(D=10^7+1,D=10^7+2,\dots)=\sum_{i=10^7+1}^\infty p_i,
\end{eqnarray}

Throughout the topic it is often useful to get used to manipulate event sets. For example, in the previous example we exploited the fact that the event $\{D>10^7\}$ can be written as the union that the demand takes exactly each value greater than $10^7$, noting that each of these events is disjoint (it can be a value or another but not two at the same time). We could have also taken an alternative approach by calculating the probability of the \emph{complement} event, that is, that the demand is less or equal to $10^7$ and then, given that the probability of all outcomes is $1$, subtract $P(D\le 10^7)$. In other words:
\begin{eqnarray}
	\pr(D>10^7) =  1 - \pr(D\le10^7).
\end{eqnarray}
For a given number $i$, the probability $F_D(i)=\pr(D\le i )$ is known as the \emph{cumulative mass function} or cmf. The usefulness of the cmf, is related again to event calculus, for example asking the probability that the demand is between 1 and 10 million units. That is, $\pr(10^6< D\le 10^7)$ can be expressed in terms of the cmf as $F_D(10^7)-F_D(10^6)$.

Another very important quantity is known as the mean, or average. This comes from averaging all possible outcomes for the demand weighted by their probabilities. This quantity is denoted as $\mu_{D}=\ex(D)=\sum_{i=0}^\infty ip_i=0p_0+1p_1+2p_2+\dots$\footnote{Related to event calculus, the mean $\mu_{D}$ is also equal to $\sum_{i=0}^\infty \bar F_D(i)$, where $\bar F_D(i)=1-F_D(i)$. Why is this the case?} Averages are important for demand estimation, because they provide a scenario that is centered around the more likely outcomes for the demand. For example, when making a resource allocation decision, such as building a factory or a data-center, the capacity should comfortably handle the average demand for the facility. But then, the next immediate question is how much slack capacity should the factory have? 

The variance of a random variable is a measure of how much the random variable \emph{deviates} from its mean $\mu_D$. The variance is defined as $\sigma^2_D=\var(D)=\ex[(D-\mu_D)^2]$, the $\ex$ operator (or function) denotes the operation of multiplying each possibly outcome by its probability. So, $\ex[(D-\mu_D)^2]$ can be computed as $\sum_{i=0}^\infty(i-\mu_D)^2p_i$. The square root of the variance is known as the \emph{standard deviation} as is denoted by the symbol $\sigma_D$.

\subsection{Useful Probability Calculations}
So far, we have reviewed basic concepts in probability that can answer some questions about the random outcome of the demand. We can go one step further and with some knowledge of the distribution, answer other questions: For example, say that you only know the mean $\mu_D$ of the demand, and they ask you, if the average demand is 1 million units, what's the probability that the demand is less than 5 million? 

A useful relation known as \emph{Markov's Inequality} states that $\pr(D\ge a)\le \mu_D/a$, replacing $a=5\times10^6$ tells us that the probability that the demand exceeds 5 million units is at most  $\mu_D/(5\times10^6)$. If the average demand $\mu_D$ is 1 million, the probability that it is at least 5 million is at most  $\mu_D/a=10^6/(5\times10^6)=1/5=20\%$. Therefore, the probability that the demand is less than 5 million is at least $80\%$ (could be more). That being said, this estimate is rather loose as the inequality is \emph{at most} the value. In the previous example, it could be the case that the probability that the production exceeds 5 million is exactly 0.

Another useful inequality is known as \emph{Chebyschev's Inequality} stating that $\pr(|\mu_D-D|\ge k\sigma_D )\le1/k^2$. Going back to our previous example, suppose that additionally to knowing that the average demand $\mu_D$ is 1 million, we also know that its standard deviation $\sigma_D$ is 100,000 units. Applying this inequality, the probability that the demand is below 800,000 or above 1,200,000 is at most 25\% (using the inequality with $k=2$). Going back to the previous question about the probability that the demand exceeds 5 million is at most $1/40^2=0.06\%$. Even more can be said, for example the probability that the demand exceeds 2 million is less than $1/100=1\%$. The moral of this, is that with just the mean $\mu_D$ and the standard deviation of the demand $\sigma_D$ there is already a lot of useful information that can be inferred about the demand without knowing the full probability distribution of all its outcomes.

\subsection{Side Information}
\label{sec:side_info}

The demand is clearly influenced by other variables (random or deterministic): the economic cycle, interest rates, competitors prices, the demand of other products, weather and so on. Having a framework to understand how a variable affects the demand is crucial to better estimate and anticipate the demand. 

Suppose you are in charge of the production of the Airpods wireless earphones. Intuitively, the demand of Airpods is \emph{positively} related to the demand of other products in the iOS ecosystem (such as laptops, iPhones, etc) and \emph{negatively} to competitor products, such as Android headphones. Moreover, the demand for Airpods should not be affected by the demand of apples or oranges.

Ideally, it would be great to have a mathematical function $g$, that given an outcome of a related random variable $X=x$ gives us the best estimate of the random demand $D$. You can imagine $X$ as the demand for iPhones in the AirPods example (as a significant portion of AirPods users own iPhones). Then, we would like to find the function $g$ that in average gives the closest prediction of the demand $D$. One such formulation is

\begin{eqnarray}
\min_g \ex \left[(g(X)-D)^2\right],
\end{eqnarray}

This way of thinking is a strategy that we will repeatedly use in the exposition of the topic. That is, to ask what's the best possible data-driven way to make guesses about the random nature of the demand. Assuming $X$ can take $x_1,x_2,\dots$ different values. The expected value $\ex \left[(g(X)-D)^2\right]$, taking $p_{ix}=\pr(D=i,X=x)$ can be written as:
\begin{eqnarray}
	\ex \left[(g(X)-D)^2\right]=\sum_{i=0}^\infty\sum_{x=x_1,\dots}(g(x)-i)^2p_{ix},
\end{eqnarray}
The strategy to solve a seemingly complex problem like this is to break down a complex problem into smaller and easier problems to solve (a well known paradigm in Computer Science known as divide-and-conquer). Then, the first observation is that the solution only depends on the values of the random variable $X$. That is, for a given $X=x$, the function $g(x)$ is a single value guessing the most likely outcome of the demand $D$. Taking only a single value of $X=x$ gives the expression $\sum_{i=0}^\infty(g(x)-i)^2p_{ix}$. This expression has a single unknown $g(x)$, which is a convex function\footnote{A sum of convex functions is convex.}. Taking derivatives and equating to 0 leads to the following expression for $g(x)$:
\begin{eqnarray*}
	g(x)=\sum_{i=0}^\infty i \frac{p_{ix}}{p_x}=\ex(D|X=x),
\end{eqnarray*}
This expression intuitively says, that the best predictor for the demand $D$ when the side variable $X$ is equal to $x$ is the expected value of the demand given that $X=x$. This is the so-called \emph{conditional expectation}, which can be seen as simply taking the well known expected value with new probabilities $p_{ix}/p_x$ representing the likelihood of each outcome of the demand given a certain outcome for the side variable $X$. The denominator $p_x$ scales the new probabilities $p_{ix}$ of the demand so that they add up to 1, because instead of considering all outcomes of $X$, only the case where $X=x$ is considered.

\begin{exm}[{\bf Hydroelectric Plant}]
Suppose you are managing a hydroelectric plant supplying electric energy to a major city. $X$ represents whether or not it is an extremely hot day ($X=1$ means it is extremely hot, and $X=0$ means it is not), extremely hot days occur rarely, say $\pr(X=1)=p_h=5\%$. But when they occur, the demand for electricity is very high (due to the usage of Air Conditioner and other appliances). The conditional distribution of $(D|X)$ is given by $\pr(D=5\text{Mw}|X=0)=\pr(D=15\text{Mw}|X=1)=1$. Given this distribution, what is the average demand $\mu_D$ of the power plant? The average demand of the power plant is:
\begin{eqnarray}
	\mu_D=\ex[D]=\ex[\ex[D|X]]=5\text{Mw}(1-p_h)+15\text{Mw}p_h=5.5\text{Mw}.
\end{eqnarray}
This relationship clearly shows one of the main ideas of this course: While the demand in average is relatively stable at a level of 5.5Mw most days, with a small probability it can jump to almost a 300\% increase if it's a very hot day. If the power plant does not have the capacity to service the demand, it can be a catastrophic event for the city, since many crucial services including hospitals and other vital supply chains rely on electricity. Hedging against rare events is the whole point of a resilient supply chain. See Figure \ref{fig:energy_demand_final}. 

\begin{figure}[htbp]
    \centering
    \begin{tikzpicture}
        \begin{axis}[
            width=14cm, 
            height=8cm,
            xlabel={Time (Days)},
            ylabel={Energy Consumption (Mw)},
            xmin=0.5, xmax=9.5,
            ymin=0, ymax=18, 
            xtick={1,2,3,4,5,6,7,8,9},
            ytick={0, 5, 10, 15},
            grid=both,
            grid style={dashed, gray!30},
            legend pos=north west, 
            legend cell align={left},
            legend style={nodes={scale=0.9, transform shape}}, 
        ]
            \addplot[
                color=blue,
                mark=circle,
                mark size=2.5pt,
                line width=1.2pt,
            ] coordinates {
                (1,5)(2,5)(3,5)(4,5)(5,5)(6,15)(7,5)(8,5)(9,5)
            };
            \addlegendentry{Realized Demand ($D$)}
            
            \addplot[
                color=red,
                dashed,
                thick,
                domain=0.5:9.5
            ] {5.5};
            \addlegendentry{Average Demand ($\mu_D = 5.5$)}
            
            \draw[stealth-, thick] (axis cs:6,15) -- (axis cs:7,16) 
                node[anchor=west] {Extremely Hot Day ($X=1$)};
                
            \node[anchor=north, font=\small, text=gray] at (axis cs:3, 4.8) {Normal Days ($X=0$)};

        \end{axis}
    \end{tikzpicture}
    \caption{Energy consumption over time. The spike represents a rare event where temperature ($X$) affects the demand ($D$).}
    \label{fig:energy_demand_final}
\end{figure}
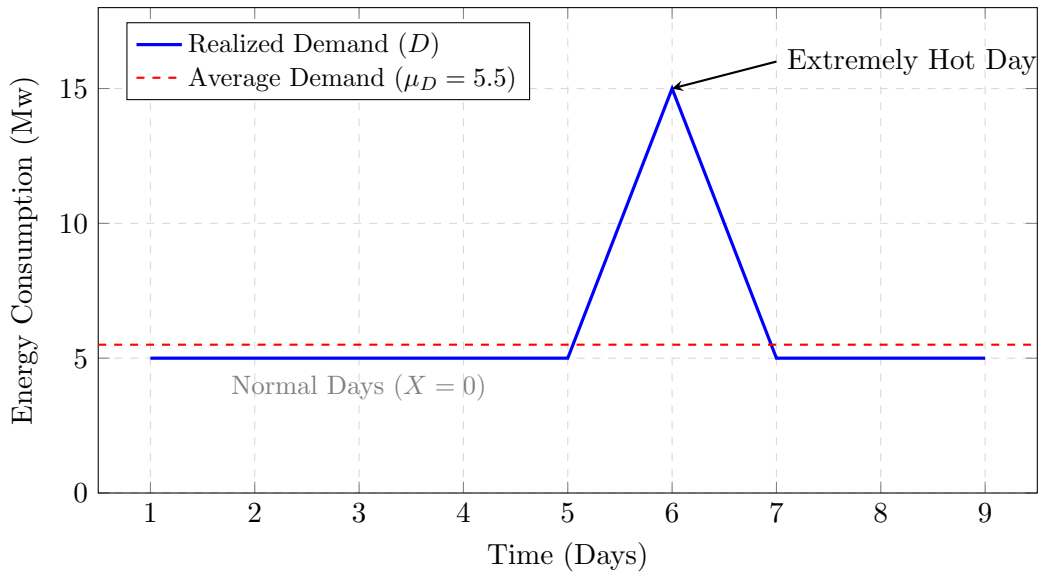

\end{exm}

\section{Demand Evolution over Time}
\label{sec:demand_modeling}

So far, we have presented a static view of the demand as random variable. While useful as a first order tool of analysis, it is easy to imagine that the demand of many products is in fact time-dependent: Ice cream sales peak in the summer, airline traffic peaks around holidays and restaurants demand has modes around lunch and dinner. Similar to last example of the previous section, averaging without side information can generate misleading estimates for the demand. Therefore, if the demand is time-dependent it makes sense to include time in its modeling.

Time is represented as a variable $t=0,1,2,\dots$ where each value of $t$ represents a moment in time. The unit of time can be taken as seconds, hours, days or years. Depending on the context of the specific supply chain decision: For example, in analysing the demand behaviour in an online store, it makes sense to have relatively small units of time. Conversely, when planning capacities or resources for facilities, it makes sense to consider the demand over years. Later it will also become evident that there is a trade-off between the historical amount of data available and the quality of the estimation of the demand.

Incorporating time into the demand estimation leads to indexing the demand random variable by time as $D_t$, denoting the random demand at time $t$. In most scenarios, we have access to the demand up to some point in time in the past. For example, suppose that you are given the monthly historical demand for the last 5 years. Chronologically, these are 60 datapoints $d_0,d_1,d_2,\cdots,d_{59}$. Data like this can be denoted compactly as $\{d_s\}_{s=0}^t$ (imagine a spreadsheet with a single column and a row for each entry of the demand at every date). The distinction between the capital and lowercase symbols for the demand denote the random variable $D_t$ and certain outcome for it $d_t$. This is important, as even though the past is a realization of the demand, the outcome was a result of the probabilistic nature of the demand (for example, a month of low demand could be caused by an unexpected shock, and it might provide information for future outcomes of the demand).

Random variables that evolve over time are referred as \emph{Stochastic Processes}. We present a couple of examples to get a feel for stochastic processes. One of the simplest models is the following:
\begin{eqnarray}
	\label{eq:base}
	D_{t+1}=\alpha+\beta D_t+V_t,
\end{eqnarray}
We call this the \emph{Base Model} as it will allow us to model some common behaviours of the demand over time. In this model $\alpha,\beta\in\mathbb{R}$ are two constants while $V_t$ is an independent random variable. In plain words, the equation says that the demand at time $t+1$ is equal to a constant $\alpha$ plus a term multiplying the previous period demand $D_t$, plus a random quantity $V_t$ (normally this quantity expresses what's commonly referred as an error or simply unexplained variation of the demand). This model is related to the previous one (independent of time) when $\beta=0$\footnote{If $V_t$ is a random variable so is $\alpha+V_t$. In fact, constants are random variables with their probability mass concentrated at a single value.}.

\begin{figure}[ht!]
    \centering
    
    \begin{minipage}{0.48\textwidth}
        \centering
        \begin{tikzpicture}
            \begin{axis}[
                width=\textwidth,
                height=5.5cm,
                grid=major,
                title={Mean Reverting ($\beta = 0.7$)},
                xlabel={Time ($t$)},
                ylabel={Demand ($D_t$)},
                xmin=0, xmax=30,
                ymin=0, ymax=130,
                title style={font=\bfseries}
            ]
            \addplot[blue, thick] coordinates {
                (0,100)(1,81.3)(2,69.8)(3,66.4)(4,53.0)(5,46.1)(6,47.9)(7,41.2)(8,35.5)(9,43.1)(10,44.8)(11,39.0)(12,28.5)(13,25.2)(14,36.6)(15,41.5)(16,33.8)(17,38.6)(18,42.1)(19,32.8)(20,37.1)(21,34.5)(22,28.9)(23,37.2)(24,39.4)(25,29.5)(26,23.9)(27,34.1)(28,35.8)(29,42.1)(30,45.1)
            };
            \addplot[red, dashed, thick] coordinates {(0, 33.3) (30, 33.3)};
            \node[red, anchor=south] at (axis cs: 15, 33.3) {\footnotesize Mean $\approx 33.3$};
            \end{axis}
        \end{tikzpicture}
    \end{minipage}
    \hfill
    \begin{minipage}{0.48\textwidth}
        \centering
        \begin{tikzpicture}
            \begin{axis}[
                width=\textwidth,
                height=5.5cm,
                grid=major,
                title={Random Walk ($\beta = 1.0$)},
                xlabel={Time ($t$)},
                ylabel={Demand ($D_t$)},
                xmin=0, xmax=30,
                ymin=50, ymax=450,
                title style={font=\bfseries}
            ]
            \addplot[green!60!black, thick] coordinates {
                (0,100)(1,113.8)(2,131.2)(3,138.8)(4,135.5)(5,131.1)(6,149.3)(7,165.7)(8,162.2)(9,168.6)(10,183.1)(11,190.4)(12,187.9)(13,205.1)(14,219.0)(15,220.3)(16,237.1)(17,245.5)(18,259.9)(19,270.2)(20,285.5)(21,292.1)(22,309.8)(23,325.3)(24,336.8)(25,350.1)(26,362.4)(27,361.0)(28,377.2)(29,390.1)(30,405.8)
            };
            \end{axis}
        \end{tikzpicture}
    \end{minipage}

    \vspace{0.8cm}

    \begin{minipage}{0.48\textwidth}
        \centering
        \begin{tikzpicture}
            \begin{axis}[
                width=\textwidth,
                height=5.5cm,
                grid=major,
                title={Explosive ($\beta = 1.15$)},
                xlabel={Time ($t$)},
                ylabel={Demand ($D_t$)},
                xmin=0, xmax=30,
                ymin=0, ymax=11000,
                title style={font=\bfseries}
            ]
            \addplot[red!80!black, thick] coordinates {
                (0,100)(1,120.3)(2,153.1)(3,189.9)(4,225.0)(5,264.4)(6,321.1)(7,375.5)(8,440.0)(9,519.8)(10,612.4)(11,708.9)(12,825.1)(13,955.5)(14,1111.4)(15,1282.9)(16,1485.1)(17,1711.2)(18,1982.5)(19,2291.1)(20,2640.4)(21,3042.8)(22,3509.8)(23,4032.5)(24,4652.8)(25,5360.1)(26,6171.4)(27,7099.9)(28,8171.2)(29,9398.5)(30,10811.2)
            };
            \end{axis}
        \end{tikzpicture}
    \end{minipage}

    \caption{Stochastic demand simulations ($\alpha=10, D_0=100, V_t \sim \text{U}(-15, 15)$) illustrating the impact of $\beta$. The \textbf{Mean Reverting} ($\beta=0.7$) case centers around $\alpha/(1-\beta)$, the \textbf{Random Walk} ($\beta=1$) exhibits non-stationary drift, and the \textbf{Explosive} ($\beta=1.15$) case follows a steep exponential trajectory.}
    \label{fig:demand_simulations}
\end{figure}
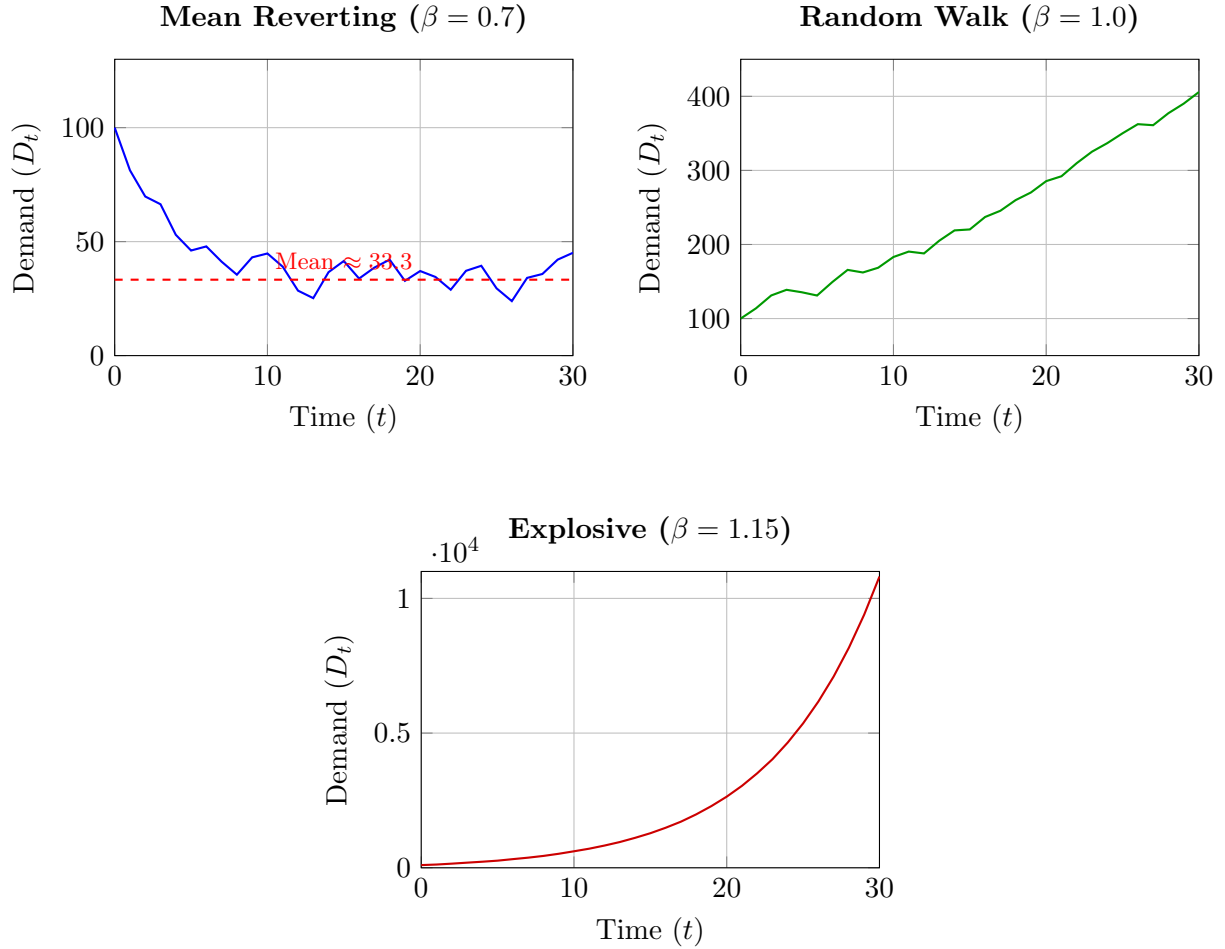

To see the dependence over time, assume that at time $0$, there is an initial demand $D_0\ge 0$. $D_1=\alpha+\beta D_0+V_0$, and $D_2=\alpha+\beta D_1+V_1=\alpha+\beta (\alpha+\beta D_0+V_0)+V_1=\alpha+\beta\alpha+V_1+\alpha V_0+\beta^2D_0$. Continuing this iterative process leads to the following expression:
\begin{eqnarray}
	D_{t+1}=\sum_{j=0}^{t-1}(\alpha+V_{t-j})\beta^{j}+\beta^t D_0,
\end{eqnarray}
Assuming $V_t$ for $t=0,1,\dots$ has identical and independent distribution with $\ex(V_t)=0$. Taking expectations yields the following:
\begin{eqnarray}
\ex(D_{t+1}) =\alpha \frac{1-\beta^t}{1-\beta}+\beta^t\ex(D_0),
\end{eqnarray}
The first term of the right hand side comes from the finite sum of a geometric series\footnote{Take $s_t = 1 + \beta + \cdots+ \beta^t$ and $\beta s_t = \beta+\dots+\beta^{t+1}$, subtracting the two equations yields the result shown.}. Then, the expected value of the demand depends on the expected value of the initial demand, plus a term that converges to $\alpha/(1-\beta)$ at $t$ goes to $\infty$. Moreover, this convergence only occurs if $|\beta|<1$. If $|\beta|>1$ the demand explodes as time goes by (think why this is the case). When $\beta=1$ with $\pr(V_t= 1)=\pr(V_t= -1)=1/2$. The model retrieved is called a \emph{Random Walk}. Among other interesting properties, random walks drift away indefinitely to $\pm \infty$ and can also be negative. See Figure \ref{fig:demand_simulations} to see the different models. With a small modification, of taking $\ln (D_t)$ instead of $D_t$ a \emph{Geometric Random Walk} can be used to model the demand.

\subsection{Modeling Seasonality and Time-varying Factors}

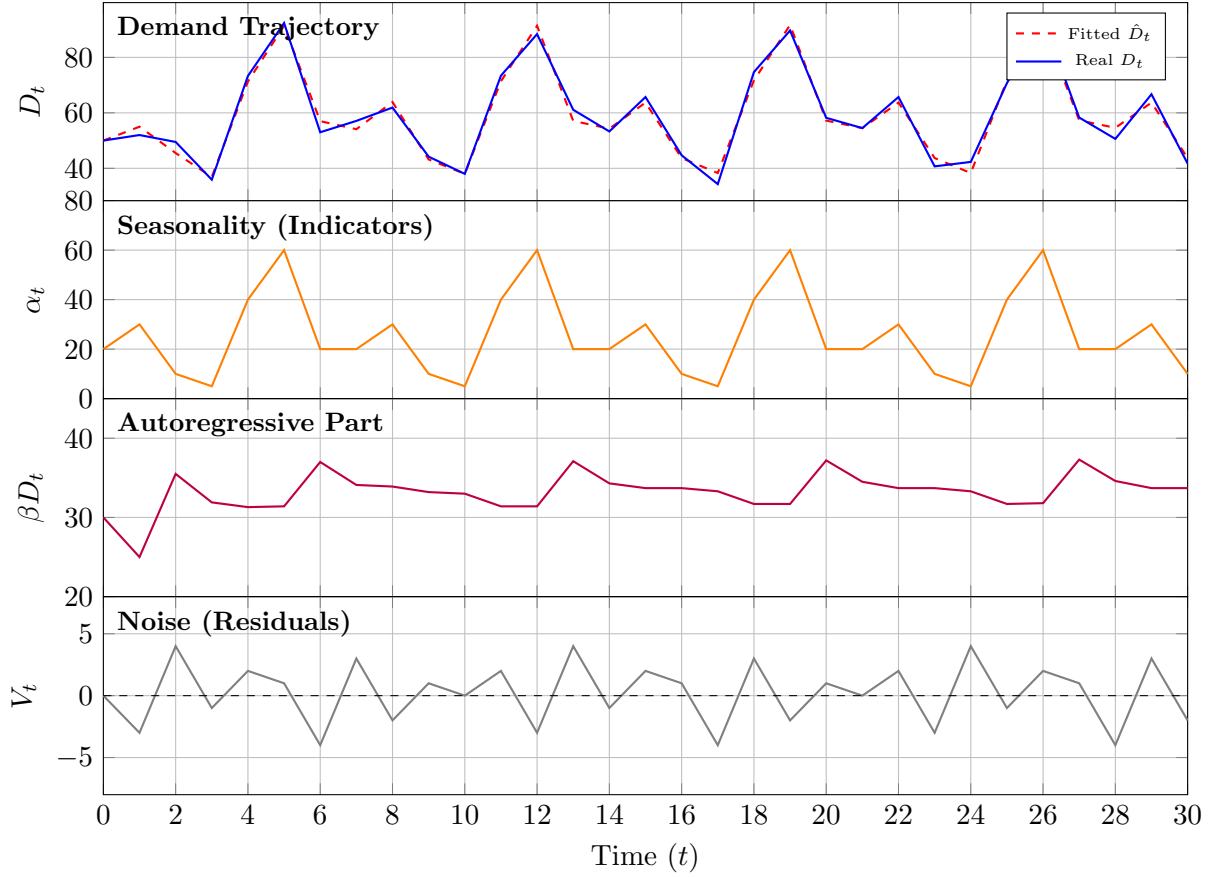
\begin{figure}[ht!]
    \centering
    \begin{tikzpicture}
        \begin{groupplot}[
            group style={
                group size=1 by 4,
                vertical sep=0pt, 
                x descriptions at=edge bottom,
            },
            width=\textwidth,
            height=4.2cm,
            xmin=0, xmax=30,
            grid=major,
            every axis title/.style={below right, at={(0,1)}, font=\bfseries\small, inner sep=5pt}
        ]

        \nextgroupplot[
            ylabel={$D_t$}, 
            title={Demand Trajectory},
            legend style={at={(0.98,0.95)}, anchor=north east, font=\tiny}
        ]
            \addplot[red, dashed, thick] coordinates {
                (0,50.0)(1,55.0)(2,45.5)(3,36.9)(4,71.3)(5,91.4)(6,57.0)(7,54.1)(8,63.9)(9,43.2)(10,38.0)(11,71.4)(12,91.4)(13,57.1)(14,54.3)(15,63.7)(16,43.7)(17,38.3)(18,71.7)(19,91.7)(20,57.2)(21,54.5)(22,63.7)(23,43.7)(24,38.3)(25,71.7)(26,91.8)(27,57.3)(28,54.6)(29,63.7)(30,43.7)
            };
            \addlegendentry{Fitted $\hat{D}_t$}
            
            \addplot[blue, thick] coordinates {
                (0,50.0)(1,52.0)(2,49.5)(3,35.9)(4,73.3)(5,92.4)(6,53.0)(7,57.1)(8,61.9)(9,44.2)(10,38.0)(11,73.4)(12,88.4)(13,61.1)(14,53.3)(15,65.7)(16,44.7)(17,34.3)(18,74.7)(19,89.7)(20,58.2)(21,54.5)(22,65.7)(23,40.7)(24,42.3)(25,70.7)(26,93.8)(27,58.3)(28,50.6)(29,66.7)(30,41.7)
            };
            \addlegendentry{Real $D_t$}

        \nextgroupplot[ylabel={$\alpha_t$}, title={Seasonality (Indicators)}, ymin=0, ymax=80]
            \addplot[orange, thick] coordinates {
                (0,20)(1,30)(2,10)(3,5)(4,40)(5,60)(6,20)(7,20)(8,30)(9,10)(10,5)(11,40)(12,60)(13,20)(14,20)(15,30)(16,10)(17,5)(18,40)(19,60)(20,20)(21,20)(22,30)(23,10)(24,5)(25,40)(26,60)(27,20)(28,20)(29,30)(30,10)
            };

        \nextgroupplot[ylabel={$\beta D_t$}, title={Autoregressive Part}, ymin=20, ymax=45]
            \addplot[purple, thick] coordinates {
                (0,30.0)(1,25.0)(2,35.5)(3,31.9)(4,31.3)(5,31.4)(6,37.0)(7,34.1)(8,33.9)(9,33.2)(10,33.0)(11,31.4)(12,31.4)(13,37.1)(14,34.3)(15,33.7)(16,33.7)(17,33.3)(18,31.7)(19,31.7)(20,37.2)(21,34.5)(22,33.7)(23,33.7)(24,33.3)(25,31.7)(26,31.8)(27,37.3)(28,34.6)(29,33.7)(30,33.7)
            };

        \nextgroupplot[ylabel={$V_t$}, xlabel={Time ($t$)}, title={Noise (Residuals)}, ymin=-8, ymax=8]
            \addplot[gray, thick] coordinates {
                (0,0.0)(1,-3.0)(2,4.0)(3,-1.0)(4,2.0)(5,1.0)(6,-4.0)(7,3.0)(8,-2.0)(9,1.0)(10,0.0)(11,2.0)(12,-3.0)(13,4.0)(14,-1.0)(15,2.0)(16,1.0)(17,-4.0)(18,3.0)(19,-2.0)(20,1.0)(21,0.0)(22,2.0)(23,-3.0)(24,4.0)(25,-1.0)(26,2.0)(27,1.0)(28,-4.0)(29,3.0)(30,-2.0)
            };
            \addplot[black, thin, dashed] coordinates {(0,0) (30,0)};
        \end{groupplot}
    \end{tikzpicture}
    \caption{Demand decomposition into different components: Seasonal, autoregressive, and noise. Corresponding to each term in the specification of the demand.}
    \label{fig:demand_decomposition}
\end{figure}

Another common feature of the demand is \emph{seasonality}. Seasonality is defined as the fluctuation of the demand caused by the passage of time or certain time-dependent events (such as holidays or weather variations). Think of holidays and demand that changes due to weather or temperature. The most common way to model this is to follow the method in the weather example of the previous section, by adding an extra variable codifying the side-information of interest. In the case of seasonality, this can be achieved by adding binary variables that are equal to 1 when $t$ is equal to that time. For example, think of daily data where we can to capture day of week effect on the demand (think of a restaurant or retail business). A way to achieve this, is to modify the $\alpha$ term in Equation (\ref{eq:base}) as:
\begin{eqnarray}
	\alpha_t = \sum_{d=0}^{s-1}\ind\{t\,\text{mod}\,s=d\}\gamma_d,
\end{eqnarray}
In this case, the indicator $\ind\{t\,\text{mod}\,s=d\}$ is only one every $s$ observations. In the days of the week example, $s=7$ as there are 7 days in a week. Putting this expression back in the base model yields:
\begin{eqnarray}
	D_{t+1}=\sum_{d=0}^{s-1}\ind\{t\,\text{mod}\,s=d\}\gamma_d+\beta D_t+V_t,
\end{eqnarray}
The model is quite similar to the original one (that had only one constant term $\alpha$ while this one has a constant term for every day of the week). In general, we can also include time-varying factors, such as seasons (also encoded in indicators or numerical variables. A clean way to include these factors is to flatten these variables into a vector $h_t$. For example, letting $h_t=(1,t)'$ retrieves a model with a constant plus a linear trend allows to express the original $\alpha_t$ term as $\alpha_t=A h_t$ where $A$ in this case is a vector of parameters with the same dimension as the number of time-varying variables, that is $A\in\mathbb{R}^{|h_t|}$. This allows to re-write the base model as:
\begin{eqnarray}
	D_{t+1}=\alpha_t+\beta D_t+V_t=A h_t+\beta D_t+V_t,
\end{eqnarray}
This instance of the base model allows for arbitrary specifications of time-dependent variables. See Figure \ref{fig:demand_decomposition} for an example of the decomposition of the demand into different components.

\subsection{Other Modeling Specifications}
It is easy to imagine that a similar strategy can be employed to add arbitrary random variables or other stochastic processes $X_t$ of variables that are related to the demand $D_t$. In general, one can imagine that any arbitrary real-valued function can be calibrated to solve the following problem:
\begin{eqnarray}
	\min_g\ex[(g(D_{t},D_{t-1}\dots,X_t,X_{t-1},\dots)-D_{t+j})^2].
\end{eqnarray}
It should be evident that this is exactly the same problem in Section [\ref{sec:side_info}] where the optimal solution is $g(D_{t},D_{t-1}\dots,X_t,X_{t-1},\dots)=\ex(D_{t+j}|D_{t},D_{t-1}\dots,X_t,X_{t-1},\dots)$. Then, speaking strictly from an optimization point of view, what matters is that the family of functions $g$ can express the conditional expectation of the demand $D_{t+j}$ as a function of the other variables. Whether the model is a simple linear regression or a sophisticated neural network, what matters is that it can capture the shape of the conditional expectation. As we will see next, there is another tension related to the fact that in general the probabilistic distribution of the demand is unknown and only a finite sample of the variables has to be used for estimation.

\section{Demand Estimation in Practice}

So far, we have built a methodological theory to study the evolution of the demand over time. The core message is that the best possible guess for the demand is the conditional expectation of the demand taking into account information provided by other relevant variables. In practice, the situation is significantly different: At best there are limited observations of the historical demand (or no information at all in the worst case). The side-variables that affect the demand are unknown or unobservable and in many cases they are endogenous (competitors set prices reactively to their own demand and compete against each other for market domination). The latter situation is more complex and will be studied in a later chapter. In this section, we focus on the case where we have limited data and want to predict the demand in the \emph{short-term}. Here, short-term is used loosely as the foreseeable future.

The main estimation problem for the demand is that its true probabilistic distribution is unknown. Then, the goal is to learn the demand from \emph{sample} observations.  Suppose, a historical yearly demand of salmon in a country. This demand should be stable and more or less proportional to the population of the country, to account for the population fluctuation, consider the \emph{per capita} demand. This historical demand is represented by the observations $\{d_s\}_{s=1}^{t}$. Estimating the average demand from samples is calculated as:
\begin{eqnarray}
	\mu_D(t) = \frac{1}{t}\sum_{s=1}^td_s,
\end{eqnarray}
This average converges to a stable (constant value) in various cases: For example, when the demand samples are \emph{independent and identically distributed}, or i.i.d. Meaning that they come from the same probabilistic distribution, and the previous values do not affect the probability of new outcomes of the demand. One can see getting observations of the demand as sampling the outcome of an experiment repeatedly and measuring the demand after each trial. 

It is important to determine how reliable this estimate is, especially when managing a supply chain as making resource allocation decisions in the future can be extremely costly. Probability theory has two answers regarding how trusted this estimator can be: The \emph{Law of Large Numbers} and the \emph{Central Limit Theorem}.

The Law of Large Numbers (LLN) states that for a large enough $t$, the average demand (sampled i.i.d. with finite mean and variance) $\mu_D(t)$ converges to the true average value $\mu_D$. That is, $\mu_D(t)\rightarrow \mu_D$ as $t\rightarrow \infty$. As mentioned before, this is important because there is a probabilistic guarantee that sampling the demand i.i.d. and averaging its values over time recovers the true average demand.

The Central Limit Theorem (CLT) goes a step further and provides a probability distribution for the randomness of the demand. Given an i.i.d. sample of the demand $\{d_s\}_{s=1}^{t}$. The random variable defined by $\frac{\sqrt{t}(\mu_D(t)-\mu_D)}{\sigma_D}$ converges in distribution to a random variable $Z$ distributed as a standard normal random variable\footnote{A standard normal random variable has probability density function (pdf) $f(x)=e^{-x^2/2}/\sqrt{2\pi}$. A pdf is the continuous analogous to the pmf in the sense that for a random variable $X$ with density $f(x)$ the probability $\pr[X\in (x,x+\delta_x)]=f(x)\delta_x$ for a small $\delta_x>0$.}.


\begin{exm}[{\bf Rice Storage Levels}]
The government of Japan asks you to conduct a study to determine what level of rice reserves they should hold. Their first question is how many years $t$ of data they should look at to get a 95\% confidence estimate of the average demand $\mu_D$ that is within 10\% of the true value, using the sample estimate $\mu_D(t)$. Suppose you are told that the standard deviation of the demand $\sigma_D$ is 30\% of the mean $\mu_D$. That is, $\sigma_D=\mu_D30\%$.

Using the CLT we get that the sample average demand is equal in distribution to:
\begin{eqnarray}
	\mu_D(t)\,\eqd\,\frac{\sigma_D}{\sqrt{t}}Z+\mu_D,
\end{eqnarray}

Dividing by $\mu_D$ yields the percentage error  $\frac{\sigma_D}{\mu_D\sqrt{t}}Z+1$. Where $\frac{\sigma_D}{\mu_D\sqrt{t}}Z$ is the error in terms of the demand, which can be positive or negative. Then, for this expression to within 10\% of the true value with 95\% confidence we need the following to be true:
\begin{eqnarray}
\pr\left(\left\lvert\frac{\sigma_D}{\mu_D\sqrt{t}}Z\right\rvert\le10\%\right)=95\%,
\end{eqnarray}
Using the fact that $\sigma_D/\mu_D=30\%$, the probability can be rearranged to $\pr(|Z|\le \sqrt{t}/3)=1-\pr(|Z|> \sqrt{t}/3)=95\%$ or $\pr(Z>\sqrt{t}/3)=5\%/2=2.5\%$. That is, the tail cdf of the standard normal distribution $\bar\Phi(\sqrt{t}/3)=2.5\%$. Applying the inverse\footnote{An inverse function $f^{-1}(x)$ is a function such that $f^{-1}(f(x))=x$. For example, the inverse function of $f(x)=e^x$ is $f^{-1}(x)=\ln(x)$.} tail cdf yields:
\begin{eqnarray}
	\sqrt{t}/3=\bar\Phi^{-1}(2.5\%)=1.96.
\end{eqnarray}
With this, we get our answer that the number of observations $t$ needed for the sample average to be within 10\% of the mean is at least $t=(3(1.96))^2=34.5$ years.
\end{exm}

In general, the above example is the best possible case. In probability textbooks, a commonly encountered rule is that the sample size for estimation should be at least 30 observations. As both the LLN and the CLT relied on the i.i.d. property of the samples for convergence.

To illustrate this, consider a random demand that follows the following demand process $D_{t+1}=D_t+c+\sigma Z$, with $D_0=0$. Note that this process is an instance of the base model with $\alpha=c$, $\beta=1$ and $V_t=\sigma Z$. In this case, $c+\sigma Z$ is called a normal random variable with mean $c$ and variance $\sigma^2$ (likewise, a normal random variable can be transformed in a standard normal by subtracting its mean and dividing by $\sigma$. This operation is useful to make computations always in terms of standard normal random variables).

The demand process is equal to $D_{t}=ct+\sigma\sum_{i=1}^tZ_t$. The expected value of the demand $\ex(D_{t})=tc$, since each independent standard normal $Z_t$ has mean 0. Likewise, the sample mean (extracted from historical observations of the demand) will be approximately this expected value, that is, $\mu_D(t)\approx tc$. As $t$ increases indefinitely, so will the sample mean, that is, $\mu_D(t)\rightarrow\infty$ as $t\rightarrow\infty$. So, how can we estimate the demand in these cases?

Most textbooks in Time Series spend a significant amount of their exposition in studying \emph{stationary} stochastic processes. That is, processes that converge to a stable mean in the long-term and do not explode indefinitely. The textbook solution is to apply a transformation to the original series, such that the transformation is stationary. In the previous example, a solution is to consider a new series taking the differences of the demand. For example, consider the difference process $\Delta_t=D_t-D_{t-1}=c+\sigma Z$, which as discussed previously is a normal random variable with mean $c$ and variance $\sigma^2$. The sample mean of historical observations of $\Delta_t$ should be the same as sampling standard normal observations. 

\subsection{Data-driven Demand Estimation}
In practice we cannot expect to know beforehand the actual transformation that makes the demand stationary, since we only have access to a historical demand series $\{d_s\}_{s=1}^{t}=\{d_1,d_2,\dots,d_t\}$ that is sampled from an unknown stochastic process (with unknown mean, variance, seasonality and functional form).

The strategy that we will implement consists of the following steps:
\begin{enumerate}
\item Propose a functional form rich enough to capture the demand within a well-defined family of functional forms.
\item Calibrate the parameters of the function that fit the demand closest via optimization.
\item Evaluate the performance of the proposed function out-of-sample.
\end{enumerate}

For Step 1, we can use some of the modeling techniques that we have discussed in Section \ref{sec:demand_modeling}. For example, we can use for example the model presented in the previous section, $\hat D_{t+1}=\hat A h_t+\hat\beta D_t$. Here, the hats above the expression denote quantities that are calibrated from data. Recalling the definition of conditional expectation in Section \ref{sec:side_info}, a reasonable approach to calculate the parameters for the demand is to solve the following optimization problem:
\begin{eqnarray}
\label{eq:pop_demand}
	\min_{\hat A,\hat\beta} \ex [(\hat D_{t}(\hat A,\hat\beta)-D_{t})^2],
\end{eqnarray}
This expression can be seen as an instance of the conditional expectation defined as the best predictor for a random variable. The solution of this problem, is in fact an approximation of the conditional expectation of the demand, that is, $\hat D_t(A^*,\beta^*)\approx \ex [D_t|D_{t-1},\dots]$. The quality of the approximation depends crucially on two factors: That the functional form $\hat D_t$ allows to approximate the conditional expectation. And also, on the sample size $t$ of historical observations for the demand. Moreover, note that when the true parameters of the demand $A$ and $\beta$ are used, the difference $\hat D_t(A,\beta)-D_t\,\eqd \,V_t$. In general, $\ex (V_t)$ should be equal to 0 (in terms of the base model, if $V_t$ has non-zero mean $c$, define new variables $\alpha'=\alpha+c$ and $V'_t=V_t-c$, then $\ex(V'_t)=0$). Therefore, this gives us another clue that $\ex(\hat D_t(A^*,\beta^*)-D_t)$ should be approximately $0$ when solving the optimization problem of fitting the demand.

To solve the optimization problem in Step 2, we solve the data-driven analogous of the problem in (\ref{eq:pop_demand}) with a sample of demand observations $\{d_s\}_{s=0}^t$ (in this case with the particular choice of demand function in Step 1):
\begin{eqnarray}
\min_{\hat A,\hat\beta} \frac{1}{t}\sum_{s=1}^{t}(\hat A h_{s-1}+\hat\beta d_{s-1}-d_{s})^2,
\end{eqnarray}
This problem is a convex optimization program with a unique solution, that can be solved with any off-the-shelf quadratic optimization solver\footnote{Linear least squares is an instance of Quadratic Programming. Note that in matrix form linear least squares has the objective $||\beta X-y||_2^2=2(\beta XX^T\beta/2-\beta Xy^T+yy^T/2)$ which is a quadratic program by letting $Q=XX^T/2$ and $c^T=- Xy^T$. Then, $\min_\beta ||\beta X-y||_2^2=2\min_\beta(\beta Q\beta^T+\beta c^T) + yy^T$.}. 

It is appealing to use a family of functions that can express very rich functional forms (for example, for a closed interval it is well-known that a Fourier basis can approximate any continuous function). The formalism would be to express the fitted demand $\hat D_t(\theta)$ as a function of parameters $\theta$, for example, the weights of a deep neural network, or the collection of splits in a random forest. The approach would be the same, that is, to solve the optimization problem $\min_{\theta} \ex [(\hat D_{t}(\theta)-D_{t})^2]$. 

The issue with this approach is that if the family of functions is too rich, the sample minimization would be too optimistic (potentially 0 error as $\hat D_t(\theta^*)=d_t$, for all observations in the sample $\{d_s\}_{s=1}^t$) and the performance for future predictions could be arbitrarily off. The two main reasons for this phenomenon are: First, in practice we only have access to \emph{finite} samples that might not  represent the full geometry of the conditional distribution (a well-known phenomenon in statistical analysis). Second, real-life stochastic process are in general changing in distribution over time (a so-called \emph{doubly stochastic process}, where the original process changes distribution based on another unknown random variable). In the context of demand estimation, this could be a shock to the demand or a technological change. For example, the demand for horse carriages catastrophically decreased with the invention and mass production of mechanical engine cars (think that the measure of power of engines is called \emph{horsepower}) in the early twentieth century. Likewise, the introduction of new competitors, or shocks to the demand or supply can change the trajectory of the demand.

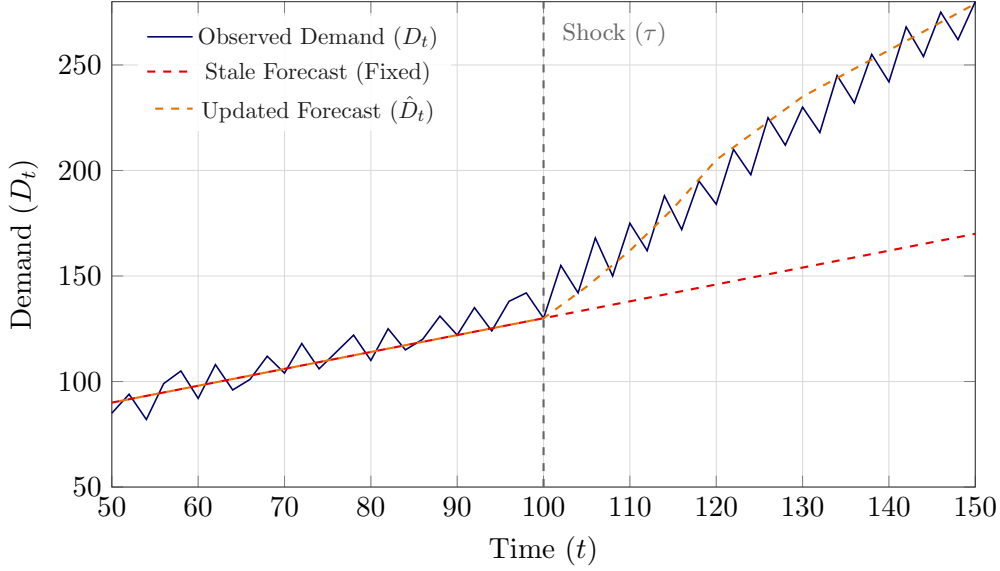
\begin{figure}[ht]
    \centering
    \begin{tikzpicture}
    \begin{axis}[
        xlabel={Time ($t$)},
        ylabel={Demand ($D_t$)},
        legend pos=north west,
        grid=both,
        grid style={line width=.1pt, draw=gray!15},
        major grid style={line width=.2pt, draw=gray!30},
        width=13cm, height=8cm,
        xmin=50, xmax=150, 
        ymin=50, ymax=280,
        axis line style={black!80},
        legend style={nodes={scale=0.8}, fill=white, fill opacity=0.9, draw=none, row sep=2pt}
    ]

    \addplot[blue!40!black, semithick] coordinates {
        (50,85) (52,94) (54,82) (56,99) (58,105) (60,92) (62,108) (64,96) (66,101) (68,112)
        (70,104) (72,118) (74,106) (76,114) (78,122) (80,110) (82,125) (84,115) (86,120) (88,131)
        (90,122) (92,135) (94,124) (96,138) (98,142) (100,130) 
        (102,155) (104,142) (106,168) (108,150) (110,175) (112,162) (114,188) (116,172) (118,195) (120,184)
        (122,210) (124,198) (126,225) (128,212) (130,230) (132,218) (134,245) (136,232) (138,255) (140,242)
        (142,268) (144,254) (146,275) (148,262) (150,280)
    };
    \addlegendentry{Observed Demand ($D_t$)}

    \addplot[red!90!black, thick, dash pattern=on 3pt off 3pt] coordinates {
        (50,90) (100,130) (150,170)
    };
    \addlegendentry{Stale Forecast (Fixed)}

    \addplot[orange!90!black, thick, dash pattern=on 3pt off 3pt, dash phase=3pt] coordinates {
        (50,90) (75,110) (100,130) 
        (105,145) (110,162) (115,182) (120,205) (125,220) (130,235) 
        (135,246) (140,257) (145,268) (150,279)
    };
    \addlegendentry{Updated Forecast ($\hat{D}_t$)}

    \draw[black!60, dashed, thick] (axis cs:100,50) -- (axis cs:100,280);
    \node[anchor=north west, font=\footnotesize, color=black!60] at (axis cs:101,275) {Shock ($\tau$)};

    \end{axis}
    \end{tikzpicture}
    \caption{The impact of structural shocks on forecast divergence. Prior to the shock ($t<100$), both models track the same pre-shock parameters; an alternating dash pattern (red/orange) is used to display both distinct forecasting paths. After the shock ($\tau$), the stale model (red) fails to account for the change, while the updated model (orange) recalibrates and captures the new demand trajectory.}
    \label{fig:fixed_forecast}
\end{figure}

While this problem cannot be fully mitigated, if there is enough data we can use what we have learned to come up with a partial solution: For the first problem of over-fitting the demand, we can sub-divide our full sample of historical observation into multiple sub-samples where by setting up an experiment, we can use the LLN and the CLT to get an idea of the performance of our demand estimation strategy out-of-sample. For the second problem, imagine that distribution changes at random with a throw of a coin with small probability $p$, then, the time $\tau$ in between changes in distribution is in average $\ex(\tau)=1/p$. Then, if instead of training our model with all of the historical observations, we train the model with $k$ observations, such that $k<\ex(\tau)$, with probability $(1-p)^k\approx 1$ we are in the current distribution of the demand (another approach would be to weight the observations with a decaying weight). To see the point of using $k$ observations imagine a sudden financial crisis that impacts the demand of luxury goods, such as designer bags. Using the whole historical demand will make the estimates to continue growing if the historical data includes most of the period before the crisis. Using the last $k$ observation makes the estimates adjust faster to the new distribution of the demand. For another example, imagine a positive shock to the demand, in Figure \ref{fig:fixed_forecast} a fixed model cannot update to changes in the distribution of the demand.

Then, for Step 3 we take the following approach: Given historical observations of the demand $\{d_s\}_{s=1}^T$, a estimation window size $k$ and a prediction horizon $j$. Define the sample optimization problem as $\ell(t,t')=1/(t'-(t+1))\sum_{s=t+1}^{t'}(A h_{s-1}+\beta d_{s-1}-d_{s})^2$. Likewise, define the absolute error $|\varepsilon_{t+j}(A,\beta)|=|\hat D_{t+j}(A,\beta)-d_{t+j}|$. The strategy is going to be to estimate multiple observations of the error out-of-sample, and average them as a measure of how the model performs out of sample. The pseudo-code routine to perform the \emph{backtest} can be written as:

\begin{algorithm}
	\caption{Backtesting Demand Prediction}
	\label{alg:findmax}
	\begin{algorithmic}[1]
		\State \textbf{Input:} Historical observations for the demand $d_1,\dots,d_T$. A training window $k$, and a prediction horizon $j$. And a model specification for the demand $\hat D_{t+j}(\cdot)$
		\State \textbf{Output:} The average out-of-sample absolute error $\ex(|\varepsilon_{t+j}|)$
		\Procedure{Backtest}{$\{d_s\}_{s=1}^T,k,j$}
		\State $\mu_{|\varepsilon|}=0$
		\For{$t = k$ \textbf{to} $T-j$}
		\State $A^*,\beta^* \gets \text{argmin} \,\,\ell(t-k+1,t)$ 
		\State $|\varepsilon_{t+j}(A^*,\beta^*)|=|\hat D_{t+j}(A^*,\beta^*)-d_{t+j}|$
		\State $\mu_{|\varepsilon|} = \mu_{|\varepsilon|}+|\varepsilon_{t+j}(A^*,\beta^*)|/(T-k-j)$
		\EndFor
		\State \textbf{return} $\mu_{|\varepsilon|}$
		\EndProcedure
	\end{algorithmic}
\end{algorithm}
With this procedure, given two models (say the linear one that we have been exploring) and a more sophisticated one like a Deep Neural Network model. The estimates of the mean absolute error $\mu_\varepsilon$ can be used to compare the performance of the model out-of-sample. This procedure is commonly used in finance to backtest algorithms. The key concept is that the measurement of the error is done outside of the training data of the model, thus, mimicking the behavior of the model in real life. Likewise, the choosing of the training window size $k$ can be evaluated with running this procedure multiple times with different sizes of $k$. Of course, to make the comparison fair, some exclusion of data might be needed as a larger window size needs more observations, and the number of observations for the error $\varepsilon_{t+j}$ will differ depending on the window size. Lastly, this procedure can be easily modified to estimate the variance of the error or other metrics of interest.

The errors are themselves objects of interest. As discussed previously $V_{t+j}\,\eqd\, \hat D_{t+j}(A,\beta)-D_{t+j}=\varepsilon_{t+j} (A,\beta)$. Then, the mean $\mu_{\varepsilon(t+j)}$ and the variance of the errors $\sigma^2_{\varepsilon(t+j)}$ can be used to approximate the error of predictions. For example, approximating $V_{t+j}=\mu_{\varepsilon(t+j)}+\sigma_{\varepsilon(t+j)}Z$ allows to construct confidence intervals for the demand. For example, the probability that the forecasted demand at time $t+j$ exceeds the quantity $x$ is equal to:
\begin{eqnarray}
\pr( D_{t+j}>x)=\pr(\hat D_{t+j}(A,\beta)+V_{t+j}>x)=\pr\left(Z>\frac{x-\hat D_{t+j}(A,\beta)-\mu_{\varepsilon(t+j)}}{\sigma_{\varepsilon(t+j)}}\right).
\end{eqnarray}
Which is the well-known tail cdf of a standard normal distribution $\bar\Phi\left(\frac{x-\hat D_{t+j}(A,\beta)-\mu_{\varepsilon(t+j)}}{\sigma_{\varepsilon(t+j)}}\right)$. It is important to note that the prediction horizon $j$ is crucial, because in general the longer the prediction horizon, the larger the errors. That being said, there is no free lunch in prediction since the estimation of the mean error and its variance is in itself a random variable. Moreover, there is an implicit assumption that the distribution of the stochastic process is in some sense stationary or fixed within a time interval. In the real-world there is no guarantee that this is the case.

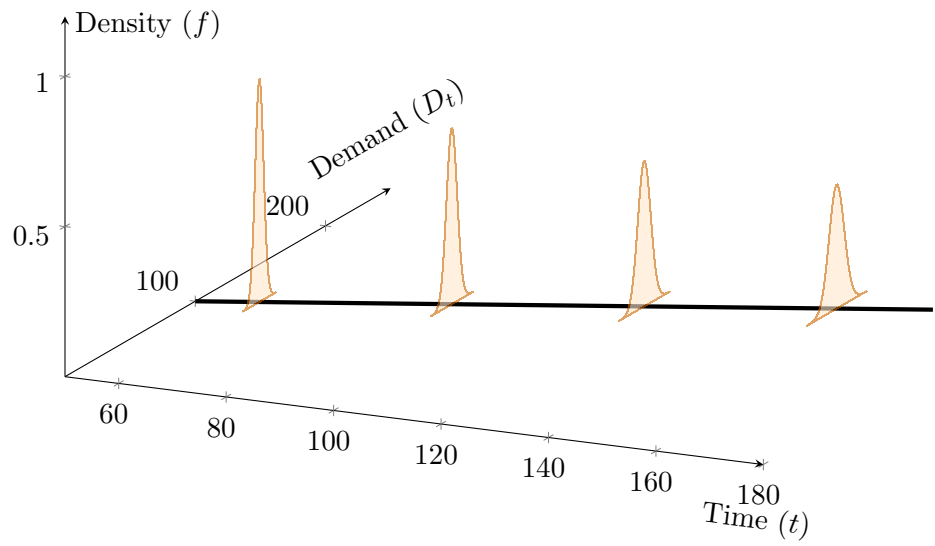
\begin{figure}[ht]
    \centering
    \begin{tikzpicture}
    \begin{axis}[
        width=0.95\textwidth,
        height=10cm,
        view={25}{30}, 
        xlabel={Time ($t$)},
        ylabel={Demand ($D_t$)},
        zlabel={Density ($f$)},
        zmin=0, zmax=1.2,
        xmin=50, xmax=180,
        ymin=0, ymax=250,
        ytick={0, 100, 200},
        axis lines=center,
        grid=both,
        grid style={line width=.1pt, draw=gray!15},
        clip=false,
        xlabel style={sloped like x axis, anchor=north, yshift=-12pt},
        xticklabel style={anchor=north, yshift=-3pt},
        ylabel style={sloped like y axis, anchor=south, xshift=10pt, yshift=10pt},
        yticklabel style={anchor=south east, xshift=-0.1pt, yshift=0.1pt},
    ]

    \addplot3 [
        black, 
        ultra thick, 
        samples y=0,
        domain=50:165,
    ] (x, {0.8*x + 60}, 0);

    \pgfplotsinvokeforeach{60, 90, 120, 150}{
        \addplot3 [
            fill=orange!20, 
            draw=orange!80!black, 
            opacity=0.6,
            samples=40,
            domain=({0.8*#1 + 60 - 4*(1.5 + #1/35)}):({0.8*#1 + 60 + 4*(1.5 + #1/35)})
        ] (
            {#1}, 
            {x}, 
            { (1/((1.5 + #1/35)*sqrt(2*3.1415))) * exp(-0.5*((x-(0.8*#1 + 60))/(1.5 + #1/35))^2) * 6 }
        ) -- (axis cs:#1, {0.8*#1 + 60}, 0) -- cycle;
    }

    \end{axis}
    \end{tikzpicture}
    \caption{The prediction trajectory $\hat{D}_{t+j}$ is plotted against Time ($t$), with the vertical Z-axis representing the probability density function. The increasing spread of the normal distributions illustrates the growing uncertainty $\sigma_{\varepsilon(t+j)}$ over the prediction horizon.}
    \label{fig:3d_probabilistic_final_clean}
\end{figure}


\chapter{Inventory Management Models}
\label{sec:inventory_management}
\section{Introduction}
Even with a good grasp of the magnitude and direction of the randomness affecting the future demand, there is always a decision behind predicting the future demand. Allocating resources to match random demand is a difficult problem. The difficulty lies in that resources are often fixed and not flexible: Determining the production capacity of a building, ordering a fixed amount of materials to manufacture products, or scheduling workers ahead of time to service a facility or store. What all these decision have in common is that there is an inherent trade-off between allocating too many resources or allocating too little. On one hand, allocating too many resources is costly and inefficient, and on the other, allocating too little, leads to congestion and lost sales. Both are undesirable outcomes that have to be balanced in an optimal manner. In this chapter, we present strategies to balance optimally these trade-offs in a data-driven way.

\section{The Newsvendor Model}
\label{sec:newsvendor}

Following our exposition of characterizing the demand as a Stochastic Process $\{D_t\}_{t=0}^\infty$, one of the most popular production models is known as the Newsvendor Model \cite{arrow1951optimal}. In this model, a production or inventory decision quantity $q$ has to be taken at time $t=0$ to satisfy a random demand $D_T$ at horizon time $T$. The goal of the model is to determine what is the optimal quantity of product that should be produced/ordered. Product is sold at a price $p$ with a unit cost $c$. The objective function of the seller is their expected profit, defined as:
\begin{eqnarray}
\max_q \ex [(D_T\wedge q)p-qc],
\end{eqnarray}
The revenue of the firm is the minimum operator (denoted by $a\wedge b=\min(a,b)$) between the production/inventory quantity $q$ and the random demand $D_T$ multiplied by the price $p$. As the production cost of each unit is $c$, the total variable cost is $qc$. Subtracting the cost from the revenue gives the profit (note that fixed costs do not affect the production quantity as they do not depend on $q$).%

\begin{figure}[htbp]
    \centering
    \begin{tikzpicture}
        \begin{axis}[
            no markers, 
            domain=-3:3, 
            samples=100,
            ymin=0, 
            axis lines*=left, 
            xlabel={Demand ($D_T$)}, 
            ylabel={Density},
            height=6cm, 
            width=10cm,
            xtick={0, 1}, 
            xticklabels={$\mu$, $q^*$},
            ytick=\empty,
            enlargelimits=upper,
            clip=false
        ]
            
            \addplot [name path=curve, very thick, blue!70!black] {exp(-x^2/2)};
            \path [name path=axis] (axis cs:-3,0) -- (axis cs:3,0);
            
            \addplot [cyan!20] fill between [
                of=curve and axis, 
                soft clip={domain=-3:1}
            ];
            
            \draw [dashed, thick, red!80!black] (axis cs:1,0) -- (axis cs:1,0.606);
            
            \draw [dotted] (axis cs:0,0) -- (axis cs:0,1);
            
            \node[anchor=west, align=left, font=\small] at (axis cs:1.2, 0.5) {
                $q^* = \mu + \sigma \Phi^{-1}\left(\frac{p-c}{p}\right)$
            };

        \end{axis}
    \end{tikzpicture}
    \caption{Optimal inventory level $q^*$ in the Newsvendor Model. The shaded area represents the probability of meeting all demand, corresponding to the critical ratio $(p-c)/p$.}
    \label{fig:newsvendor_normal}
\end{figure}
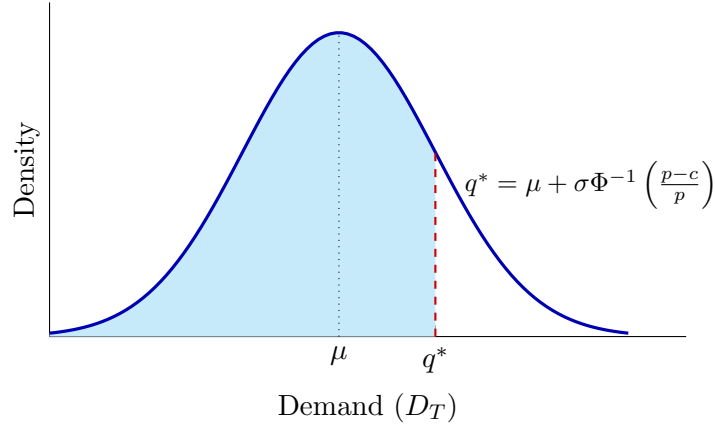

The first step to solve this problem is to note for any numbers $a,b\in\mathbb{R}$ we have $a\wedge b=a-(a-b)^+$ where $(x)^+=\max(x,0)$. With this fact, the profit can be re-written as $(D_T\wedge q)p-qc=[q-(q-D_T)^+]p-qc=q(p-c)-p(q-D_T)^+$. Then, the profit can be re-written as:
\begin{eqnarray}
	\label{eq:newsvendor_obj}
	\max_q q(p-c)-p\ex(q-D_T)^+,
\end{eqnarray}
To solve the problem we need to calculate what is the critical point $q^*$ such that the first derivative of the profit is equal to 0. We also need to check that the function is concave, as this will ensure that the critical point of the unique optimal solution maximizing the profit. Taking derivate and equating to 0, yields $\frac{\partial \ex (q-D_T)^+}{\partial q}=\frac{p-c}{p}$. We can use our definition of expected value (in Section \ref{sec:demand_rv}) to write $\frac{\partial \ex (q-D_T)^+}{\partial q}$ as follows:
\begin{eqnarray}
	\label{eq:trick_nv}
	\frac{\partial}{\partial q}\sum_{i=0}^{\infty} (q-i)^+p_i=\frac{\partial}{\partial q}\sum_{i=0}^{\infty} (q-i)p_i\ind\{i\le q\}=\frac{\partial}{\partial q}\sum_{i=0}^{q} (q-i)p_i=\pr(D_t\le q)=F_{D_T}(q).
\end{eqnarray}
The first equality is just a re-writting of the event $(x)^+=x\ind\{x\ge 0\}$. The second one just changes the limit of the sumation where the inner term is non-negative. The third equality exploits the fact that the derivative is a linear operator\footnote{The derivative of a sum is the sum of the derivatives of its terms.} and that the derivative of $qp_i$ is $p_i$ with respect to $q$. The second derivative is equal to 0, thus making the profit function concave and ensuring that the critical point in indeed the quantity that maximizes the profit. The optimal quantity is:
\begin{eqnarray}
q^*=F^{-1}_{D_T}\left(\frac{p-c}{p}\right),
\end{eqnarray}
$F^{-1}_{D_T}$ denotes the inverse cdf of the distribution of the demand. Connecting this to the demand distribution Section \ref{sec:demand}, for example by taking the data-driven estimate of the demand as the approximation $D_{T}=\hat D_T(\theta)+V_T\,\eqd\,\mu+\sigma Z$ with $\mu=\hat D_T(\theta)$ (in general, an arbitrary functional form for the demand. With $\theta=A,\beta$ we can take this as an instance the \emph{base model}) and $\sigma=\sqrt{\var(V_T)}$. With this approximation, the optimal quantity becomes:
\begin{eqnarray}
q^*=\mu+\Phi^{-1}\left(\frac{p-c}{p}\right)\sigma.
\end{eqnarray}
Where $\Phi^{-1}$ is the inverse cdf of the normal distribution. The interpretation of this expression is that the optimal quantity $q^*$ is the predicted value from the model $\mu=\hat D_T(\theta)$, plus a buffer proportional to the standard deviation of the demand $\sigma=\sqrt{\var(V_T)}$. The multiplier $\Phi^{-1}\left(\frac{p-c}{p}\right)$ balances how aggressive the provisioning of inventory should be. First, the price $p$ is normally higher than the cost $c$, that is, $p>c$. To get a feel for the multiplier, suppose the price is 3 times the cost, or $p=3c$ (a typical pricing strategy in the absence of information), then $(p-c)/p=2/3$. Evaluating the inverse at this value gives $\Phi^{-1}(2/3)\approx 43\%$, meaning that the extra buffer above the predicted demand $\mu$ is $43\%$ of its standard deviation $\sigma$. The larger the difference between the price and cost is $p-c$, the greater the amount of provisioning should be (because the loss of profit from a sale lost is orders of magnitude higher than the cost). As we will see next, we can modify this model to more complex situations where there are costs associated with unsold product. See Figure \ref{fig:newsvendor_normal} for a visual representation.

\begin{exm}[{\bf Avocado Toast Restaurant}]
	A national restaurant chain serves a signature avocado toast dish in their menu. The restaurant procures avocados from a wholesale producer and can purchase avocados at a price $c=\$0.68$ per unit. Every avocado dish is sold at a price $p=\$12$, other variable costs add up to $c_a=\$3$ per dish (other ingredients and labor). Every unused avocado has to be disposed at a cost of $c_d=\$0.20$. The distribution of the demand at time $T$ is estimated as a normal random variable with mean $\mu=70$ and standard deviation $\sigma=10$ dishes for a specific location of the restaurant. Also, when customers arrive and the avocado toast is sold-out, with probability $p_r=1\%$ they never come back (imagine a binary random variable $B_r$ with $\pr(B_r=1)=p_r=1-\pr(B_r=0)$, independent of the demand). The average lifetime profit per customer is estimated to be $\ell=\$200$. What is the optimal quantity of avocados that maximizes the expected profit for this location? What happens when the lifetime profit per customer increases to $\ell=\$2,000$?
	
	The expected profit for the restaurant can be written as:
	\begin{eqnarray}
		\max_q \ex[(D_T\wedge q)p-q(c+c_a)-(q-D_T)^+c_d-(D_T-q)^+B\ell],
	\end{eqnarray}
	Rearranging terms (noting $a\wedge b=a-(a-b)^+$ as in Equation (\ref{eq:newsvendor_obj})) leads to this expression:
	\begin{eqnarray}
		\max_q q(p-c-c_a)-(c_d+p)\ex(q-D_T)^+-p_r\ell\ex(D_T-q)^+,
	\end{eqnarray}
	Similar to the previous example, the strategy to solve the problem is to take derivatives with respect to $q$ and equate to 0. The derivative $\frac{\partial \ex(D_T-q)^+}{\partial q}$ is equal to $-\bar F_{D_T}(q)=F_{D_T}(q)-1$ (the trick is the same as in Equation (\ref{eq:trick_nv}), noting that the expression is non zero when the demand is between $q$ and $\infty$). Grouping terms leads to the following optimal demand quantity:
	\begin{eqnarray}
		q^*=F^{-1}_{D_T}\left(\frac{p-c-c_a-p_r\ell}{p+c_d-p_r\ell}\right).
	\end{eqnarray}
	In the case of this example, this leads to $q^*=\mu+\sigma\Phi^{-1}\left(\frac{p-c-c_a-p_r\ell}{p+c_d-p_r\ell}\right)\approx 76.03$ avocados. When the lifetime profit per customer $\ell$ is $\$2,000$ the optimal quantity rises to $q^*\approx 81.72$ avocados.
\end{exm}
\subsection{Empirical Distributions}
The Newsvendor model relies heavily on the accuracy of the inverse cdf $F^{-1}_{D_T}$. In real life, the true distribution of the demand is rarely known, and the estimation has to be done on the empirical distribution from a sample. The empirical distribution of the demand is constructed by first sorting observations of the demand in ascending order. A sample of i.i.d. observations of the demand $\{d_i\}_{i=1}^N$ (here, the big assumption is that these are all independently sampled observations of the demand $D_T$). Let $d_{(i)}$ be the $i$-th smallest element of the sample, such that $d_{(1)}\le d_{(2)}\le \cdots \le d_{(N)}$. The \emph{Empirical probability measure} assigns each sample a equal probability $1/n$. The sample or empirical mean is then defined as $\frac{1}{n}\sum_{i=1}^nd_i$ which by the LLN converges to the true mean as we discussed in the previous chapter. Likewise, the empirical cdf is defined as $\hat F_n(q)= i/n$ for $i=\text{argmax} \{i:d_{(i)}\le q\}$.  Similarly, the empirical inverse cdf is defined as $\hat F_n^{-1}(p)=d_{(i)}$ for $i=\text{argmin} \{i:i\ge pn\}$. A natural question is how far off is the Newsvendor inventory quantity when using the empirical cdf in lieu of the real one.

A useful bound for the divergence of the empirical cdf and the actual one was developed in two papers \cite{dvoretzky1956asymptotic,massart1990tight} and is commonly known as the \emph{DKW bound} (due to the authors of the first paper, which was later refined by the second one). The bound is stated as follows:
\begin{eqnarray}
	\label{eq:dkw}
\pr\left(\sup_q |F(q)-\hat F_n(q)|> \lambda \right)\le2 e^{-2n\lambda^2},
\end{eqnarray}
What this bound states is that the probability that the deviation from the empirical and the actual cdf exceeds a value $\lambda$ (for any value $q$) is less than $2e^{-2n\lambda^2}$. As $n$ increases, this probability decays to $0$. Meaning that as the sample size increases, the accuracy of the empirical cdf converges to the true one. Interestingly, we can use this inequality to build confidence intervals on the empirical cdf. For example, with probability $1-\alpha$ the true value of the cdf $F(q)$ is within $\hat F_n(q)\pm\sqrt{\frac{1}{2n}\ln\left(\frac{2}{\alpha}\right)}$ (hint: equate the right hand side of (\ref{eq:dkw}) to $\alpha$) for any $q$.

Going back to our topic of interest, we want to know how far is the solution of the empirical Newsvendor model $\hat q_n^*=F_n^{-1}(\gamma)$ with respect to the real solution $q^*=F^{-1}(\gamma)$ (here, $\gamma$ denotes the optimal constant in the Newsvendor model, for example $\gamma=\frac{p-c}{c}$ in the plain vanilla formulation). In Section \ref{sec:demand} we emphasized that event calculus is a powerful problem solving tool. That is, to write complex events in terms of events that we know well probabilistically. In this case, we know about bounding the empirical cdf in Equation (\ref{eq:dkw}) and need to find the right event to express the relation to the empirical inverse cdf $F_n^{-1}(\gamma)$.

\begin{figure}[htbp]
    \centering
    \begin{tikzpicture}
        \begin{axis}[
            width=10cm, height=7cm,
            axis lines=left,
            xmin=0, xmax=1.1,
            ymin=0, ymax=10,
            xlabel={Probability ($p$)},
            ylabel={Quantity ($q = \hat{F}_n^{-1}(p)$)},
            xtick={0.4, 0.7, 1.0},
            xticklabels={$\gamma-\epsilon$, $\gamma$, $\gamma+\epsilon$},
            ytick={3.4136, 5.4623, 8.7355}, 
            yticklabels={$q_L$, $\hat{q}_n^*$, $q_U$},
            xticklabel style={font=\small},
            yticklabel style={font=\small},
            clip=false,
            declare function={invcdf(\x) = 1.8 * exp(1.6 * \x);}
        ]

            \addplot [very thick, blue!70!black, domain=0.15:1.02, samples=100] {invcdf(x)};
            \node [blue!70!black, anchor=west, font=\small] at (axis cs:0.75, 9.5) {$\hat{F}_n^{-1}(p)$};

            \draw [red, thick] (axis cs:0.4, 0) -- (axis cs:0.4, 3.4136) -- (axis cs:0, 3.4136);
            \filldraw [red] (axis cs:0.4, 3.4136) circle (1.5pt);
            
            \draw [red, thick] (axis cs:1.0, 0) -- (axis cs:1.0, 8.7355) -- (axis cs:0, 8.7355);
            \filldraw [red] (axis cs:1.0, 8.7355) circle (1.5pt);
            
            \draw [dashed, gray] (axis cs:0.7, 0) -- (axis cs:0.7, 5.4623) -- (axis cs:0, 5.4623);
            \filldraw [black] (axis cs:0.7, 5.4623) circle (2pt) node[anchor=south east, font=\small] {$(\gamma, \hat{q}_n^*)$};

            \draw [thick] (axis cs:0.4, -0.2) -- (axis cs:1.0, -0.2);
            \draw [thick] (axis cs:0.4, -0.1) -- (axis cs:0.4, -0.3);
            \draw [thick] (axis cs:0.7, -0.1) -- (axis cs:0.7, -0.3);
            \draw [thick] (axis cs:1.0, -0.1) -- (axis cs:1.0, -0.3);

        \end{axis}
    \end{tikzpicture}
    \caption{Mapping the DKW confidence interval. The symmetric interval $\gamma \pm \epsilon$ results in a non-symmetric confidence interval $[q_L, q_U]$ due to the nonlinearity of the inverse CDF.}
    \label{fig:dkw_non_symmetric}
\end{figure}
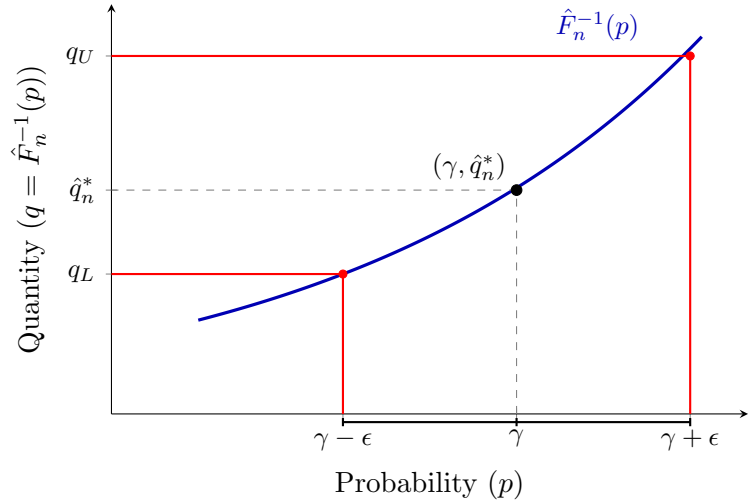

The optimal quantity $q^*=F^{-1}(\gamma)=\inf \{q:F(q)\ge \gamma\}$. Now we know that $F(q)$ is within $\hat F_n(q)\pm\sqrt{\frac{1}{2n}\ln\left(\frac{2}{\alpha}\right)}$ with probability at least $1-\alpha$. Putting these two facts together we get that $q^*\in \inf \{q:\hat F_n(q)\pm\sqrt{\frac{1}{2n}\ln\left(\frac{2}{\alpha}\right)}\ge \gamma\}$ or $q^*\in F_n^{-1}\left(\gamma\pm \sqrt{\frac{1}{2n}\ln\left(\frac{2}{\alpha}\right)}\right)$ with probability $1-\alpha$. Note that this bound is not symmetric and not centered at $F^{-1}_n(\gamma)=q_n^*$. See Figure \ref{fig:dkw_non_symmetric}.

Often, it is also a natural question what should be the average profit using a quantity $\hat q$ (for example, the optimal sample inventory $q_n^*)$ given a series of $n$ historical observations of the demand. The sample average profit $\mu_\pi(n)$ of the plain vanilla Newsvendor model is:
\begin{eqnarray}
\mu_\pi(n) = \frac{1}{n}\sum_{i=1}^n \hat q(p-c) - p (\hat q-d_i)^+,
\end{eqnarray}
We know that because of the CLT, this mean will converge to the true average profit $\mu_\pi=\hat q(p-c) - p \ex(\hat q-D_T)^+$. The variance of the profit is $\var (p[(\hat q-D_T)^+])=p^2\var(\hat q -D_T)^+$, for any distribution of the demand $D_T$, this variance is at most $q^2/4$\footnote{This is the so-called Popovicius's variance inequality. Let $X$ be a random variable such that $\pr(X\in[a,b])=1$, then, $\var(X)\le (b-a)^2/4$.}. Then, we have the following probabilistic relationship between the sample and true profit:
\begin{eqnarray}
	\mu_\pi(n)\,\eqd \,\mu_\pi+\frac{p\sqrt{\var(\hat q -D_T)^+}}{\sqrt{n}}Z \le \mu_\pi+\frac{p\hat q}{2\sqrt{n}}Z,
\end{eqnarray}
The first equality is the equality in distribution of the CLT. The second is the result of bounding the variance with the inequality above. The result says that for any distribution of the demand, in the i.i.d. case solving the plain vanilla Newsvendor problem converges to the true demand as $n$ grows, and that the noise is bounded by the product of the price $p$ and the quantity $\hat q$ divided by a factor of $2\sqrt{n}$. For example, to obtain a confidence interval with probability $1-\alpha$ of the true average profit $\mu_\pi$ by using quantity $\hat q$, we have that the error $|\mu_\pi-\mu_\pi(n)|\,\eqd\, \frac{p\hat q}{2\sqrt{n}}|Z|$. Then, we want to find the value $\epsilon$ such that $\pr\left (|Z|p\hat q/(2\sqrt{n})>\epsilon\right)=\alpha$, which is the same as $\pr (Z>2\epsilon  \sqrt{n}/(p\hat q))=\alpha/2$. Applying the inverse tail cdf of the normal distribution, we get solving for $\epsilon$:
\begin{eqnarray}
\mu_\pi\in \mu_\pi(n)\pm  \bar\Phi^{-1}(\alpha/2)\frac{p\hat q}{2\sqrt{n}}\text{    w.p.   }1-\alpha,
\end{eqnarray}
For example, choosing $\alpha=5\%$, $\bar\Phi^{-1}(\alpha/2)=1.96$, meaning that the true profit is with 95\% probability within $1.96\frac{p\hat q}{2\sqrt{n}}$ of the sample profit in the plain vanilla Newsvendor model. This technique can also be used for all variations of the Newsvendor problem in the i.i.d. case, noting that the complement variance $\var (D_T-\hat q)^+$ is also bounded.

These bounds are very useful because the do not depend on the distribution of the demand and allow to infer the true performance of the model as a function of the sample size. We illustrate this in the following example.

\begin{exm}[{\bf Launching a new product}]
	A supermarket chain is planning to release a new bottled orange juice. They perform a pilot of 30 days ($n=30$), observing an average daily profit of \$1,500. The price of the juice is \$3.49 and they produce per store $\hat q =100$. Assuming the distribution of the demand of juices is i.i.d. in the plain vanilla Newsvendor model. The shelf opportunity cost of the supermarket (the minimum profit that a product makes in that shelf area) is \$1,450 per day. With this data, the supermarket wants to conduct a hypothesis test to decide whether to include the orange juice permanently or continue the pilot. What's the probability that the daily profit is less than \$1,450?
	
With our analysis we know the true profit $\mu_\pi$ is bounded in distribution to $\mu_\pi(n)+\frac{p\hat q}{2\sqrt{n}}Z=\$1,500+(\$349/2\sqrt{30})Z$. Then $\pr(\$1,500+(\$349/2\sqrt{30})Z<\$1,450)$ simplifies to $\pr(Z<-\$50(2\sqrt{30})/\$349)$, which is $\pr(Z<-1.56)=\Phi(-1.56)=5.9\%$. Which means that for a 95\% confidence one sided test, the pilot should continue (as the probability is above 5\%).
\end{exm}

\section{Inventory Replenishment - The $(s,S)$ rule}

In the previous section, we proposed a model for optimizing a production quantity such that profit is maximized in a static setting. The methodology can be extended to optimize multiple products and locations. As well as optimizing production accross time. In this section, we deal with a common inventory policy taking into account the evolution of the demand over time.

In the previous chapter, we characterized the randomness of the demand as a Stochastic Process $\{D_t\}_{t=0}^\infty$. Assuming the distribution of the demand is known, define the inventory level at period $t$ by $I_t$. The inventory can be defined as the quantity of a commodity or product that is stored to service a random demand. The simplest rule to replenish the inventory is known as the $(s,S)$ rule. This rule works as follows: at the end of period $t$, if the inventory $I_t$ is below a level $s$. An order (assumed to fulfill instantaneously) replenishes the inventory to a level $S$. Then, the inventory $I_{t+1}$ is equal to $S-D_{t+1}$ when $I_t\le s$, and it is $I_t-D_{t+1}$ otherwise (in the case when $I_t>s$). Clearly, the time between inventory replenishments is random and depends on the distribution of the demand.

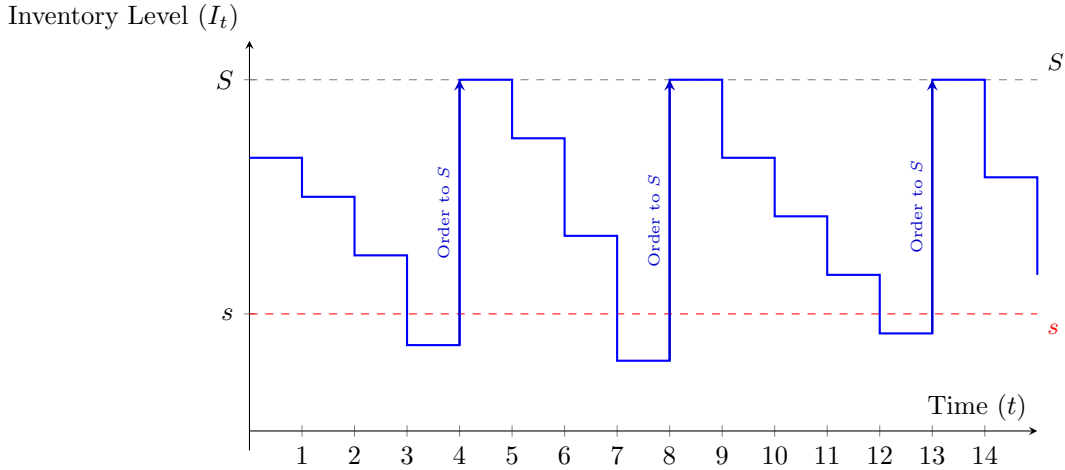
\begin{figure}[htbp]
    \centering
    \begin{tikzpicture}
        \begin{axis}[
            width=12cm, height=7cm,
            axis lines=middle,
            xlabel={Time ($t$)},
            ylabel={Inventory Level ($I_t$)},
            ylabel style={at={(axis description cs:0,1)}, anchor=south east},
            xmin=0, xmax=15,
            ymin=-0.5, ymax=10,
            xtick={0,1,2,3,4,5,6,7,8,9,10,11,12,13,14},
            ytick={3, 9},
            yticklabels={$s$, $S$},
            tick label style={font=\small},
            label style={font=\small},
            grid=none,
            clip=false
        ]
            \draw [dashed, gray] (axis cs:0,9) -- (axis cs:15,9); 
            \draw [dashed, red] (axis cs:0,3) -- (axis cs:15,3);  
            
            \addplot[thick, blue, const plot] coordinates {
                (0,7) (1,6) (2,4.5) (3,2.2) 
                (4,9) (5,7.5) (6,5) (7,1.8)  
                (8,9) (9,7) (10,5.5) (11,4) (12,2.5) 
                (13,9) (14,6.5) (15,4)
            };

            \draw[stealth-, blue!80!black, thick] (axis cs:4,9) -- (axis cs:4,2.2);
            \node[rotate=90, anchor=south, font=\tiny, blue!80!black] at (axis cs:4, 5.6) {Order to $S$};

            \draw[stealth-, blue!80!black, thick] (axis cs:8,9) -- (axis cs:8,1.8);
            \node[rotate=90, anchor=south, font=\tiny, blue!80!black] at (axis cs:8, 5.4) {Order to $S$};

            \draw[stealth-, blue!80!black, thick] (axis cs:13,9) -- (axis cs:13,2.5);
            \node[rotate=90, anchor=south, font=\tiny, blue!80!black] at (axis cs:13, 5.75) {Order to $S$};

            \node[anchor=south west, font=\footnotesize] at (axis cs:15,9) {$S$};
            \node[anchor=north west, font=\footnotesize, red] at (axis cs:15,3) {$s$};

        \end{axis}
    \end{tikzpicture}
    \caption{Inventory level $I_t$ over time under the $(s,S)$ policy. Replenishment is triggered the moment $I_t < s$, instantly bringing the stock back to the maximum level $S$.}
    \label{fig:sS_inventory_path}
\end{figure}

Under the $(s,S)$ rule, the inventory is in itself a random variable that depends on the demand. The conditional distribution of the inventory is given by:
\begin{eqnarray}
	\pr(I_{t+1}=j|I_t=i)=\pr(D_{t+1}=S\ind\{i<s\}+i\ind\{s\le i\le S\}-j),
\end{eqnarray}
When the distribution of the demand is i.i.d., that is, $D_t\,\eqd D$ for all $t$. We can express the inventory process $I_t$ as a Markov Chain\footnote{A Markov Chain is a special class of Stochastic Processes, where the distribution of the process $X_t$ depends only on the realization of the process the previous period $X_{t-1}$.} with transition probabilities $p_{ij}=\pr(D=S\ind\{i<s\}+i\ind\{s\le i\le S\}-j)$. Note that the inventory can take negative values as the demand can exceed the inventory quantity $I_t$. Conversely, under the $(s,S)$ rule, the inventory never exceeds $S$ by construction. See Figure \ref{fig:sS_inventory_path}.

With these elements we have enough building blocks to express problems of interest: For example, suppose that everytime there is a replenishment of inventory, there is a fixed cost $K$ plus a variable cost of $c$ per unit. Likewise, suppose there is a holding cost of $c_h$ per day the inventory is stored. Finally, when the inventory is negative (meaning that there was not enough inventory to service the demand), those sales are lost, a loss of a net revenue $\tilde p$ per unit. We can summarize these costs as a function $c(i)$ for each possible state of the inventory:

\begin{eqnarray}
\label{eq:inv_cost}
c(i)= \begin{cases}
    ic_h, & \text{if } s\le i \le S, \\
    ic_h+K+(S-i)c,   & \text{if }  0\le i < s, \\
    K+Sc-i\tilde p, & \text{if } i < 0.
\end{cases}
\end{eqnarray}

A Markov Chain transitions between states forever (even when it gets absorbed into a single state). When the transitions can be expressed with a fixed matrix $P$ with elements $p_{ij}$ denoting the probability that the chain at state $i$ transitions to state $j$. The so-called limiting distribution is given by a vector $\pi$ as:
\begin{eqnarray}
\lim_{t\rightarrow \infty} P^t=\pi\,\,\,\text{ or }\,\, \pi P=\pi,
\end{eqnarray}
Verbally, the vector $\pi$ denotes the distribution of the states long-term. For example, $\pi(i)=\lim_{t\rightarrow\infty}\sum_{s=1}^t\ex \ind\{X_s=i\}/t$ shows what proportion of time (out of $t$) is spent in state $i$ for a Markov Chain $\{X_t\}$.

\begin{figure}[htbp]
    \centering
    \begin{tikzpicture}
        \begin{axis}[
            width=13cm, height=8cm,
            ybar=0pt, 
            bar width=6pt,
            ylabel={Magnitude},
            xlabel={Inventory State ($i$)},
            xtick={-2, -1, 0, 1, 2, 3, 4, 5, 6},
            xmin=-2.5, xmax=6.5,
            ymin=0, ymax=1.2,
            legend style={at={(0.5,1.1)}, anchor=south, legend columns=3, font=\small},
            ymajorgrids=true,
            axis lines*=left,
            clip=false
        ]
            \addplot[fill=blue!40, draw=blue!60] coordinates {
                (-2, 0.04) (-1, 0.06) (0, 0.08) (1, 0.10) 
                (2, 0.18) (3, 0.16) (4, 0.14) (5, 0.12) (6, 0.12)
            };
            \addlegendentry{$\pi(i)$}

            \addplot[fill=red!40, draw=red!60] coordinates {
                (-2, 1.00) (-1, 0.85) 
                (0, 0.60) (1, 0.50)  
                (2, 0.10) (3, 0.15) (4, 0.20) (5, 0.25) (6, 0.30) 
            };
            \addlegendentry{$c(i)$}

            \addplot[fill=purple!60, draw=purple!80] coordinates {
                (-2, 0.04) (-1, 0.051) (0, 0.048) (1, 0.05) 
                (2, 0.018) (3, 0.024) (4, 0.028) (5, 0.03) (6, 0.036)
            };
            \addlegendentry{$\pi(i)c(i)$}

            \draw[dashed, gray, thin] (axis cs:1.5, 0) -- (axis cs:1.5, 1.1);
            \draw[dashed, gray, thin] (axis cs:-0.5, 0) -- (axis cs:-0.5, 1.1);
            
            \node[font=\tiny, gray, align=center] at (axis cs:-1.5, 1.15) {Lost Sales\\$i < 0$};
            \node[font=\tiny, gray, align=center] at (axis cs:0.5, 1.15) {Replenish\\$0 \le i < s$};
            \node[font=\tiny, gray, align=center] at (axis cs:4, 1.15) {Holding\\$i \ge s$};

        \end{axis}
    \end{tikzpicture}
    \caption{Long-term inventory cost analysis. The purple bars represent the weighted cost contribution of each state to the total expected cost $\sum \pi(i)c(i)$. Optimization aims to choose $s$ and $S$ such that the high-probability states ($\pi$) align with low-cost states ($c$).}
    \label{fig:pi_c_distribution}
\end{figure}
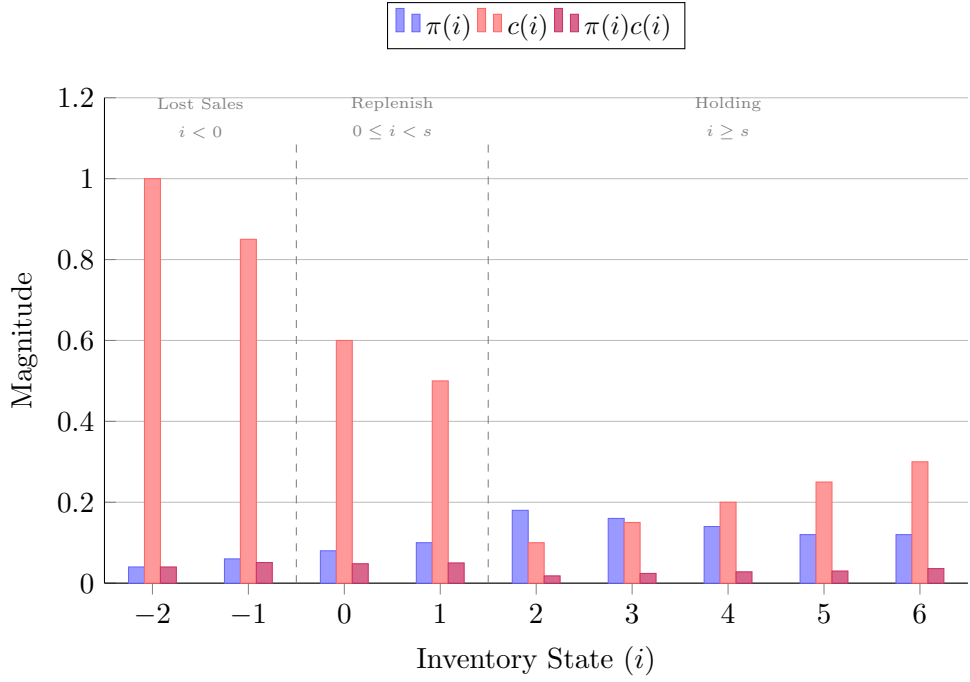

In our case, this limiting distribution $\pi$ is useful because we know we can express our inventory model as a Markov Chain with transition probabilities $p_{ij}=\pr(D=S\ind\{i<s\}+i\ind\{s\le i\le S\}-j)$. Then, the long-term cost of the $(s,S)$ inventory rule is given by:

\begin{eqnarray}
\sum_{i=s-M}^S \pi(i)c(i),
\end{eqnarray}
Where $M$ is given by a value such that $\pr(D\le M)=1$. That is, the largest the demand can be (why is the lower limit of this Markov Chain $s-M$?). To compute the limiting distribution we just need to solve the system of equations $\pi P=\pi$ numerically. We spend the rest of this section discussing how to find actual values for $(s,S)$. See Figure \ref{fig:pi_c_distribution} for an example analysis.

\subsection{Finding the $(s,S)$ values in the i.i.d. case}
The most direct approach to find a pair of values $0\le s<S$ that is optimal in the sense of minimizing is to enumerate all possible combinations of values of $s,S$ and calculate the stationary distribution by solving $\pi P(s,S)=\pi$, where $P(s,S)$ is the transition matrix under the policy $(s,S)$ (recall, the transition probabilities depend on both $s$ and $S$, as $p_{ij}(s,S)=\pr(D=S\ind\{i<s\}+i\ind\{s\le i\le S\}-j)$). We summarize the procedure in the following algorithm:

\begin{algorithm}
	\caption{Finding Optimal $(s^*,S^*)$ Inventory Policy - i.i.d. case}
	\label{alg:findmax_inv}
	\begin{algorithmic}[1]
		\State \textbf{Input:} A distribution of the demand $\pr(D=i)$. A cost function $c_\theta(i)$ such as the one in Equation (\ref{eq:inv_cost}) with parameters $\theta$. And upper limit of inventory capacity $S_{max}$, and for the demand $M$.
		\State \textbf{Output:} The optimal inventory policy parameters $(s^*,S^*)$.
		\Procedure{OptimalInventory}{$F_D,c_\theta,S_{max}$}
		\State $S,s^*,S^*=S_{max}$
		\State $c^*,\bar c=\infty$
		\For{$S = 1$ \textbf{to} $S_{max}$}
		\For{$s = 0$ \textbf{to} $S-1$}
		\State $P\gets [p_{ij}]=[\pr(D=S\ind\{i<s\}+i\ind\{s\le i\le S\}-j)]\text{ for }i,j\in\{s-M,\cdots,S\}$.
		\State Solve $\pi$ such that $\pi=\pi P$ and $|\pi|=1$.
		\State Update and compute $\bar c\gets c_\theta'\pi$.
		\If{$\bar c<c^*$}
		\State  Update $c^*\gets \bar c,s^*\gets s,S^*\gets S$.
		\EndIf
		\EndFor
		\EndFor
		\State \textbf{return} $s^*,S^*,c^*$.
		\EndProcedure
	\end{algorithmic}
\end{algorithm}

Given the structure of the transitions of the inventory, this is not the most efficient algorithm to solve the inventory policy, but it is the easiest to implement as a computational routine. Next we examine the case to find policies for time-varying demands.

\subsection{The time-varying demand optimal inventory}
The key observation in the previous subsection was that the inventory process is a Markov Chain, that is, that the transition probabilities are time-homogeneous. As we explored in Section \ref{sec:demand}, demand is often seasonal and time-dependent. In this section, we present a useful technique to solve problems that have a time-dependent structure.

In the time-dependent setting, the demand distribution $F_{D_t}(i)=\pr(D_t\le i)$ changes from period to period. Suppose we want to derive an inventory policy, that at each period tells us the optimal ordering level that minimizes some cost (for example, the cost function in Equation (\ref{eq:inv_cost})). Given a time horizon of $T$ periods in the future. The objective is to minimize a cost over the time horizon as:
\begin{eqnarray}
	\ex\left [\sum_{t=1}^Tc_t(I_t)\right],
\end{eqnarray}

The trajectory of the inventory $I_t$ depends on the previous periods and is not time homogeneous, making the problem seemingly hard to solve. A useful technique to solve sequential problems is to look one step ahead (or behind) and try to define a recursive relation that allows to solve the problem going to the last period $T$ (where the decision is normally trivial) and solve the problem \emph{backwards}, in a recursive way. The easiest way to set up such relation is to define the \emph{action} that can be taken separated from the \emph{state} of the system, so that at each step the problem reduces to find the optimal action at that state in that time step. The state of the system in this case is already defined as the inventory level $I_t$, an integer (potentially negative) denoting the inventory amount at the end of period $t$. In the previous section, for the $(s,S)$ policy, the action was coded in an indicator function (as $I_{t+1}=S\ind\{I_t<s\}+I_t(1-\ind\{I_t<s\})-D_{t+1}$ or $I_{t+1}=I_t+(S-I_{t})\ind\{I_t<s\}-D_{t+1}$). Inspired by this, we can express the state and the action of the system $I_{t+1}=j-D_{t+1}$ such that $j\ge I_t=i$. In this case, $j$ represents the new level of replenishment of the inventory (or rather, the inventory at the beginning of $t+1$, before the demand $D_{t+1}$ is subtracted). The strategy is then to write a function $\bar C_t(i)$ that represents the optimal expected cost from time $t$ to $T$, that is, $\bar C_{t}(i)=\ex\left[\sum_{s=t}^Tc^*_s(I_s)|I_t=i\right]$, where $c^*$ is the optimal cost of taking the best inventory actions. The strategy is then, to write down a one-step recursion that allows us to break down the problem to a single replenishment decision at every time $t$. With our representation, we can express the optimal cost between states and actions by the following relation:
\begin{eqnarray}
	\label{eq:bellman_inv}
	\bar C_t(i)=\min_{j\ge i}\ex \left[c(i,j)+\bar C_{t+1}(j-D_{t+1})\right],
\end{eqnarray}
The interpretation of this equation is that the optimal cost $\bar C_t(i)$ at time $t$ when the last inventory $I_t=i$ is the expected cost of that state $i$ and the optimal action $j$, for that period in the function $c(i,j)$. Plus, the expected optimal cost of the new inventory at a level $j-D_{t+1}$. The minimum denotes choosing the action $j\ge i$ that minimizes this cost.

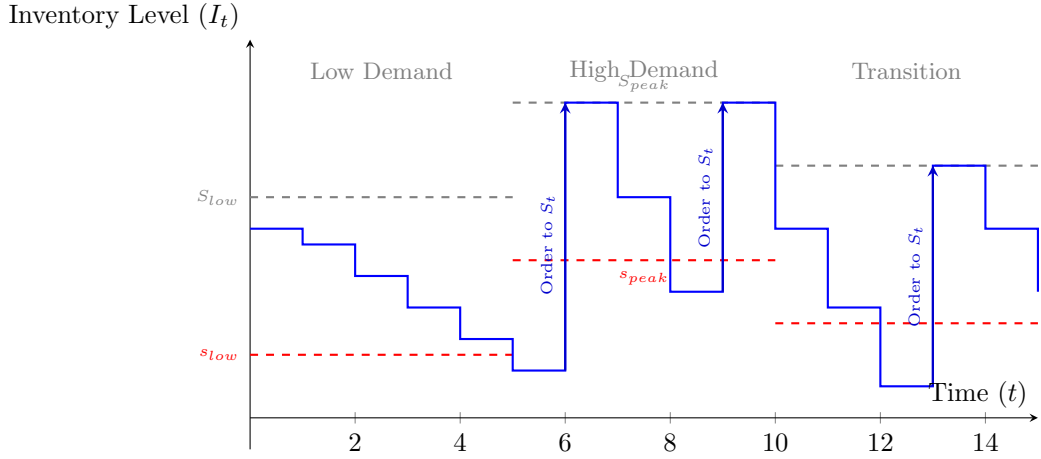
\begin{figure}[htbp]
    \centering
    \begin{tikzpicture}
        \begin{axis}[
            width=12cm, height=7cm,
            axis lines=middle,
            xlabel={Time ($t$)},
            ylabel={Inventory Level ($I_t$)},
            ylabel style={at={(axis description cs:0,1)}, anchor=south east},
            xmin=0, xmax=15,
            ymin=-1, ymax=12,
            xtick={0,2,4,6,8,10,12,14},
            ytick=\empty,
            tick label style={font=\small},
            label style={font=\small},
            grid=none,
            clip=false
        ]
            \draw[dashed, gray, thick] (axis cs:0,7) -- (axis cs:5,7) node[pos=0, left, font=\tiny] {$S_{low}$};
            \draw[dashed, red, thick]  (axis cs:0,2) -- (axis cs:5,2) node[pos=0, left, font=\tiny] {$s_{low}$};
            
            \draw[dashed, gray, thick] (axis cs:5,10) -- (axis cs:10,10) node[midway, above, font=\tiny] {$S_{peak}$};
            \draw[dashed, red, thick]  (axis cs:5,5)  -- (axis cs:10,5)  node[midway, below, font=\tiny] {$s_{peak}$};
            
            \draw[dashed, gray, thick] (axis cs:10,8) -- (axis cs:15,8);
            \draw[dashed, red, thick]  (axis cs:10,3) -- (axis cs:15,3);

            \addplot[thick, blue, const plot] coordinates {
                (0,6) (1,5.5) (2,4.5) (3,3.5) (4,2.5) (5,1.5) 
                (6,10) (7,7) (8,4) 
                (9,10) (10,6) (11,3.5) (12,1) 
                (13,8) (14,6) (15,4)
            };

            \draw[stealth-, blue!80!black, thick] (axis cs:6,10) -- (axis cs:6,1.5);
            \node[rotate=90, anchor=south, font=\tiny, blue!80!black] at (axis cs:6, 5.5) {Order to $S_t$};

            \draw[stealth-, blue!80!black, thick] (axis cs:9,10) -- (axis cs:9,4);
            \node[rotate=90, anchor=south, font=\tiny, blue!80!black] at (axis cs:9, 7) {Order to $S_t$};

            \draw[stealth-, blue!80!black, thick] (axis cs:13,8) -- (axis cs:13,1);
            \node[rotate=90, anchor=south, font=\tiny, blue!80!black] at (axis cs:13, 4.5) {Order to $S_t$};

            \node[font=\footnotesize, gray] at (axis cs:2.5, 11) {Low Demand};
            \node[font=\footnotesize, gray] at (axis cs:7.5, 11) {High Demand};
            \node[font=\footnotesize, gray] at (axis cs:12.5, 11) {Transition};

        \end{axis}
    \end{tikzpicture}
    \caption{Inventory process with time-varying demand. As the depletion rate (slope) increases during peak periods, the Dynamic Programming solution adjusts the optimal $s_t$ and $S_t$ levels to minimize expected costs over the horizon $T$.}
    \label{fig:time_varying_sS}
\end{figure}

With this recursion, it is possible to solve for all $C_{t}(i)$ for all $i=-M,\dots,S_{max}$ and times $t=0,1,\dots,T$, starting from the end at time $T$. And populating entries backwards (it is easier to imagine a matrix $\bar C$ with entries $c_{i,t}=\bar C_t(i)$). This problem solving technique is known as \emph{Dynamic Programming} and Equation (\ref{eq:bellman_inv}) is a \emph{Bellman Equation}\footnote{A equation representing a recursive relation of the optimal decision of a dynamic problem.}. Storing the optimal replenishment levels $j^*_{t,i}=\text{argmin}_{j\ge i}\ex \left[c(i,j)+\bar C_{t+1}(j-D_{t+1})\right]$ is also a good idea, as they will be a \emph{dictionary} to the optimal replenishment level when the inventory is $i$ at time $t$. See Figure \ref{fig:time_varying_sS} to visualize the dynamic inventory policy.

The modeling aspect is captured in the cost function $c(i,j)$. The goal is to modify this function to correctly capture the underlying business problem that we want to solve. We illustrate this with the following example:

\begin{exm}[{\bf Oxygen Concentrator Distribution Center}]

An oxygen concentrator company that delivers oxygen concentrators (a device used for patients that need a constant flow of oxygen while they are at home) rents a warehouse to be a distribution center for the concentrators in a major city. They have a random seasonal daily demand for concentrators $D_t$, that among other features, depends on the atmospheric humidity level $X_t$ (there is a correlation between respiratory diseases and the level of humidity in the air). The cost structure of running the warehouse works as follows: orders that are fulfilled with stock from the warehouse have a unit cost of $c$ dollars (using an in-house delivery system). Orders that are fulfilled when the demand exceeds the stock cost $\tilde c$ dollars using an external courier that picks up a concentrator in another facility and delivers it to the client, this cost is higher than the in-house system, that is $\tilde c>c$. Stock replenishment orders have a fixed cost $K$ per truck plus a unit cost $c_r$ (the lead time is overnight, and each truck can store up to $L$ concentrators). Holding the concentrators in the facility has an associated daily cost $c_h$ (including maintenance and the internal logistics of the warehouse that are proportional to the number of concentrators). How to model and minimize this cost over the next month ($T=30$, assuming today is $t=0$) assuming a known high-fidelity forecast of the humidity levels $X_t$ for next month $\{x_t\}_{t=0}^T$?

The first step is to write the cost function $c(i,j)$ to reflect these costs. For example, $c(i,j)=i^+ c_h+\tilde c(D_{t+1}-j)^++\ind\{j>i\}(\lceil (j^+-i^+)/L\rceil K+(j^+-i^+)c_r)$, note the $i^+=\max(i,0)$ reflects that when the inventory ends up negative, these orders are not \emph{back-ordered} as they still end up fulfilled at the higher cost $\tilde c$. $(j^+-i^+)$ denotes the order quantity. The quantity $\lceil (j^+-i^+)/L\rceil$ denotes the number of trucks needed ($\lceil a \rceil$ is the closest integer above $a$, typically referred in programing as the ceiling function\footnote{A definition is $\lceil a \rceil=\min\{x\in \mathbb{Z}:x\ge a\}$.}).

The second step is to incorporate the humidity forecast $\{x_t\}_{t=0}^T$. To do this, recall in Section \ref{sec:demand} that a way to incorporate side-information is to analyze the conditional random variable, that is, $D_t|X_t$. Calibrating this conditional distribution results in conditional probabilities $\pr(D_{t}=i|\{X_{t}=x_{t}\}_{t=0}^T)$. Solving the problem with this side information amounts to calculate the expectations with these conditional probabilities (instead of the unconditional ones). Replacing both the new cost $c(i,j)$ and taking expectations conditional to $X_t=x_t$ in (\ref{eq:bellman_inv}). Explicitly, this yields the recursion for $\bar C_t(i)=\min_{j\ge i}\ex [i^+ c_h+\tilde c(D_{t+1}-j)^++\ind\{j>i\}(\lceil (j^+-i^+)/L\rceil K+(j^+-i^+)c_r)+\bar C_{t+1}(j-D_{t+1})|\{X_{t}=x_{t}\}_{t=0}^T]$. 

\end{exm}


\chapter{Network Fulfilment Models}
\label{sec:networks}
\section{Introduction}

So far, we have described mathematical models to match the random demand in the time dimension. In this Chapter, we change the pace by thinking of problems satisfying the demand physically in the space dimension. Most if not all problems in supply chain analysis deal with some form of transportation problem as the demand requires moving goods or supplies between the producers and customers.

Matching the demand in space also adds a non-trivial time delay that needs to be accounted for when matching resources to satisfy the demand.  Large companies like Amazon have to optimize the location of their warehouses and orchestrate a large chain of distribution networks to deliver their products in record time. Likewise, ridesharing marketplaces (such as Uber or Doordash) have to balance a demand that changes in space and time with a supply side of workers and drivers that will service different geographical locations as they finish their deliveries and rides.

We will explore a series of common problems and solutions to supply chain management problems in space, as well as algorithms and heuristics to solve them.

\section{Warehouse consolidation and positioning}

Imagine a supplier wanting to build a warehouse to service demand in a city or geographical location. The supplier would want to build the warehouse in a place that is economically viable and resource efficient to service this demand.

While the physical world that we inhabit is 3-dimensional, it is often useful to think of abstractions that simplify the problem: For example, most geographical maps use a 2-dimensional representation of the Earth (this makes sense as the Earth is more or less a manifold locally, i.e. locally, we live in a flat world). Perhaps, one of the most powerful abstractions for physical spaces are graphs. A graph $G$ is defined as a set of vertices $V$ and edges $E$.

\begin{figure}[htb]
    \centering
    
    \begin{subfigure}[b]{0.48\textwidth}
        \centering
        
        \includegraphics[width=\textwidth]{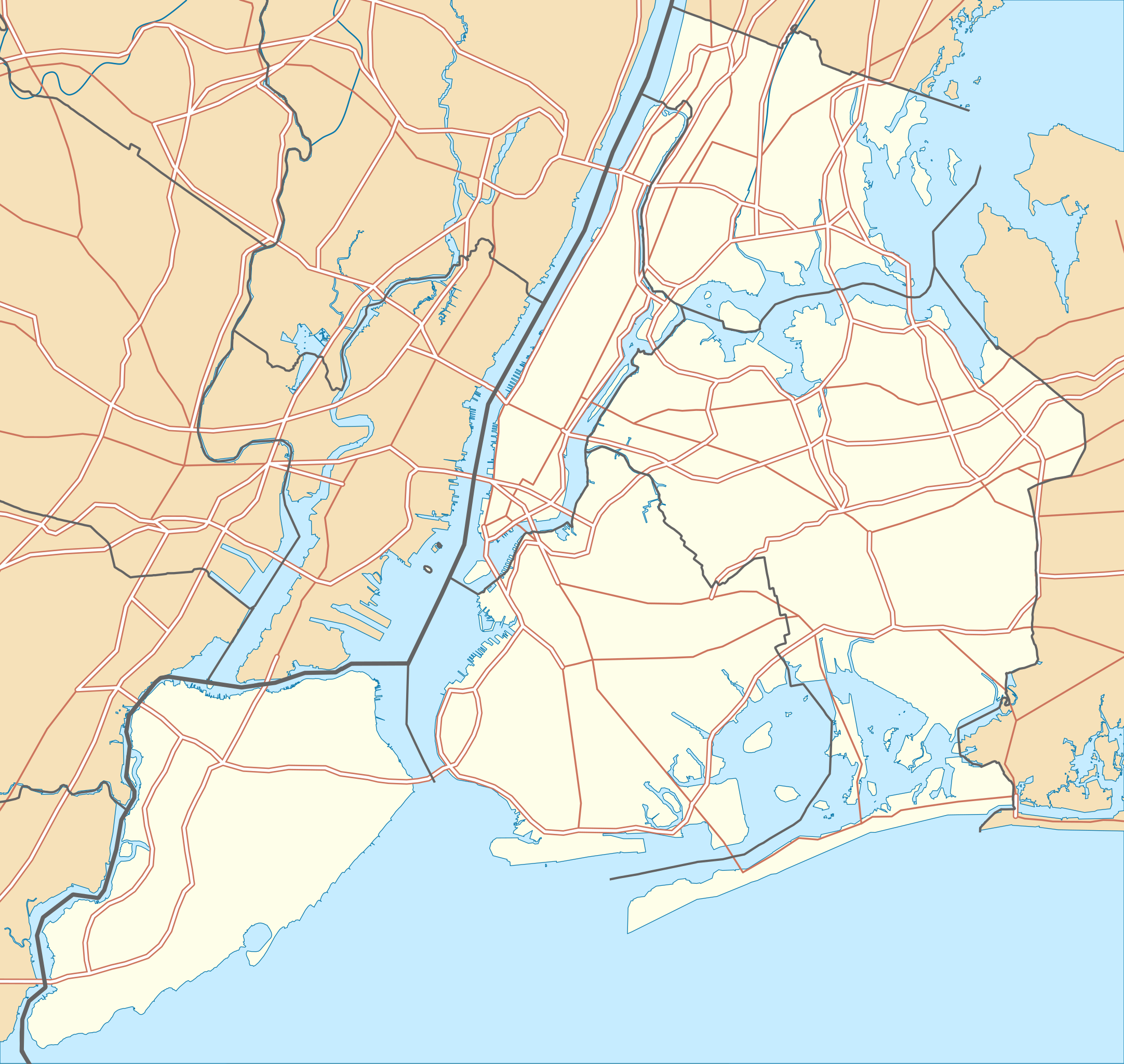}
        
        \caption{Geographic Map of NYC's Five Boroughs.}
        \label{subfig:nyc_map_img}
    \end{subfigure}%
    \hfill
    \begin{subfigure}[b]{0.48\textwidth}
        \centering
        \begin{tikzpicture}[
            scale=1.5,
            node distance=2.5cm, 
            main/.style={rectangle, draw=black, fill=gray!30, minimum width=2.5em, minimum height=1.5em, font=\bfseries\small},
            edge_style/.style={very thick, black} 
        ]
            \node[main] (BX) at (0, 3) {The Bronx};
            \node[main] (MN) [below left of=BX] {Manhattan};
            \node[main] (QN) [below right of=BX] {Queens};
            \node[main] (BK) [below right of=MN] {Brooklyn};
            \node[main] (SI) [below of=BK] {Staten Island};

            
            \draw[edge_style] (BX) -- (MN);
            
            \draw[edge_style] (BX) -- (QN);
            
            \draw[edge_style] (MN) -- (QN);
            
            \draw[edge_style] (QN) -- (BK);
            
            \draw[edge_style] (BK) -- (MN);
            
            \draw[edge_style] (SI) -- (BK);
            
        \end{tikzpicture}
        \caption{Graph abstraction of NYC.}
        \label{subfig:nyc_graph_full_accessibility}
    \end{subfigure}    

    \caption{Graph abstraction of NYC representing each of the five Boroughs as vertices, and their roads connections as edges.}
    \label{fig:nyc_map_graph_final}
\end{figure}
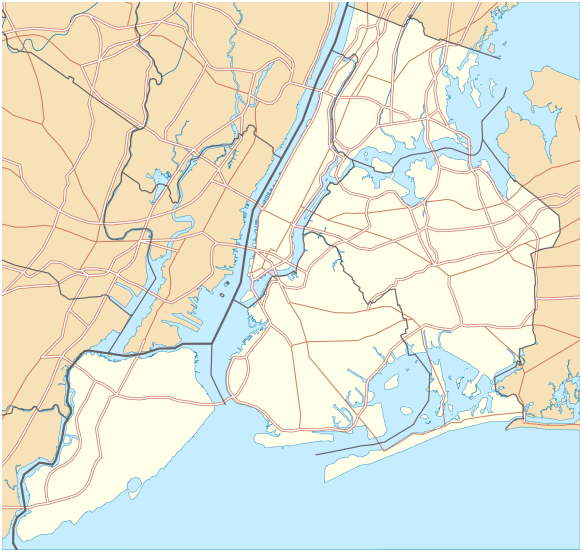
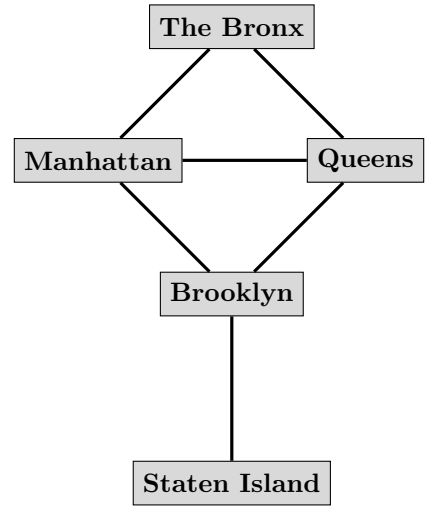

Graphs can be seen as discretization of a continuous space into a disjoint partition, connecting with edges all physically adjacent partitions. For example, take a map such as in the left side of Figure \ref{fig:nyc_map_graph_final} with the five Boroughs of New York City (Manhattan, The Bronx, Queens, Brooklyn and Staten Island) and make each Borough a vertex, drawing an edge if there is a road connection between two Boroughs (vertices).

Modeling a physical space as a graph has the natural trade-off of having tractable algorithms to solve a wide-range of problems while losing some of the accuracy and resolution of the underlying modeled real-world system. For example, in the graph of the figure, direct travels between Manhattan and Queens can occur via two different bridges, but there is only one edge representing their connection. In modeling traveling times, there is either the choice of creating more nodes and egdes to represent the different routes and to attach numerically the respective travel time to each edge, or to \emph{average} into the single edge all the possible routes between Queens and Manhattan. Clearly, the tradeoff is between accuracy of the abstraction and tractability of the solution, as most optimization algorithms scale proportional to the size of the Graph. In general, mathematical modeling is the art of choosing an appropriate abstraction of a real-life system that preserves the behavior of interest while being algorithmically tractable.

Going back to the problem of choosing a location for a warehouse, start with a graph $G$, composed by $N$ vertices $V=\{1,2,\dots,N\}$ and edges $E\subseteq\{(i,j)|(i,j)\in V^2\}$. For example, the graph in Figure \ref{fig:nyc_map_graph_final} is composed by $N=5$ vertices $V=\{1,2,3,4,5\}$, each vertex corresponds to a Borough: $1\rightarrow\text{Manhattan}$, $2\rightarrow\text{The Bronx}$, $3\rightarrow\text{Queens}$, $4\rightarrow\text{Brooklyn}$ and $5\rightarrowtail\text{Staten Island}$. And edges $E=\{(1,2),(1,3),(1,4),(2,3),(3,4),(4,5)\}$. 

Vertices can have different unique attributes, for example, in the model of this section each vertex has an associated random demand $D_i$ for $i=1,\dots,N$. This demand, is the usual demand random variable and denotes the random realization of the demand (in this case, there is not yet a temporal element to the demand and it can be taken as a static model with i.i.d. realizations of the demand that repeat every day). Likewise, every edge $(i,j)$ for $i,j\in\{1,\dots,N\}$ also has a random variable $T_{ij}$ denoting the random traveling time of a vehicle starting in vertex $i$ to vertex $j$.

\subsection{Single Location Warehouse}
With these elements, we can now start to model the problem of choosing a location to build a warehouse. Let $c_i$ be the cost of running a warehouse in location $i$ with enough capacity to service the demand $\sum_{i=1}^ND_i$ with high probability (see Chapters \ref{sec:demand} and \ref{sec:inventory_management} on how this number can be calculated, here it is taken as given\footnote{In practice, one way would be to calculate a capacity $\epsilon$ such that $\pr (\sum_{i=1}^ND_i>\epsilon)\le 1-\alpha$ and then quote or calculate the actual cost of running a warehouse of that capacity for the value of $c_i$. You can also set up a Newsvendor-like model to calculate a capacity that is also sensitive to missed sales or holding costs.}). In this scenario, if it's not possible to run a warehouse in some vertex, then let the cost be $c_i=\infty$ for such vertex (or exclude it from the set of possible locations to run the warehouse). Likewise, let the cost of transports between vertices $(i,j)$ be equal to\footnote{$a\propto b$ denotes that $b$ is proportional to $a$ up to a constant, say $b=2a$.} $c_{ij}\propto  T_{ij}$, for example the cost of the labor of a courier plus the cost of gas that is proportional to the random traveling time $T_{ij}$, that is $c_{ij}=\gamma_{ij}T_{ij}$, clearly $c_{ij}$ is also a random variable. 

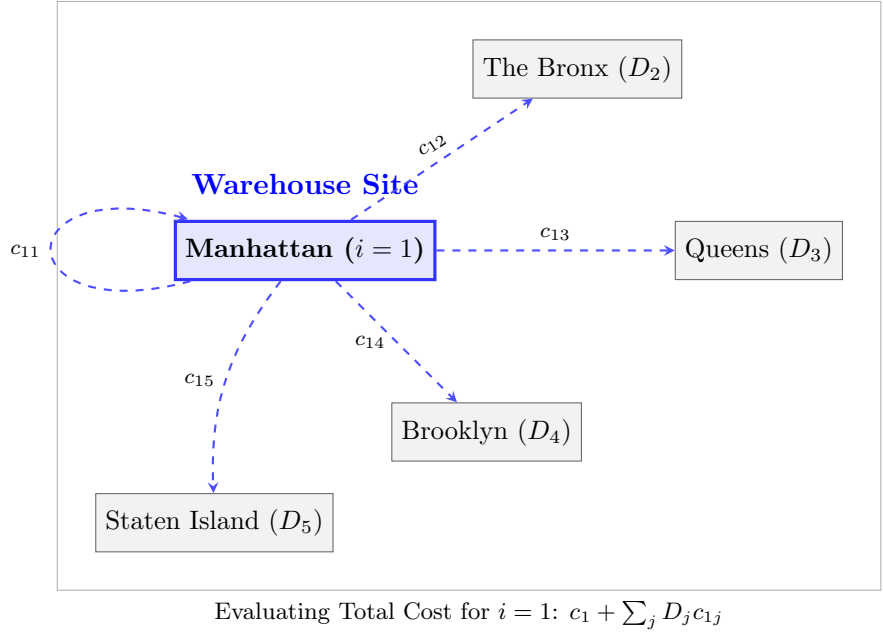
\begin{figure}[htb]
    \centering
    \begin{tikzpicture}[
        scale=1.2,
        node distance=2cm,
        warehouse/.style={rectangle, draw=blue!80, fill=blue!10, very thick, minimum width=3em, minimum height=2em, font=\bfseries\small},
        demand/.style={rectangle, draw=black!60, fill=gray!10, minimum width=3em, minimum height=2em, font=\small},
        flow/.style={->, >=stealth, thick, blue!70, dashed}
    ]

        \node[demand] (BX) at (0, 3) {The Bronx ($D_2$)};
        \node[warehouse] (MN) at (-3, 1) {Manhattan ($i=1$)};
        \node[demand] (QN) at (2, 1) {Queens ($D_3$)};
        \node[demand] (BK) at (-1, -1) {Brooklyn ($D_4$)};
        \node[demand] (SI) at (-4, -2) {Staten Island ($D_5$)};

        \node[blue, above=0.2cm of MN] {\textbf{Warehouse Site}};

        \draw[flow] (MN) -- node[above, sloped, black, font=\scriptsize] {$c_{12}$} (BX);
        \draw[flow] (MN) -- node[above, black, font=\scriptsize] {$c_{13}$} (QN);
        \draw[flow] (MN) -- node[left, black, font=\scriptsize] {$c_{14}$} (BK);
        \draw[flow] (MN) to [bend right=20] node[left, black, font=\scriptsize] {$c_{15}$} (SI);
        
        \path[flow] (MN) edge [loop left] node[black, font=\scriptsize] {$c_{11}$} (MN);

        \begin{pgfonlayer}{background}
            \node[draw=black!30, fill=white, fit=(BX) (MN) (QN) (BK) (SI), inner sep=0.5cm, label=below:{\footnotesize Evaluating Total Cost for $i=1$: $c_1 + \sum_{j} D_j c_{1j}$}] {};
        \end{pgfonlayer}

    \end{tikzpicture}
    \caption{Visualization of the enumeration step for $i=1$ (Manhattan). To solve the problem, the supplier calculates the expected total cost (warehouse overhead $c_i$ plus weighted transportation costs to all demand nodes $D_j c_{ij}$) for each potential vertex and selects the minimum.}
    \label{fig:warehouse_enumeration}
\end{figure}

For a warehouse built on location $i$, the cost is equal to $c_i+\sum_{j=1}^ND_jc_{ij}$. As this number is a random variable, a first approximation is to look at the expected value and choose over all possible locations, that is, solving the enumeration problem:

\begin{eqnarray}
\min_{i=1,\dots,N} \ex\left[c_i+\sum_{j=1}^ND_jc_{ij}\right],
\end{eqnarray}
It is an enumeration problem because it amounts to calculate the cost for each location and choose the minimum. The emphasis in writing the costs as random variables instead of the more common practice of writing them as constants is that there is a conscious choice to solve the problem only in expected value or to consider other moments of their distribution. For example, it could be the case that running the warehouse in one location has the lowest cost in expected value, but the variation (either of the demands and the trips) is higher than another location that has a slightly higher cost in expected value, but less variance (in this case, variance denotes the daily variation of running the warehouse satisfying random demands every day or period). See Figure \ref{fig:warehouse_enumeration} to get a feel for evaluating the cost of building a single warehouse in Manhattan.

In general, for problems involving a random amount of money (either a cost or a payoff) there is an amount of risk aversion that the decision-maker is willing to take, and it needs to be made explicit in the problem formulation, or otherwise be acknowledged by the modeler. For example, given two lotteries: one that pays $X_1=\$100,000B_1+(-\$100,000)(1-B_1)$ and another $X_2=\$100B_2+(-\$100)(1-B_2)$ with two independent coins $B_1,B_2$ with $\pr(B=1)=\pr(B=0)=1/2$. That is, two loteries that with 50\% chance pay a positive amount and with the same chance the player loses the same amount, in the first case \$100,000 and in the other \$100. While both lotteries pay in average \$0, that is $\ex(X_1)=\ex(X_2)=\$0$ they are clearly not the same game. Depending on their risk aversion, some players would choose to play lottery 1, while others would play lottery 2 (or even choose not to play at all). The quantitative difference between the two games is their variance. We have that the variance of the first lottery is $\var(X_1)=2^2(\$100,000)^2\var(B_1)=(\$100,000)^2$ and the second one $\var(X_2)=(\$100)^2$. The risk ratio between the two lotteries is $\var(X_1)/\var(X_2)=10^6$. In other words, playing the first lottery has a million times higher variance than the second. While exaggerated, this example shows the potential literal risk of only modeling decisions based on their expected value. A quantitative approach to evaluate each lottery could be to \emph{weight} the variance or the standard deviation, for example $\ex(X_1)+\lambda \sqrt{\var(X_1)}$ weights the standard deviation (in dollars). A weight $\lambda=0$ denotes risk indifference (thus, only caring for the expected return of the lottery). $\lambda<0$ denotes risk aversion, while $\lambda>0$ implies risk appetite. Note that the units of $\ex(X_1)+\lambda \sqrt{\var(X_1)}$ is dollars. 

Another justification coming from the study of risk measures also justifies this by defining the risk measure of a random variable $X$ as $R_\lambda(X)=\frac{1}{\lambda}\ln{\ex(e^{\lambda X})}$. This risk measure is used in risk management as it captures risk aversion, is convex and coherent and it is sensitive to tail risks. In fact, we show that our strategy of weighting the risk is approximately equivalent to this risk measure. To see why, note that $e^u=1+u+u^2/2+\dots$. Then, $\ex(e^{\lambda X})=\ex(1+\lambda X+\frac{\lambda^2X^2}{2}+\dots)=1+\lambda\ex(X)+\frac{\lambda^2}{2}\ex(X^2)+\dots=1+v$ where $v=\lambda\ex(X)+\frac{\lambda^2}{2}\ex(X^2)+\dots$, also note that $\ln(1+v)=v-\frac{v^2}{2}+\dots$. Using these relations and noting that terms with $v^3$ and higher vanish, we get:
\begin{eqnarray}
 \ln{\ex(e^{\lambda X})}&=&\lambda\ex(X)+\frac{\lambda^2}{2}\ex(X^2)-\frac{1}{2}\left[\lambda\ex(X)+\frac{\lambda^2}{2}\ex(X^2)\right]^2+\dots\\
              &=&\lambda\ex(X)+\frac{\lambda^2}{2}\ex(X^2)-\frac{\lambda^2}{2}[\ex(X)]^2+\dots\nonumber\\
              &=&\lambda\ex(X)+\frac{\lambda^2}{2}[\ex(X^2)-\ex(X)^2]+\dots\nonumber\\
              &=&\lambda\ex(X)+\frac{\lambda^2}{2}\var(X)+O(\lambda^3)\nonumber
\end{eqnarray}

Noting that the term $O(\lambda^3)$ vanishes to 0 and dividing both sides by $\lambda$ yields the desired relation. 

This shows that optimizing the objective $\ex(X)+\lambda\var(X)$ is a principled way of handling the risk. An important aspect in practice is calibrating the risk aversion parameter $\lambda$. A well-known principle to do this, is known as the \emph{Certainty Equivalent}. This is simply an amount of money that the decision maker is willing to pay (or receive) to be indiferent between the lottery and a certain payoff. Solving for that amount would yield the implied risk aversion coefficient. For example, suppose there is a lottery $X=-\$100B+\$0(1-B)=-\$100B$, that is, losing a \$100 with probability $\pr(B=1)$ or losing 0 with probability $\pr(B=0)$. Suppose $P(B=1)=1/2$, we have $\ex(X)=-\$50$ and $\var(X)=(\$100)^2/4$. A person is willing to pay \$40 to exit the lottery, this implies a Certainty Equivalent (CE) given by $-40=-50+\frac{\lambda}{8}100^2$ or $\lambda=80/100^2$. In general, giving the stakeholders similar lotteries and with their CEs is a standard way to learn the risk aversion of a stakeholder, we will illustrate numerically how the trade-off between risk and expected cost plays out later.

Going back to the warehouse location problem, the problem incorporating variance can be stated as follows:
\begin{eqnarray}
\min_{i=1,\dots,N} \ex\left[c_i+\sum_{j=1}^ND_jc_{ij}\right]+\lambda\sqrt{\var\left(\sum_{j=1}^ND_jc_{ij}\right)},
\end{eqnarray}
With $\lambda>0$ this optimization problem models a risk-averse decision-maker that penalizes variance in their cost (caused by both variance in the demand and in the transportation times). This is again an enumeration problem that amounts to calculate the objective for each potential warehouse location. It is also possible to consider other risk measures that penalize variation.

In the case where the cost of the facilities $c_i$ is fixed. The variance of the cost $\var [c_k+\sum_{j=1}^ND_jc_{kj}]$ can be calculated as $\var [\sum_{j=1}^ND_jc_{kj}]=\bm{1}\Sigma_k\bm{1}^\intercal$ where the entry $i,j$ of the matrix $\Sigma_k$ is $\cov(D_ic_{ki},D_jc_{kj})=\gamma_{ki}\gamma_{kj}\cov(D_iT_{ki},D_jT_{kj})$. For example, in the case of the diagonal we have $\cov(D_ic_{ki},D_ic_{ki})=\gamma^2_{ki}\var(D_iT_{ki})$ if the demand and the travel times are independent then this is equal to $\gamma^2_{ki}\var(D_i)\var(T_{ki})$. For the other elements off-diagonal we use the definition of covariance\footnote{$\cov(X,Y)=\ex(XY)-\ex(X)\ex(Y)$.}, then, $\cov(D_iT_{ki},D_jT_{kj})=\ex(D_iT_{ki}D_jT_{kj})-\ex(D_iT_{ki})\ex(D_jT_{kj})$. With the assumption that while the demands and the travel times are independent, but not the demands and the travel times with each other. For example, it makes sense that travel jams in Manhattan could affect traffic in The Bronx, and also that the demand for a product is Brooklyn is correlated with the demand for the same product in Queens. With this assumption, we have that $\ex(D_iT_{ki}D_jT_{kj})$ is equal to $\ex(D_iD_j)\ex(T_{ki}T_{kj})$ and $\ex(D_iT_{ki})\ex(D_jT_{kj})$ is equal to $\ex(D_i)\ex(T_{ki})\ex(D_j)\ex(T_{kj})$. In summary, the covariance matrix $\Sigma_k$ has entries $i,j$ equal to:
\begin{eqnarray}
	\Sigma_k(i,j)=\begin{cases}
		\gamma^2_{ki}\var(D_i)\var(T_{ki}), & \text{for } i=j \\
		\gamma_{ki}\gamma_{kj}[\ex(D_iD_j)\ex(T_{ki}T_{kj})-\ex(D_i)\ex(D_j)\ex(T_{ki})\ex(T_{kj})], & \text{for } i \neq j.
	\end{cases}
\end{eqnarray}

\subsection{Multiple Location Warehouses}

Sometimes it can be more efficient to choose multiple locations to satisfy the geographical demand of multiple locations than a single warehouse location.  This can be because the cost structure is more favorable (closer deliveries are cheaper), due to economies of scale, or simply ease of logistical burden by breaking down the demand to multiple difference warehouses. The problem can also be interpreted as choosing providers to satisfy demand for different geographical locations.

\begin{figure}[htb]
    \centering
    \begin{tikzpicture}[
        scale=1.1,
        warehouse/.style={circle, draw=blue!80, fill=blue!5, thick, minimum size=25pt, font=\small},
        target/.style={circle, draw=gray!80, fill=gray!5, thick, minimum size=25pt, font=\small},
        active_edge/.style={->, >=Stealth, thick, blue!80},
        alt_edge/.style={->, >=Stealth, thin, gray!30, dashed}
    ]

    \node[warehouse] (W1) at (0, 3) {$W_1$};
    \node[warehouse] (W2) at (0, 1.5) {$W_2$};
    \node[warehouse, draw=red!60, fill=red!5] (W3) at (0, 0) {$W_3$}; 
    \node[warehouse] (W4) at (0, -1.5) {$W_M$};
    
    \node[left=0.3cm of W1, font=\footnotesize\itshape] {Candidates};
    \node[left=0.3cm of W3, font=\footnotesize, red] {Inactive ($x_i=0$)};

    \node[target] (T1) at (4, 3) {$D_1$};
    \node[target] (T2) at (4, 1.5) {$D_2$};
    \node[target] (T3) at (4, 0) {$D_3$};
    \node[target] (T4) at (4, -1.5) {$D_N$};
    
    \node[right=0.3cm of T1, font=\footnotesize\itshape] {Targets};

    \foreach \i in {W1, W2, W4}
        \foreach \j in {T1, T2, T3, T4}
            \draw[alt_edge] (\i) -- (\j);

    \draw[active_edge] (W1) -- node[above, pos=0.3, black, font=\scriptsize] {$q_{11}$} (T1);
    \draw[active_edge] (W1) -- node[above, pos=0.7, black, font=\scriptsize] {$q_{12}$} (T2);
    \draw[active_edge] (W2) -- node[below, pos=0.3, black, font=\scriptsize] {$q_{22}$} (T2);
    \draw[active_edge] (W2) -- node[above, pos=0.7, black, font=\scriptsize] {$q_{23}$} (T3);
    \draw[active_edge] (W4) -- node[below, pos=0.3, black, font=\scriptsize] {$q_{M3}$} (T3);
    \draw[active_edge] (W4) -- node[below, pos=0.3, black, font=\scriptsize] {$q_{MN}$} (T4);

    \draw [decorate, decoration={brace, amplitude=10pt, mirror}, xshift=-10pt] 
        (-0.5,-2) -- (-0.5,3.5) 
        node [black, midway, xshift=-1.2cm, align=center, font=\footnotesize] 
        {$2^M$ possible\\warehouse subsets};

    \end{tikzpicture}
    \caption{Bipartite representation of the Multiple Location Warehouse problem. Solid blue lines represent an active assignment $q_{ij}$ from selected warehouses ($x_i=1$). The light dashed lines represent the combinatorial complexity: the model must select a subset of $M$ locations and then optimize the flows to $N$ targets.}
    \label{fig:bipartite_mip}
\end{figure}
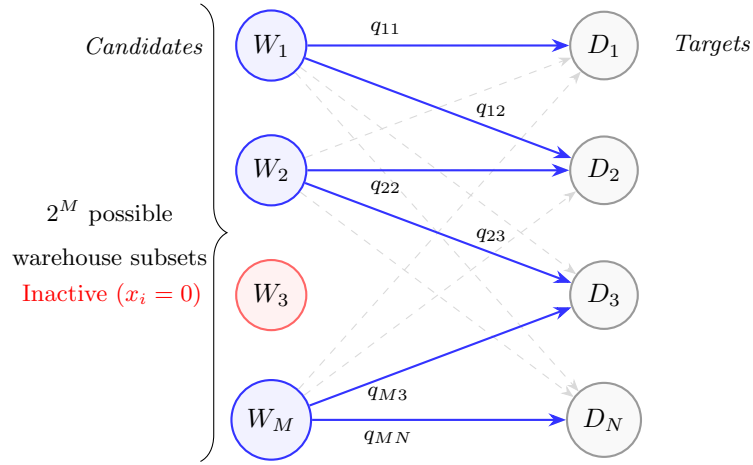

To modify the previous section formulation, now each candidate location has a finite capacity $C_i$ to service demand. Likewise, we modify the cost $c_i$ to be a fixed (and given) cost of running the warehouse $i$. The second and more important modification is that since multiple locations can be chosen to service the demand, brute enumeration is in general not a great strategy, except for small instances of the problem. To see why, suppose there are $M$ candidate locations, all possible choices of subsets of locations is $2^M$, which is a number that grows exponentially in the number of possible locations\footnote{Given a set of $M$ elements, each can be either included or not included in a subset. Since this is true for all elements, the total number of subsets is $2^M$, also called the \emph{power set}.}. Combinatorial Optimization studies algorithms to solve problems where the solution takes the form of subsets by exploiting the internal structure of the problem, not to be forced to enumerate all possible $2^M$ candidate solutions. A very powerful technique to solve this problems is called Mixed Integer Programming, that is very similar to Linear Programming (both in structure and algorithmically) by allowing integer variables (normally encoding either counting quantities or codyfing yes-no decisions as variables). See Figure \ref{fig:bipartite_mip} for a typical representation of the problem of matching possible $M$ warehouses locations with $N$ demand target locations.

A Mixed Integer Programming formulation of this problem entails adding binary variables $x_i\in\{0,1\}$ denoting whether or not the location is chosen. Also, the variables $q_{ij}\ge 0$ denote the amount of demand (quantity) satisfied from warehouse $i$ to location $j=1,\dots,N$. The next modeling step is to relate variables $x_i$ and $q_{ij}$ to each other in an equation (that normally becomes a constraint of the problem). For each warehouse it must be true that when it is selected ($x_i=1$), the total quantity that can be serviced from the location cannot exceed its capacity $C_i$. In other words, $\sum_{j=1}^Nq_{ij}\le C_i$ when $x_i=1$, and conversely when $x_i=0$ it should be $\sum_{j=1}^Nq_{ij}=0$. MIP formulations normally involve tricks manipulating real-world conditions into equations (see \cite{aimms2018aimms} for some of the most common tricks). In this case, the trick is to write the constraint as $\sum_{j=1}^Nq_{ij}\le x_iC_i$. It is immediate that this constraint achieves the goal of encoding both situations.

Likewise, it must be true that the random demand is satisfied to some level. For example, if looking at the problem in expected value, that the demand is satisfied in average is the constraint $\sum_{i=1}^Mq_{ij}\ge \ex (D_j)$. Another common technique to deal with this type of constraints is to write a constraint as $\pr(\sum_{i=1}^Mq_{ij}\le D_j)\le \eta$, meaning that with probability not exceeding $\eta$ (say a small percentage like 5\%) the random demand will be satistied. In the case where the demand can be approximated as a normal random variable, that is, $D_j\,\eqd\, \ex (D_j)+\sqrt{\var(D_j)}Z$. The constraint $\pr(\sum_{i=1}^Mq_{ij}\le D_j)\le \eta$ can then be written as $\pr(Z\le [\sum_{i=1}^Mq_{ij}-\ex (D_j)]/\sqrt{\var(D_j)})\ge 1-\eta$, or the constraint $\sum_{i=1}^Mq_{ij}\ge \ex (D_j)+\sqrt{\var(D_j)}\Phi^{-1}(1-\eta)$, the similarity to the Newsvendor model studied in Chapter \ref{sec:inventory_management} should be evident, as the contraint conveys that the quantity delivered should be the mean of the demand plus a buffer proportional to the standard deviation of the demand. The last constraint to be enforced is that $q_{ij}=0$ when $x_i=0$ and positive otherwise. Letting $\bar d_j$ be an upper bound for the demand such that $\pr(D_j>\bar d_j)=0$, this can be achieved with the constraint $q_{ij}\le x_i\bar d_j$
With these constraints, we can write a first fomulation of the problem in expected value as (let $W\subseteq V$ be the set of possible locations for the warehouses to make the notation cleaner):
\begin{eqnarray}
\label{eq:warehouse_ex_cost}
	&\underset{q_{ij}\ge 0,x_i\in\{0,1\}}{\min}& \ex\left\{ \sum_{i\in W}c_ix_i+\sum_{i\in W}\sum_{j\in V}c_{ij}q_{ij}\right\}\\
	&\text{s.t. }&\sum_{i\in W}q_{ij}\ge \ex (D_j) \,\,\,\text{for}\,\,\,j\in V, \nonumber\\
	&&\sum_{j\in V}q_{ij}\le x_iC_i \,\,\,\text{for}\,\,\,i\in W,\nonumber\\	
	&&q_{ij}\le x_i\bar d_j\,\,\,\text{for}\,\,\,i,j\in W\times V. \nonumber
\end{eqnarray}
As the costs $c_i$ are fixed, the objective can be written as $\sum_{i\in W}c_ix_i+\sum_{i\in W}\sum_{j\in V}\ex(c_{ij})q_{ij}$. This formulation is the typical warehouse optimization found in textbooks and can be solved with off-the-shelf MIP solvers. 

\begin{figure}[htb]
    \centering
    \includegraphics[width=0.8\textwidth]{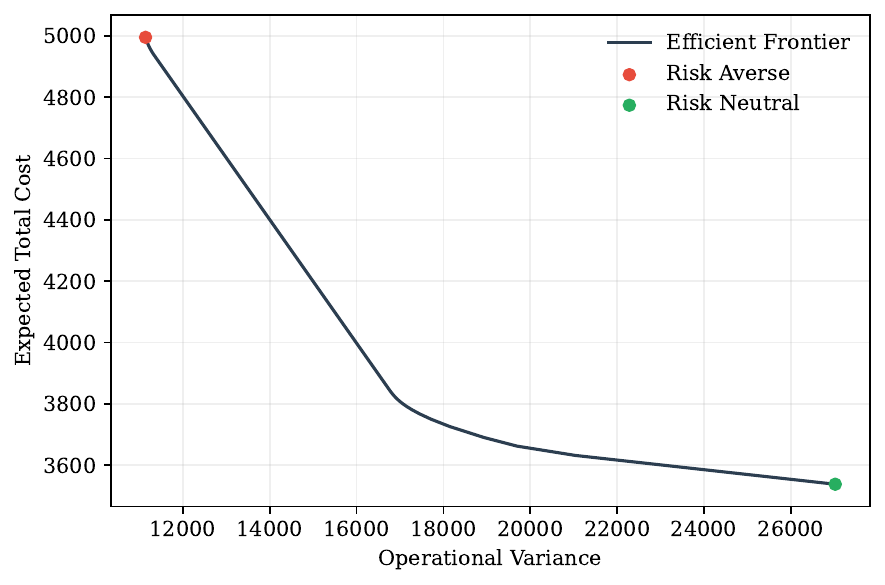}
    \caption{The Efficient Frontier for the warehouse location problem. The curve illustrates the trade-off between minimizing expected operational costs and reducing cost variance (risk). The green marker represents a risk-neutral approach ($\lambda=0$), while the red marker indicates increasing risk aversion.}
    \label{fig:efficient_frontier}
\end{figure}

 Similar to the discussion of the single warehouse model, as is, the model is only minimizing the cost in expected value (both in the objectuve and the first constraint). Variations in the travel times $T_{ij}$ will cause fluctuations in the cost. We can take a similar strategy to hedge this variation using a similar reasoning as in the single warehouse case. To do so, note that the variance $\var[\sum_{i\in W}\sum_{j\in V}c_{ij}q_{ij}]$ is equal to $\sum_{i\in W}\sum_{j\in V}\sum_{k\in W}\sum_{l\in V}q_{ij}\cov(c_{ij},c_{lk})q_{lk}$  which can be written in matrix form as $\bm{q}\Sigma\bm{q}^{\intercal}$ where $\Sigma$ has entries:
 \begin{eqnarray}
 	\cov(c_{ij},c_{lk})=\begin{cases}
 		\gamma^2_{ij}\var(T_{ij}), & \text{for } i,j=k,l \\
 		\gamma_{ij}\gamma_{kl}[\ex(T_{ij}T_{kl})-\ex(T_{ij})\ex(T_{kl})], & \text{for } i,j \neq k,l.
 	\end{cases}
 \end{eqnarray}
 Also using the chance constraint for the demand, the risk aware formulation of the problem can be written as a Mixed Integer Quadratic Program (MIQP):
 \begin{eqnarray}
 	&\underset{q_{ij}\ge 0,x_i\in\{0,1\}}{\min}& \ex\left\{ \sum_{i\in W}c_ix_i+\sum_{i\in W}\sum_{j\in V}c_{ij}q_{ij}\right\}+\lambda\var\left(\sum_{i\in W}\sum_{j\in V}c_{ij}q_{ij}\right)\\
 	&\text{s.t. }&\sum_{i\in W}q_{ij}\ge \ex (D_j)+\sqrt{\var(D_j)}\Phi^{-1}(1-\eta) \,\,\,\text{for}\,\,\,j\in V, \nonumber\\
 	&&\sum_{j\in V}q_{ij}\le x_iC_i \,\,\,\text{for}\,\,\,i\in W,\nonumber\\	
 	&&q_{ij}\le x_i\bar d_j\,\,\,\text{for}\,\,\,i,j\in W\times V. \nonumber
 \end{eqnarray}
 This formulation of the problem can be solved with a MIQP solver. Note that in this case the objective uses the variance instead of the standard deviation as in the single warehouse section. The problem can still be solved as a Mixed Integer Second-Order Conic Program (MISOCP). The formulation is also implicitly assuming that the travel times and the demands are independent from each other. In Figure \ref{fig:efficient_frontier} we can see the trade-off we mentioned earlier considering different values of $\lambda$.  
 
 \section{Distributionally Robust Warehouse Location Optimization}
 
 Readers of the first sections should be feeling slightly uneasy regarding the sensitivity of these optimization problems. By sensitivity, we mean how the optimal solution changes when the input data changes. For example, there is a lot of expected values for the demand $\ex(D_j)$, that as we have seen, when there is a small number of observations or it is in reality time-varying, the sample mean $\sum_{t=1}^Td_t$ might be within ax range of the true value. Compounding this effect across all expected values and variances (which are in turn expected values as well) creates an effect where the optimal solution of solving the warehouse problem with a sample observations $\hat q^*_{ij}$ might be substantially different from the true optimal solution $q^*_{ij}$.
 
 Given a fixed solution $q^*_{ij}$ we have explored previously how to build confidence intervals using the CLT (see Chapter \ref{sec:inventory_management} and the empirical distributions section to see some of these techniques). To get a feel for the interplay between the sample size $n$ and the number of vertices $N$ we use a different but useful way of building confidence intervals, we illustrate this with an example:
 \begin{exm}[{\bf Local fulfillment centers}]
 	A last-mile delivery company wants to know the optimal location to build their fulfillment centers to serve the New York metropolitan area. Geographically, this area is composed by 23 counties, which makes a graph with $N=23$ nodes. The company has a sample size of daily orders for each of the counties going over 4 years for a total of around $n=1,460$ observations. Their data-scientist solves the basic warehouse location problem in Equation (\ref{eq:warehouse_ex_cost}) obtaining a solution $x^*_i,q^*_{ij}$ for the optimal locations of the warehouses and the quantity that each should deliver to each county. The management of the company seems skeptical of the solution and ask the data-scientist for confidence intervals for their estimation. Specifically, they want to know with $1-\alpha=95\%$ confidence what is a range for the variable cost component of the cost (involving the quantities $q_{ij}$), as they know with certainty how much it costs to run each of the facilities once active.
 	
 	To solve this problem, we use a very useful result known as \emph{Hoeffding's Inequality} which states the following for a sum of independent random variables $S_n=X_1+\cdots+X_n$, where all variables are bounded, that is, $\pr(X_i\in [\ell,u])=1$. Then, we have\footnote{The more general version, states that when each random variable is bounded as $\pr(X_i\in [\ell_i,u_i])=1$ the denominator of the exponential is $\sum_{i=1}^n(u_i-\ell_i)^2$ instead of $n(u-\ell)^2$.}
 	\begin{eqnarray}
 		\pr\left(S_n-\ex[S_n]\ge \varepsilon\right)\le\exp\left(\frac{-2\varepsilon^2}{n(u-\ell)^2}\right),
 	\end{eqnarray}
	The trick to use this bound is to express the objective value as a function operating on a real valued function. We initially have a random vector for the demands $\bm D=(D_1,\dots,D_{N})$, given a solution $x=(x_1,\dots,x_M)$ and $\bm q=(q_{11},\dots,q_{MN})$ with cost vector $c=(c_{11},\dots,c_{MN})$. The value that we want to bound is $\ex[c\bm q^\intercal]$. The first observation is noting that for any solution, in practice is should be the case that $\sum_{i\in W}q_{ij}=D_j$ for all nodes $j\in V$. Then, by letting the vector $\bar c=(\bar c_1,\dots,\bar c_{N})$ where $\bar c_{j}=(\sum_{i\in W}c_{ij}q_{ij})/(\sum_{i\in W}q_{ij})$ for $j\in V$. In other words, the vector $\bar c_{j}$ is the average cost to deliver to each location. With this, we have $c \bm q^\intercal=\bar c\bm D^\intercal$. Now we can use Hoeffding's inequality to build our cost bound. The function $\bar c\bm D^\intercal$ can be seen as the generic random variable $X$ in the inequality. The first step is to determine its bound. Note that all demands are bounded below by $0$, that is, $\ell=0$. Likewise, each demand is bounded by above by $\bar d_j$, then, all of them are bounded by $\bar d=\max_{j\in V}\bar d_j$. Then, it follows that the cost is bounded as $\bar c\bm D^\intercal\le\lvert\bar c\rvert \bar d$ (we abuse notation by letting $|\cdot|$ be the $L-1$ norm for vectors. For example, $\lvert\bar c\rvert=\sum_{j\in V}|\bar c_j|$). This norm $\lvert\bar c\rvert $ is bounded by $N\max\{c_{ij}\}$ where $\max\{c_{ij}\}$ is the highest possible cost, the $N$ is simply the number of vertices in the graph. Then, the random variable $\pr(\bar c\bm D^\intercal\in[0,N\max\{c_{ij}\}\bar d])=1$. 
	
	An i.i.d sample of the demand vector $\bm D$ is represented as vectors $\bm d_k$ for $k=1,\dots,n$, where $\bm d_k$ is a vector with $N$ elements, the demand for each region. The sample expected cost is then $\ex_{\sf P_n}[\bar c\bm D]=\frac{1}{n}\sum_{k=1}^n\bar c \bm d_k$ while the true cost is the expected value $\ex[\bar c\bm D]$. Using Hoeffding's inequality we get that:
	\begin{eqnarray}
	\pr\left(\ex_{\sf P_n}[\bar c\bm D]-\ex[\bar c \bm D]\ge \epsilon\right)\le \exp\left(\frac{-2n^2\varepsilon^2}{n(N\max\{c_{ij}\}\bar d)^2}\right),
	\end{eqnarray} 
\end{exm}
\noindent With $1-\alpha=95\%$ probability, we get that the true cost is bounded by (solving for $\epsilon$):
 \begin{eqnarray}
 	\ex[\bar c\bm D]\le\ex_{\sf P_n}[\bar c\bm D]+  N\max\{c_{ij}\}\bar d\sqrt{\frac{1}{2n}\ln\left(\frac{1}{\alpha}\right)}.
 \end{eqnarray}
This means that solving the problem with a sample of $n$ observations is likely underestimating the true expected cost due to the uncertainty caused by the number of geographical locations $N$. To mitigate this uncertainty requires at least a quadratic increase in the sample size $n$. In the next section, we present a strategy to deal with solving an alternative optimization problem that deals with the fact that the true distribution is in a neighborhood of the true distribution.
 \subsection{A reformulation using duality and optimal transport}
 Back in the first chapter, we defined both the actual expected value of a random variable and their sample counterpart. A specially important random variable is the indicator function $\ind_{e}$ that is equal to 1 if event $e$ occurs, and 0 otherwise. $e$ could be for example the event that today rains, or that the demand exceeds a thousand units. The expected value of an indicator $\ex(\ind_e)$ is equal to the probability of the event $\pr(e)$. And the sample estimate is the same as we defined before, that is, $p_n(e)=\frac{1}{n}\sum_{i=1}^n \ind \{e_i\}$. For example, the sample probability that the demand is equal to 10, is the sample probability $p_n(D=10)=\sum_{i=1}^n \ind \{d_i=10\}/n$. Sample probabilities behave the same as expected values, that is, they converge to their true value if the distribution is i.i.d.
 
 For a random variable $X$ taking values in the set $\calx$ (for example, in the case of the demand random variable $D$, the demand takes values in the set of the natural numbers $\mathbb{N}$), the collection of sample probabilities $p_n(X=x)=\sum_{i=1}^n\ind \{X_i=x\} $ for all $x\in \calx$, forms an \emph{empirical distribution} of $X$. This distribution is typically denoted by $\pr_n$ instead of the true probability distribution $\pr$. The reader can already imagine that as the sample size $n$ increases, the empirical distribution will converge to the true one by the LLN, that is $\pr_n\rightarrow \pr$. It is also, expected that the empirical distribution of $n$ observations will overestimate the probability of some events, while also underestimating others. In fact, since both probability distributions add up to 1, the amount of overestimation is exactly the amount of underestimation (think of why this is the case).
 
\begin{figure}[htb]
    \centering
    \begin{subfigure}[b]{0.48\textwidth}
        \centering
        \begin{tikzpicture}[scale=0.8]
            \draw[->] (0,0) -- (5,0) node[right, font=\tiny] {Demand};
            \draw[->] (0,0) -- (0,3) node[above, font=\tiny] {Prob};
            
            \draw[blue!70, fill=blue!20, thick] (0.8,0) rectangle (1.1, 2.2);
            \draw[blue!70, fill=blue!20, thick] (2.5,0) rectangle (2.8, 1.8);
            \draw[blue!70, fill=blue!20, thick] (4.0,0) rectangle (4.3, 0.8);
            
            \draw[red, thick] (1.8,0) -- (1.8, 0.1) node[below, font=\tiny] {Zero Data};
            
            \node[font=\small\bfseries] at (2.5, 3.2) {Empirical $\mathbb{P}_n$};
        \end{tikzpicture}
        \caption{Sparse \& Fragile}
    \end{subfigure}
    \hfill
    \begin{subfigure}[b]{0.48\textwidth}
        \centering
        \begin{tikzpicture}[scale=0.8]
            \draw[->] (0,0) -- (5,0) node[right, font=\tiny] {Demand};
            \draw[->] (0,0) -- (0,3) node[above, font=\tiny] {Prob};
            
            \draw[thick, red!70, dashed, fill=red!5] plot [smooth, tension=0.7] 
                coordinates {(0.5,0.4) (1,1.5) (2,0.9) (3,1.2) (4,0.7) (4.5,0.2)} 
                -- (4.5,0) -- (0.5,0) -- cycle;
            
            \draw[->, >=stealth, red] (1.8, 1.5) -- (1.8, 0.9);
            \node[font=\tiny, red, align=center] at (1.8, 1.8) {Risk recognized};

            \node[font=\small\bfseries] at (2.5, 3.2) {Robust $\mathbb{P} \in \mathcal{P}_\delta$};
        \end{tikzpicture}
        \caption{Spread \& Robust}
    \end{subfigure}

    \caption{Comparison of fulfillment modeling. The empirical approach (left) ignores events not present in the sample, leading to $q_{ij}=0$ for potentially critical regions. The DRO approach (right) considers a distribution that spreads probability mass into those gaps, forcing the optimization to allocate resources ($q_{ij} > 0$) for unobserved but plausible demand.}
    \label{fig:empirical_vs_robust_side}
\end{figure}
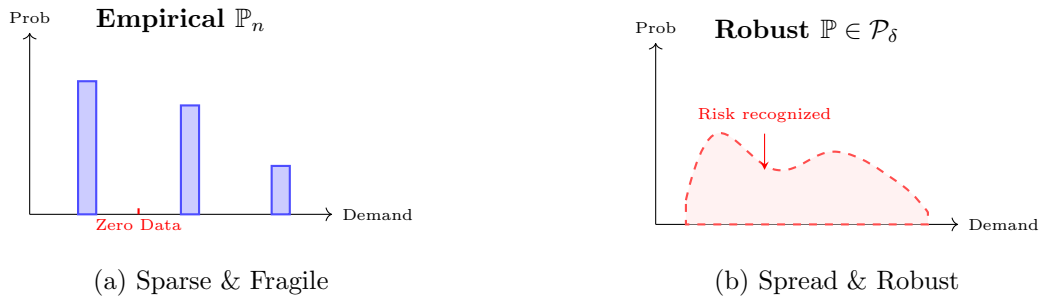

 Given these facts, the methodology known as \emph{Distributionally Robust Optimization} (DRO) aims to incorporate the idea that the true distribution $\pr$ of the data is within some neighborhood of the empirical distribution $\pr_n$. Here, neighborhood is taken as a \emph{distance metric} between the true and empirical distribution. Let $\cald$ denote a metric distance function between distributions\footnote{A distance metric is a function $\cald:\calm\times\calm\rightarrow \mathbb{R}$ that operates on a set $\calx$ with the properties: For $p_1,p_2,p_3\in \calm$ it is satisfied $\cald(p_1,p_1)=0$, $\cald(p_1,p_2)>0$ for $p_1\neq p_2$, $\cald(p_1,p_2)=\cald(p_2,p_1)$ for all $p_1,p_2\in\calm$ and $\cald(p_1,p_3)\le \cald(p_1,p_2)+\cald(p_2,p_3)$ which is the so-called triangle inequality.}. Then, we want to consider distributions $\pr$ that are around the observed empirical distribution $\pr_n$ within a distance $\delta$, that is, distributions $\pr:\cald(\pr_n,\pr)\le \delta$. In plain terms, a decision-maker acknowledges that their finite data samples might bias their estimates and wants to optimize considering the case where the true distribution is in some neighborhood of the empirical distribution $\pr_n$, in this case, neighborhood means the set of possible distributions that are in a \emph{ball} of radius $\delta$ around $\pr_n$. The decision-maker picks the radius according the their risk aversion (similar to choosing $\lambda$ in the previous examples). A larger radius $\delta$ means less trust in the empirical distribution (for example, due to a low sample size) while a small radius means that there is high confidence in the empirical distribution. See Figure \ref{fig:empirical_vs_robust_side} for a conceptual representation of the DRO approach.
 
  The DRO formulation can be set up in many ways. We will focus in a paradigm that uses duality theory and optimal transport. To simplify the problem we will only consider variation in the demand, that is, the costs $c_{ij}$ are assumed to be known and fixed. Moreover, the historical (sample) demand is the matrix $D_n=(\bm d_{k})$ for $k=1,\dots,n$, where $\bm d_k$ is a vector with the historical demands for each $i\in V$, that is, $\bm d_k$ is a vector of size $N$. Each of these vectors is an i.i.d. sample of the vector $\bm D$. Denote the set of distributions within $\delta$ distance of the empirical distributions as $\calp_\delta=\{\pr:\cald(\pr_n,\pr)\le \delta\}$. Let $x=(x_1,\dots,x_M)$. First, we re-write the warehouse problem in this setting as:
 \begin{eqnarray}
 &\min& \left\{ \tilde c x^{\intercal}+\sup_{\calp_\delta}\ex\left[\min c\bm q^{\intercal}|A\bm q^{\intercal}+Bx^{\intercal}\ge (\bm D,\bm 0)^{\intercal} \right]\right\},
 \end{eqnarray}
This formulation is a matrix form re-write of the problem of interest that is going to be useful for some algebraic manipulations that we are performing. For example, $\tilde c x^{\intercal}$ is equal to $\sum_{i\in W}c_ix_i$ and $c\bm q^{\intercal}=\sum_{i\in W}\sum_{j\in V}c_{ij}q_{ij}$. The constraints $A\bm q^{\intercal}+Bx^{\intercal}\ge (\bm D,\bm 0)^{\intercal} $ are the usual constraints of the problem, such as demand satisfaction, e.g. see the constraints in (\ref{eq:warehouse_ex_cost}), $A$ would be the matrix of coefficients multiplying the $\bm q$ variables, while the right hand side is the random demand $\bm D$ for the demand satisfaction constraints and $\bm 0$ for the constraints that only depend on variables. The interpretation of the problem is quite similar: minimize the overall cost of building the facilities, subject to the worst distribution $\calp_\delta$. Typically, the intuitive explanation is that nature chooses the worst distribution given $x$ (the $\sup$ part) and then the optimization chooses the optimal quantities $\bm q$.

\begin{figure}[htb]
    \centering
    \begin{subfigure}[b]{0.48\textwidth}
        \centering
        \begin{tikzpicture}[scale=0.8]
            \begin{axis}[
                axis lines=middle,
                xlabel={$v$}, ylabel={$y$},
                xmin=-3, xmax=3, ymin=-0.5, ymax=3,
                title={Uncertainty Cost: $\lambda |v|$},
                ytick=\empty,
                label style={font=\small},
                title style={font=\small\bfseries}
            ]
                \addplot[thick, blue, domain=-3:3, samples=100] {abs(x)};
                \node[blue, anchor=west] at (axis cs: 1, 1.8) {Risk Penalty};
                
                \draw[dashed, gray] (axis cs: -3, -1.5) -- (axis cs: 3, 1.5) node[pos=0.9, below, font=\tiny] {$\eta_N < \lambda$};
                \draw[dotted, red, thick] (axis cs: -1.5, -3) -- (axis cs: 1.5, 3) node[pos=0.8, left, font=\tiny] {$\eta_N > \lambda$};
            \end{axis}
        \end{tikzpicture}
        \caption{Cost of deviation $v$}
    \end{subfigure}
    \hfill
    \begin{subfigure}[b]{0.48\textwidth}
        \centering
        \begin{tikzpicture}[scale=0.8]
            \begin{axis}[
                axis lines=middle,
                xlabel={$\eta_N$}, ylabel={$f^*(\eta_N)$},
                xmin=-2, xmax=2, ymin=-0.5, ymax=3,
                title={Dual Constraint: $|\eta_N| \leq \lambda$},
                xtick={-1, 1}, xticklabels={$-\lambda$, $\lambda$},
                ytick={0},
                label style={font=\small},
                title style={font=\small\bfseries}
            ]
                \draw[ultra thick, blue] (axis cs:-1,0) -- (axis cs:1,0);
                \draw[ultra thick, blue, dashed, ->] (axis cs:-1,0) -- (axis cs:-1,2.5);
                \draw[ultra thick, blue, dashed, ->] (axis cs:1,0) -- (axis cs:1,2.5);
                
                \node[blue, font=\footnotesize] at (axis cs: 0, 0.3) {Feasible Region};
                \node[red, font=\footnotesize] at (axis cs: 1.5, 1.5) {$\infty$};
                \node[red, font=\footnotesize] at (axis cs: -1.5, 1.5) {$\infty$};
            \end{axis}
        \end{tikzpicture}
        \caption{Shadow Price "Discounting"}
    \end{subfigure}

    \caption{The relationship between deviation cost and shadow prices. In (a), if the shadow price $\eta_N$ (the marginal cost of demand) is steeper than the risk penalty $\lambda$, the sup-problem becomes unbounded. In (b), this translates to a hard constraint: the robust optimizer only considers shadow prices "discounted" within the $[-\lambda, \lambda]$ interval, effectively capping the system's sensitivity to demand uncertainty.}
    \label{fig:conjugate_intuition}
\end{figure}
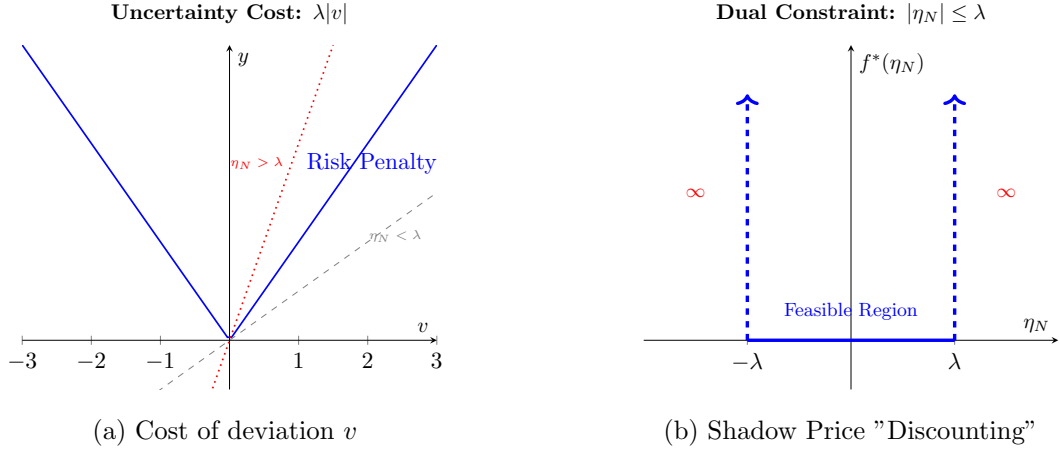

The first step to solve this problem is to re-formulate the $\sup$ by exploiting duality and results from Transport Theory. Letting $Q(x,\bm D)=\min c\bm q^{\intercal}|A\bm q^{\intercal}+Bx^{\intercal}\ge (\bm D,\bm 0)^{\intercal}$, we are going to write the Lagrangian of $\sup_{\calp_\delta}\ex[Q(x,\bm D)]$ (noting that the definition of $\calp_\delta$ is a constraint defined by the distance of the distributions), for a pair of variables $\pr,\lambda$:
\begin{eqnarray}
\sup_{\calp_\delta}\ex[Q(x,\bm D)]=\min_{\lambda\ge 0}\sup_{\sf P}\{\ex_{\sf P} [Q(x,\bm D)] + \lambda (\delta-\cald(\pr,\pr_n))\},
\end{eqnarray}
The next step is to choose a distance metric $\cald$ that makes the problem tractable. One such metric, is the so-called \emph{Wasserstein-1} distance. Choosing this distance $\cald=W_1$ has the following property for any function $g(\bm D)$:
\begin{eqnarray}
\sup_{\sf P}\left[\ex_{\sf P} [g(\bm D)]  -\lambda W_1(\pr,\pr_n))\right]=\ex_{\sf P_n}\left[ \sup_{\bm D}\{g(\bm D)-\lambda|\bm D-\bm D_k|\}\right],
\end{eqnarray} 
Where $D_k$ is a generic sample of the empirical distribution $P_n$ and $|\cdot|$ for vectors denotes the $L-1$ norm. We can let $g(\bm D)=Q(x,\bm D)$. The math behind the result is somewhat technical and outside the scope of this section. By letting $g(\bm D)=Q(x,\bm D)$ and rearranging $\min_{\lambda\ge 0}\sup_{\sf P}\{\ex_{\sf P} [Q(x,\bm D)] + \lambda (\delta-\cald(\pr,\pr_n))\}$ as $\min_{\lambda\ge 0}( \lambda\delta + \sup_{\sf P}\{\ex_{\sf P} [Q(x,\bm D)]-\lambda\cald(\pr,\pr_n)\})$ combining the two equations yields:
\begin{eqnarray}
	\sup_{\calp_\delta}\ex[Q(x,\bm D)]=\min_{\lambda\ge 0}\left\{ \lambda\delta +\ex_{\sf P_n}\left[ \sup_{\bm D}\{Q(x,\bm D)-\lambda|\bm D-\bm D_k|\}\right]\right\},
\end{eqnarray}
Which is great news as we know how to estimate expected values in the empirical distribution, combining the original objective we get:
\begin{eqnarray}
	\min_{\lambda\ge 0}\left\{ \tilde c x^{\intercal} + \lambda\delta +\frac{1}{n}\sum_{k=1}^n\sup_{\bm D}\{Q(x,\bm D)-\lambda|\bm D-\bm d_k|\}\right\},
\end{eqnarray}
Recalling $Q(x,\bm D)=\min_{\bm q\ge 0} c\bm q^{\intercal}|A\bm q^{\intercal}+Bx^{\intercal}\ge (\bm D,\bm 0)^{\intercal}$, we can calculate the dual of this problem as $\max_{\bm \eta\ge 0} \bm\eta [(\bm D,\bm 0)^{\intercal}-Bx^\intercal]\,|\,\bm\eta A\le c$. Each of the terms $\sup_{\bm D}\{Q(x,\bm D)-\lambda|\bm D-\bm d_k|\}$ can be then written as $\sup_{\bm D}\{ \bm\eta [(\bm D,\bm 0)^{\intercal}-Bx^\intercal]-\lambda|\bm D-\bm d_k|\}$. Exchanging the sup and the max and isolating the terms depending on the vector $\bm D$ yields:
\begin{eqnarray}
\max_{\bm \eta\ge 0,\bm\eta A\le c} -\bm\eta Bx^\intercal+ \sup_{\bm D}\{\bm\eta_N\bm D^{\intercal}-\lambda|\bm D-\bm d_k|\},
\end{eqnarray}
Where $\bm\eta_N$ are the first $N$ elements of $\bm\eta$ corresponding to the dual variables of the demand satisfaction constraints of the warehouse problem. Next, define the variable $\bm v_k = \bm D-\bm d_k$. The inner problem $\sup_{\bm D}\{\bm\eta_N\bm D^{\intercal}-\lambda|\bm D-\bm d_k|\}$ can be re-written as $\bm\eta_N\bm d^\intercal_k+\sup_{\bm v_k}\{\bm\eta_N\bm v_k^{\intercal}-\lambda|\bm v_k|\}$. The last term is the convex conjugate of the $L-1$ norm and is equal to 0 when $|\text{max}\, \bm\eta_N|\le \lambda $ or equal to $\infty$ otherwise. That it is equal to 0 is achieved by the constraints $-\bm 1\lambda\le \bm\eta_N\le\bm 1\lambda$. See Figure \ref{fig:conjugate_intuition} for the intuition behind these constraints. 

With this we can finally finish the rewrite of the problem as:
 \begin{eqnarray}
&\min_{\lambda,\bm\eta,x}&\left\{ \tilde c x^{\intercal} + \lambda\delta +\frac{1}{n}\sum_{k=1}^n(\bm\eta_N\bm d^\intercal_k-\bm\eta B x^\intercal)\right\}\\
	&\text{s.t. }&\bm\eta A\le c,-\bm 1\lambda\le \bm\eta_N\le\bm 1\lambda,\nonumber\\
	&&\lambda\ge 0,\bm \eta\ge 0,x\in \{0,1\}^M.\nonumber
 \end{eqnarray}
Which is almost a MIP except for the variables $\bm\eta B x^\intercal$. The standard way of dealing with this, is to define new variables $z_{ij}=\eta_ix_j$ and add the constraints $z_{ij}\le \tilde Mx_j,z_{ij}\le\eta_i,z_{ij}\ge\eta_i-\tilde M(1-x_j),z_{ij}\ge 0$ for a large enough $\tilde M$ (this is commonly known as the big-M method).

 \chapter{Scheduling Models} 
\label{sec:scheduling}
\section{Introduction}
Matching supply and demand is commonly referred to as an atemporal event,  as if markets clear instantaneously. Even in the most sophisticated financial markets, matching supply and demand is far from an instantaneous affair where order matters and makes a significant difference in the outcome. For example, the price difference between transactions can amount to very large amounts of money in a matter of milliseconds. 

In physical markets, the problem is even more pronounced. For example, in the warehouse setting considered in the previous section, demand is fulfilled in a specific order that if not properly optimized can lead to congestion, lost sales and waste of resources. Many of these problems are very simple to describe but difficult to solve: the order of packaging boxes in a warehouse, the order of manufacturing different products, or deciding the schedules of workers to fulfill a productive activity.

Similar to other chapters when we studied the demand, scheduling can be studied as a time-space matching problem. We start by studying what would be the \emph{space} variant (which could also be seen as discrete time). Then, we will move to consider the problem when time is a random variable and finally consider both problems simultaneously in future chapters.

\section{A set covering problem}
A lot of scheduling problems take the form of an abstract problem: to \emph{cover} a set at minimum cost. Here, set can be understood in a general way as a demand that extends over either time or space, and certain resources are available to satisfy it. For example, a fixed number of trucks need to be assigned to cover deliveries in a city, or workers need to be assigned tasks to complete a production process. Many of these problems are combinatorial in nature, meaning that finding optimal solutions is in general hard, except in the simplest cases.

\begin{figure}[htbp]
    \centering
    \begin{minipage}{\linewidth}
        \centering
        \begin{tikzpicture}[
            >=Stealth,
            node distance=1.2cm,
            box/.style={draw, rectangle, minimum height=0.8cm, minimum width=1.5cm, align=center, rounded corners=2pt, fill=blue!5},
            matrix style/.style={matrix of math nodes, nodes in empty cells, left delimiter=[, right delimiter=], row sep=0.1cm, column sep=0.1cm, nodes={minimum width=0.6cm, minimum height=0.6cm, anchor=center}},
            scale=0.9, transform shape 
        ]

            \matrix [matrix style, fill=green!5] (A) at (0,0) {
                a_{11} & a_{12} & \dots & a_{1J} \\
                a_{21} & a_{22} & \dots & a_{2J} \\
                \vdots & \vdots & \ddots & \vdots \\
                a_{N1} & a_{N2} & \dots & a_{NJ} \\
            };
            
            \node[above=1.1cm of A, font=\small\bfseries] {Incidence Matrix $A$ ($N \times J$)};
            
            \draw[decorate, decoration={brace, amplitude=5pt}, thick, gray] 
                ([yshift=5pt]A.north west) -- ([yshift=5pt]A.north east) 
                node[midway, above=6pt, black, font=\footnotesize] {$J$ Schedules ($\mathcal{S}$)};
                
            \draw[decorate, decoration={brace, amplitude=5pt, mirror}, thick, gray] 
                ([xshift=-5pt]A.north west) -- ([xshift=-5pt]A.south west) 
                node[midway, left=8pt, black, font=\footnotesize] {$N$ Tasks};

            \node[box, fill=orange!10, right=2.5cm of A] (x) {$\bm{x} = \begin{bmatrix} x_1 \\ x_2 \\ \vdots \\ x_J \end{bmatrix}$};
            \node[above=0.2cm of x, font=\small\bfseries, align=center] {Decision Vector $\bm{x}$ \\ {\footnotesize (Selected Schedules)}};

            \node[right=1.2cm of x] (ge) {\LARGE $\ge$};
            
            \node[box, fill=red!5, right=1.0cm of ge] (b) {$b\ex(\bm{D}) = \begin{bmatrix} b\ex(D_1) \\ b\ex(D_2) \\ \vdots \\ b\ex(D_N) \end{bmatrix}$};
            \node[above=0.2cm of b, font=\small\bfseries, align=center] {Demand / Bounds \\ {\footnotesize (Time/Space Over $\bm{D}$)}};

            \draw[->, thick, dashed] (A) -- node[midway, above=2pt, font=\footnotesize] {Multiply} (x);
            \draw[->, thick] (x) -- (ge);
            \draw[<-, thick] (b) -- (ge);

            \node[draw=red!70!black, dotted, inner sep=0.5cm, fit=(A) (x) (ge) (b), label={[red!70!black, font=\small\itshape]below:Set Covering Constraint: $Ax \ge b(\bm{D})$}] {};

        \end{tikzpicture}
    \end{minipage}
    \caption{Structural breakdown of the matrix-vector alignment in the set covering scheduling formulation, where available schedules ($\mathcal{S}$) map onto structural demands over a space-time domain.}
    \label{fig:scheduling_set_covering}
\end{figure}
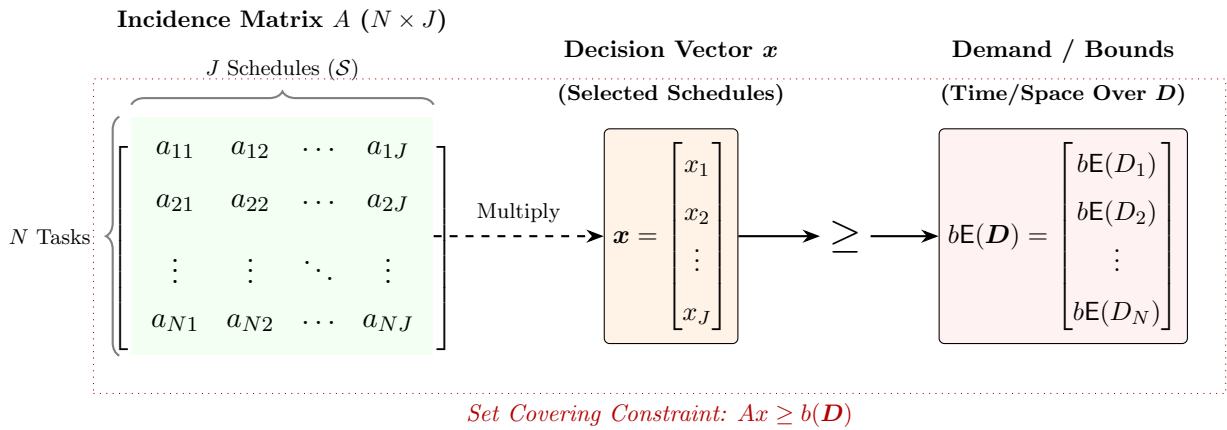

Like in previous chapters, we start with the simplest formulation that we will extend to other special cases. Consider a demand vector $\bm D$, that like in previous chapters can represent the demand for multiple locations or time-periods, the demand commonly has a task interpretation, let the dimension of the demand vector be $N$. Let $\cals$ be a set of different $J$ schedules available to fulfill some task. Moreover, define the matrix $A$ with dimensions $N\times J$, and elements $a_{ij}=1$ if task $i$ is covered by schedule $j$. Let $x$ be a decision variable (a vector of size $J$) in some space of feasible assignments $\calx$ (for example, $\calx=\{0,1\}^{J}$). The scheduling problem normally takes the form:
 \begin{eqnarray}
	&\min&f(x)\\
	&\text{s.t. }&Ax^\intercal\ge b\ex(\bm D),\nonumber\\
	&&x\in\calx.\nonumber
\end{eqnarray}
See Figure \ref{fig:scheduling_set_covering} to see the structure of the problem.

We illustrate next with a couple of examples the flexibility of this formulation.

\begin{exm}[{\bf Shift Scheduling}]
	A warehouse wants to schedule workers to process packages in a 24 hour basis. Every hour, the warehouse has a demand of packages represented by a vector $\bm D=(D_0,\dots,D_{23})$. The warehouse has a headcount of $H$ workers that it has to allocate for 8-hour consecutive shifts. That is, there are possible 24 different shifts each starting every hour. A worker can process $b$ packages per hour, therefore the need of workers is in expected value $b\ex(\bm D)$. Let $A$ be a matrix with 24 columns and 24 rows and entries $a_{tj}=1$ if schedule $j$ covers hour $t$. Moreover, let $x$ be a vector with $J=24$ elements denoting the headcount of workers assigned to each schedule.
	The manager wants to maximize the expected buffer between supply and demand across all hours of the day. 
	
The first step to solve the problem is to define the objective function $f(x)$. In this case, the gap between supply and demand at every hour is $Ax^\intercal-b\ex(\bm D)$, as $Ax^\intercal$ is the number of scheduled workers every hour. Then, the minimum buffer between supply and demand is $v=\min_t \{(Ax^\intercal)_t-b\ex(D_t)\}$. As we want to make this as large as possible, $f(x)=-v$ achieves the goal of increasing the gap between supply and demand as much as possible between supply and demand. Finally, the constraint $|x|= H$ achieves the goal of assigning all workers (recall that for vectors $|x|=\sum_{j=1}^J|x_j|$). The final formulation results in the following MIP:
\begin{eqnarray}
	&\max&v\\
	&\text{s.t. }&Ax^\intercal-b\ex(\bm D)\ge v,\nonumber\\
	&&v\in\mathbb{R},x\in\mathbb{N}^J, |x|=H.\nonumber
\end{eqnarray}

As a follow-up, the business managers want to know instead the minimum amount of workers such that they can cover their operation in expected value.

This can be achieved by modifying the function $f(x)$ to $|x|$, that we previously shown is the total number of workers. Also, the constraint that represents that the demand is covered at all hours in expected value is $Ax^\intercal\ge b\ex(\bm D)$. Leading to the alternative formulation:
\begin{eqnarray}
	&\min&|x|\\
	&\text{s.t. }&Ax^\intercal\ge b\ex(\bm D),\nonumber\\
	&&x\in\mathbb{N}^J.\nonumber
\end{eqnarray}
The last example can also be easily modified when each schedule has a cost vector $c$ associated to each entry of $x$, and minimizing $f(x)=cx^\intercal$. When the cost of all schedules is the same, or the cost is proportional to number of hours (at the same rate), then the optimal value $x^*$ is the same as the previous problem with $f(x)=|x|$.
\end{exm}

\subsection{Data Uncertainty}
We revisit a technique in the warehouse allocation problem to deal with uncertainty in the data matrix $A$. In the previous formulation. The matrix $A$ had entries $a_{ij}=1$ whenever schedule $j$ is covering the period or location $i$. The $1$ has a direct meaning that every resource allocated on that schedule is going to \emph{certainly} be allocated. In real life this is rarely the case as there are many occasions when an allocated resource might not fully work: in the case of people it can be because of a curve of productivity (that is not uniform across all hours of the shift). In the case of a machine it can be because of a failure that occurs with a small probability. Here we provide a data-driven framework to deal with these kind of uncertainty:

Consider each row of the matrix $A$ as a multi-variate normal random variable $a_i$ with mean $\bar a_i$ and covariance matrix $\Sigma_i$. For example, in the productivity curve case where each hour of the shift has their own productivity level $a_i\in [0,1]^J$ is a random vector with the different productivity levels. The mean and the variance could be estimated from a sample size of $n$ observations for example. Each constrain for the demand satisfaction can be transformed into the probabilistic constraint:
\begin{eqnarray}
\pr(a_ix^\intercal \ge b D_i)\ge 1-\eta,
\end{eqnarray}
Assuming the demand $D_i$ is another independent random variable, the linear combination $a_ix^\intercal-bD_i$ is another random variable. The mean of this random variable is $\bar a_i x^\intercal-b\ex(D_i)$ and its variance is $x\Sigma_i x^\intercal+b^2\var(D_i)$. Using the usual trick with the inverse of the normal random variable trick seen in other chapters, we have that $\pr(a_ix^\intercal\ge bD_i)$ is equal to $\pr(Z\ge -(\bar a_ix^\intercal-b\ex(D_i))/\sqrt{x\Sigma_ix+b^2\var(Di)})$ or $1-\pr(Z\le -(\bar a_ix^\intercal-b\ex(D_i))/\sqrt{x\Sigma_ix+b^2\var(D_i)})$. Making the inequality equal to $\pr(Z\le -(\bar a_ix^\intercal-b\ex(D_i))/\sqrt{x\Sigma_ix+b^2\var(D_i)})\le\eta$. Taking the inverse cdf of the normal yields the equivalent constraint:
\begin{eqnarray}
\bar a_ix^\intercal+\Phi^{-1}(\eta)\sqrt{x\Sigma_ix+b^2\var(D_i)}\ge b\ex (D_i),
\end{eqnarray}
Thematically, this constraint follows a similar pattern through the book that hedging the risk of allocating resources not only involves demand risk, but also the risk associated to allocating the resources.

Next, we need to transform this constraint into a tractable formulation as the term $x\Sigma_ix^\intercal$ is not linear in $x$, thus, not directly a MIP. A powerful paradigm for these kind of problems is known as \emph{Second-order Cone Programming} or SOCP. Note that $\sqrt{x\Sigma_ix^\intercal+b^2\var(D_i)}$ is equivalent to $||(\Sigma^{1/2}_i x^\intercal,\sqrt{\var(b D_i)})||_2$, where $\Sigma^{1/2}_i$ is a matrix such that $\Sigma^{1/2}_i(\Sigma^{1/2}_i)^\intercal=\Sigma_i$, for example, think of Cholesky decomposition. Using this fact, the original scheduling program can be rewritten as:
 \begin{eqnarray}
	&\min&f(x)\\
	&\text{s.t. }&-\Phi^{-1}(\eta)\left\|\begin{pmatrix}
\Sigma^{1/2}_i x^\intercal \\
\sqrt{\var(b D_i)}
\end{pmatrix}\right \|_2\le \bar a_ix^\intercal-b\ex(D_i)\text{   for   }i=1,\dots,N,\nonumber\\
	&&x\in\calx.\nonumber
\end{eqnarray}
When $f(x)$ is an affine function, this optimization program can be solved by a MISOCP solver. This formulation is flexible to accomodate an interesting range of applications: for example, a finite sample observation of $\bar a_i$

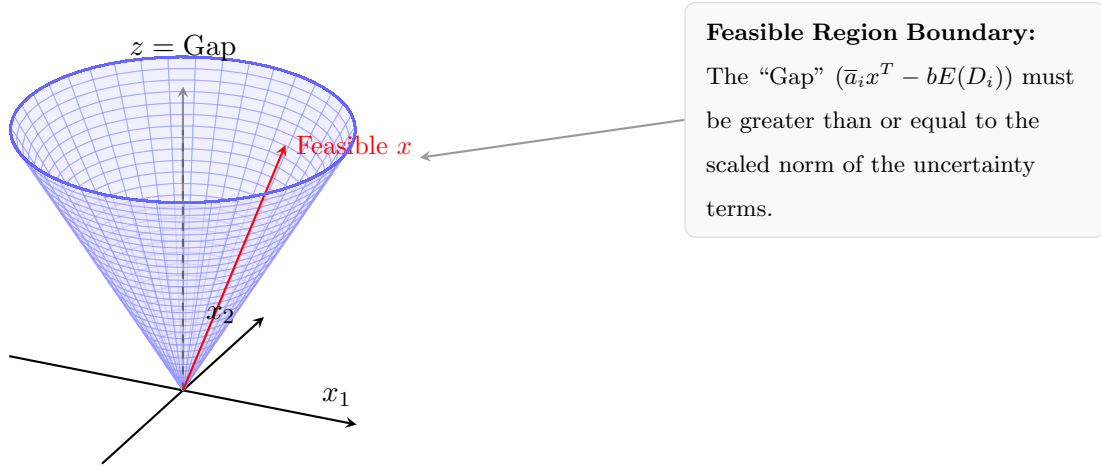
\begin{figure}[htbp]
    \centering
    \begin{tikzpicture}
        \begin{axis}[
            view={25}{25}, 
            axis lines=center,
            width=12cm, 
            height=10cm,
            axis equal image,
            xmin=-4, xmax=4,
            ymin=-4, ymax=4,
            zmin=0, zmax=7, 
            xlabel={$x_1$},
            ylabel={$x_2$},
            zlabel={$z = \text{Gap}$},
            zlabel style={at={(ticklabel* cs:1.05)}, anchor=south},
            ticks=none,
            axis line style={-stealth, thick},
            clip=false
        ]
            \addplot3 [
                surf,
                shader=flat,
                fill=blue!10,
                draw=blue!40,
                opacity=0.6,
                domain=0:6, 
                y domain=0:360,
                samples=30,
                samples y=40,
                z buffer=sort
            ] ({0.6*x*cos(y)}, {0.6*x*sin(y)}, {x});

            \draw[dashed, thick, black!60] (axis cs:0,0,0) -- (axis cs:0,0,6);

            \draw [red, thick, -stealth] (axis cs:0,0,0) -- (axis cs:1.8, 1.2, 5.5) 
                node[anchor=west, font=\small] {Feasible $x$};

            \addplot3 [blue!60, thick, domain=0:360, samples=60] ({0.6*6*cos(x)}, {0.6*6*sin(x)}, {6});
            
            \node (boundary) at (axis cs:2.2, 1.5, 4.5) {};
        \end{axis}

        \node [anchor=west, text width=5cm, font=\footnotesize, draw=gray!30, fill=gray!5, rounded corners, inner sep=8pt] (note) at (10, 5) {
            \textbf{Feasible Region Boundary:}\\
            The ``Gap" ($\overline{a}_i x^T - bE(D_i)$) must be greater than or equal to the scaled norm of the uncertainty terms.
        };
        
        \draw [-stealth, gray!80, thick] (note.west) -- (6.5, 4.5);

    \end{tikzpicture}
    \label{fig:cone}
    \caption{The Lorentz Cone representation. The vertical axis represents the supply-demand gap, while the interior represents the feasible set for the SOCP constraint.}
\end{figure}
%
%
%
%
%
%

\section{Production Scheduling}
We change pace by describing a different kind of scheduling models closely related to closest route problems, and that provide a useful : Suppose a production planning problem over a horizon of  $T$ periods in the future. We start with the simplest case when the demand is known, that is, at every period the demands $d_1,d_2,\dots,d_T$ are known in advance. Define $c(s,t)$ as the production and holding costs of periods $s,\dots,t$ for $s\le t\in\{0,\dots,T\}$ (inclusive of all periods in the range and the starting and endpoint) so that demand is satisfied exactly. 

For example, $c(0,4)$ includes the setup cost for production, plus the variable cost of producing for demands $d_1,\dots,d_4$ and the cost of storing the demand all the way to period $4$ (starting production at time 0). The key observation is that once this production is finished, starting at period 5, the system \emph{renews} and the past production decision can be ignored. As an example, this cost can be a setup cost $K_0$ of starting production at period 0, plus a variable cost $\sum_{t=1}^4v_0(d_t)$ and holding costs $h_0(d_1,\dots,d_4)$ until the respective demand period. Then, $c(0,4)$ could be:
\begin{eqnarray}
c(0,4)=K_0+\sum_{t=1}^4v_0(d_t)+h_0(d_1,\dots,d_4),
\end{eqnarray}
Nonetheless, this is not the only possible specification for the costs. The important part is that all these costs are deterministic and can be computed.

This structure of the problem allows us to use a powerful paradigm to solve the problem. Let $c(t)$ be defined as the minimum cost of meeting demand for periods $t,\dots, T$. Then, we have the following recursion:
\begin{eqnarray}
c(s)=\min_{t:t\ge s}[c(s,t)+c(t+1)],
\end{eqnarray}
This recursion implies a similar structure to a shortest-path-problem in a graph with nodes $V=\{0,1,\dots,T\}$ and edges $E=\{(s,t):s\le t\in\{0,\dots,T\}\}$ with distances $c(s,t)$. The technique to solve the problem is to start with $c(T+1)=0$, and solve recursively backwards until arriving to $c(0)$. The solution is a shortest-path between nodes $0$ and $T$ of the graph. Each node visited represents the periods where production should take place. When the demands are random and the costs have non-linearities. The model reduces to the dynamic inventory management case studied in Chapter \ref{sec:inventory_management}.

\begin{figure}[htbp]
    \centering
    \begin{minipage}{\linewidth}
        \centering
        \begin{tikzpicture}[
            >=Stealth,
            node distance=2cm,
            node/.style={circle, draw=blue!70, fill=blue!5, thick, minimum size=0.9cm, font=\bfseries},
            pathnode/.style={circle, draw=orange!90, fill=orange!10, ultra thick, minimum size=0.9cm, font=\bfseries},
            edge/.style={->, thick, color=gray!60},
            pathedge/.style={->, ultra thick, color=orange!90}
        ]

            \node[pathnode] (n0) at (0,0) {0};
            \node[node]     (n1) at (2.2,0) {1};
            \node[node]     (n2) at (4.4,0) {2};
            \node[pathnode] (n3) at (6.6,0) {3};
            \node[node]     (n4) at (8.8,0) {4};
            \node[pathnode] (n5) at (11.0,0) {5}; 

            
            \draw[edge] (n0) to[bend left=30] node[midway, above=-1pt, font=\scriptsize, color=black] {$c(0,0)$} (n1);
            \draw[edge] (n0) to[bend left=40] node[midway, above=-1pt, font=\scriptsize, color=black] {$c(0,1)$} (n2);
            \draw[pathedge] (n0) to[bend left=50] node[midway, above=1pt, font=\small\bfseries, color=orange!90!black] {$c(0,2)$} (n3);
            \draw[edge] (n0) to[bend left=60] node[midway, above=-1pt, font=\scriptsize, color=black] {$c(0,3)$} (n4);
            \draw[edge] (n0) to[bend left=70] node[midway, above=-1pt, font=\scriptsize, color=black] {$c(0,4)$} (n5);

            \draw[edge] (n1) to[bend right=30] node[midway, below=-2pt, font=\scriptsize, color=black] {$c(1,1)$} (n2);
            \draw[edge] (n1) to[bend right=40] node[midway, below=-2pt, font=\scriptsize, color=black] {$c(1,2)$} (n3);
            \draw[edge] (n1) to[bend right=50] node[midway, below=-2pt, font=\scriptsize, color=black] {$c(1,3)$} (n4);

            \draw[edge] (n2) to[bend left=30] node[midway, above=-1pt, font=\scriptsize, color=black] {$c(2,2)$} (n3);
            \draw[edge] (n2) to[bend left=40] node[midway, above=-1pt, font=\scriptsize, color=black] {$c(2,3)$} (n4);

            \draw[edge] (n3) to[bend right=30] node[midway, below=-2pt, font=\scriptsize, color=black] {$c(3,3)$} (n4);
            \draw[pathedge] (n3) to[bend right=45] node[midway, below=1pt, font=\small\bfseries, color=orange!90!black] {$c(3,4)$} (n5);

            \draw[edge] (n4) to[bend left=30] node[midway, above=-1pt, font=\scriptsize, color=black] {$c(4,4)$} (n5);

            \node[below=0.4cm of n0, font=\small\bfseries, color=orange!90!black] (startlbl) {Start ($s=0$)};
            \node[below=0.4cm of n5, font=\small\bfseries, color=orange!90!black] (termlbl) {Terminal ($c(5)=0$)};
            
            \node[draw=orange!90, fill=orange!5, dashed, rounded corners, below=2.2cm of n2.center, xshift=1.1cm, inner sep=0.3cm] (legend) {
                \begin{minipage}{11cm}
                    \centering\small\itshape\color{orange!90!black}
                    Optimal Path Solution shown in Orange: Produce at Period 0 (covers demands up to Period 2), then renew and produce at Period 3 (covers up to Period 4).
                \end{minipage}
            };

        \end{tikzpicture}
    \end{minipage}
    \caption{Production scheduling mapped as a shortest path network over a horizon where $T=4$. Edge weights $c(s,t)$ represent setup, variable production, and holding costs. Visited nodes on the optimal shortest path (highlighted in orange) signal the optimal periods where regeneration takes place and production runs must occur.}
    \label{fig:production_shortest_path}
\end{figure}
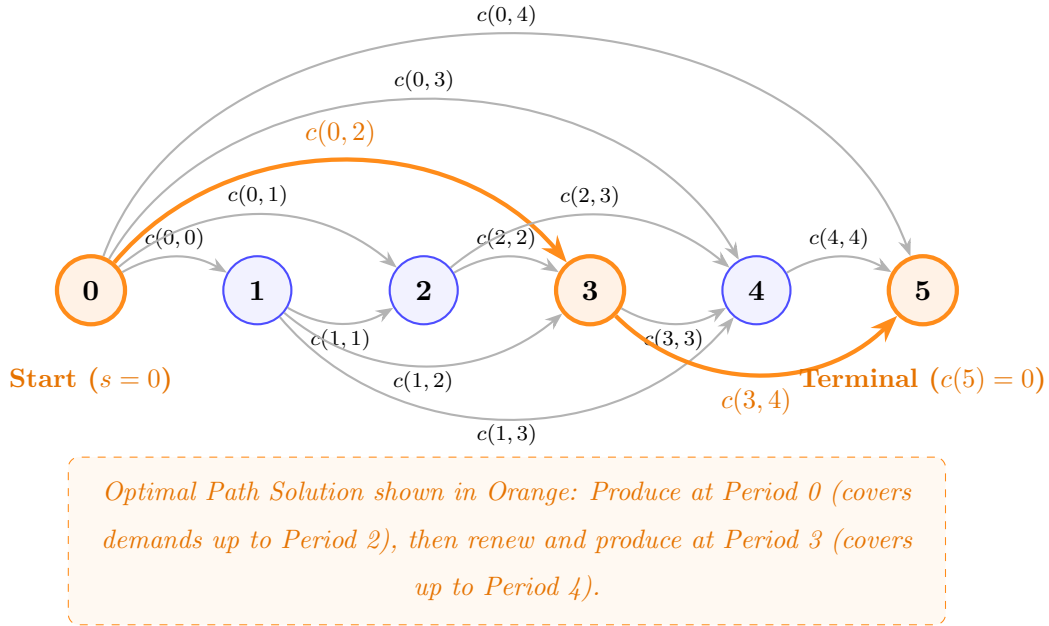

\section{A knapsack problem}
Lastly, we cover a problem that while not necessarily deals with scheduling in a nominal way, it is closely related as it deals with fitting objects in a limited space and can be interpreted in many ways: for example choosing jobs to perform in a limited time, or choosing investments subject a budget. This problem is known as the \emph{Knapsack Problem}. The most common variations of this problem can be formulated in the following setting. Imagine a knapsack (or a warehouse) where a series of objects $i=1,\dots,N$ has a profit or benefit expressed as $f_i(x_i)$ and a capacity of the knapsack consumed $c_i(x_i)\ge 0$. This capacity is assumed to be positive and increasing in $x_i$. $x_i\in\calx$ either expresses a quantity or indicator of how much of each object is used in the knapsack.  The objective is to pack the knapsack with the objects that have maximum benefit, subject to not exceeding the total capacity of the knapsack.

\begin{eqnarray}
& \max& \sum_{i=1}^Nf_i(x_i)\\
  &\text{s.t.}& \sum_{i=1}^Nc_i(x_i)\le C,\nonumber\\
  && x_i\in\calx\,\,\,\text{for}\,\,\,i=1,\dots,N.\nonumber
\end{eqnarray}

This is again a difficult combinatorial problem, as the greedy heuristics of sorting objects with respect to the metric $f_i(x_i)/c_i(x_i)$ might not lead to an optimal assortment of objects.

\begin{figure}[htbp]
    \centering
    \begin{minipage}{0.85\linewidth}
        \centering
        
        \begin{subfigure}{\linewidth}
            \centering
            \begin{tikzpicture}
                \begin{axis}[
                    width=\linewidth,
                    height=5.2cm,
                    axis lines=left,
                    xlabel={$x_i$ (Allocated Quantity)},
                    ylabel={$f_i(x_i)$},
                    xmin=0, xmax=5,
                    ymin=0, ymax=130,
                    xtick={0, 1, 2, 3, 4, 5},
                    ytick={0, 40, 80, 120},
                    grid=both,
                    grid style={line width=.1pt, draw=gray!10},
                    major grid style={line width=.2pt, draw=gray!20},
                    every axis x label/.style={at={(current axis.right of origin)}, anchor=north east, yshift=-0.45cm},
                    every axis y label/.style={at={(current axis.above origin)}, anchor=south, yshift=0.2cm},
                    title={\textbf{Non-linear Profit with Capacity Threshold Discontinuity}},
                    title style={font=\small, yshift=0.1cm}
                ]
                    \addplot[
                        domain=0:2.98, 
                        samples=100, 
                        color=blue!70!black, 
                        ultra thick
                    ] {40*x - 3*x^2};
                    
                    \addplot[
                        domain=3:5, 
                        samples=100, 
                        color=blue!70!black, 
                        ultra thick
                    ] {40*x - 3*x^2 - 35};
                    
                    \addplot[dashed, color=gray!60, thin] coordinates {(3,0) (3,93)};
                    \draw[blue!70!black, fill=white, thick] (axis cs:3,93) circle (2pt);
                    \draw[blue!70!black, fill=blue!70!black] (axis cs:3,58) circle (2pt);
                    
                    \node[red!80!black, font=\footnotesize, anchor=south] at (axis cs:3, 95) {Truck Threshold Penalty};

                    \node[circle, fill=orange, inner sep=1.5pt, label={left:\footnotesize $f_i(2)=68$}] at (axis cs:2.0, 68) {};
                \end{axis}
            \end{tikzpicture}
        \end{subfigure}
        
        \vspace{0.6cm} 
        
        \begin{subfigure}{\linewidth}
            \centering
            \begin{tikzpicture}
                \begin{axis}[
                    width=\linewidth,
                    height=5.2cm,
                    axis lines=left,
                    xlabel={$x_i$ (Allocated Quantity)},
                    ylabel={$c_i(x_i)$},
                    xmin=0, xmax=5,
                    ymin=0, ymax=16,
                    xtick={0, 1, 2, 3, 4, 5},
                    grid=both,
                    grid style={line width=.1pt, draw=gray!10},
                    major grid style={line width=.2pt, draw=gray!20},
                    every axis x label/.style={at={(current axis.right of origin)}, anchor=north east, yshift=-0.45cm},
                    every axis y label/.style={at={(current axis.above origin)}, anchor=south, yshift=0.2cm},
                    title={\textbf{Resource Capacity Overhead Curve}},
                    title style={font=\small, yshift=0.1cm}
                ]
                    \addplot[
                        domain=0:5, 
                        samples=100, 
                        color=red!70!black, 
                        ultra thick
                    ] {1.5*x + 0.2*x^2};
                    
                    \node[circle, fill=orange, inner sep=1.5pt, label={left:\footnotesize $c_i(2)=3.8$}] at (axis cs:2.0, 3.8) {};
                \end{axis}
            \end{tikzpicture}
        \end{subfigure}
        
    \end{minipage}
    \caption{Visualizing a non-linear knapsack problem where the profit function incorporates a structural step-down penalty. When the allocated allocation quantity crosses the threshold ($x_i = 3$), the sudden acquisition cost of an additional delivery vehicle drops total net profit before marginal returns resume growth.}
    \label{fig:discontinuous_knapsack_profit}
\end{figure}
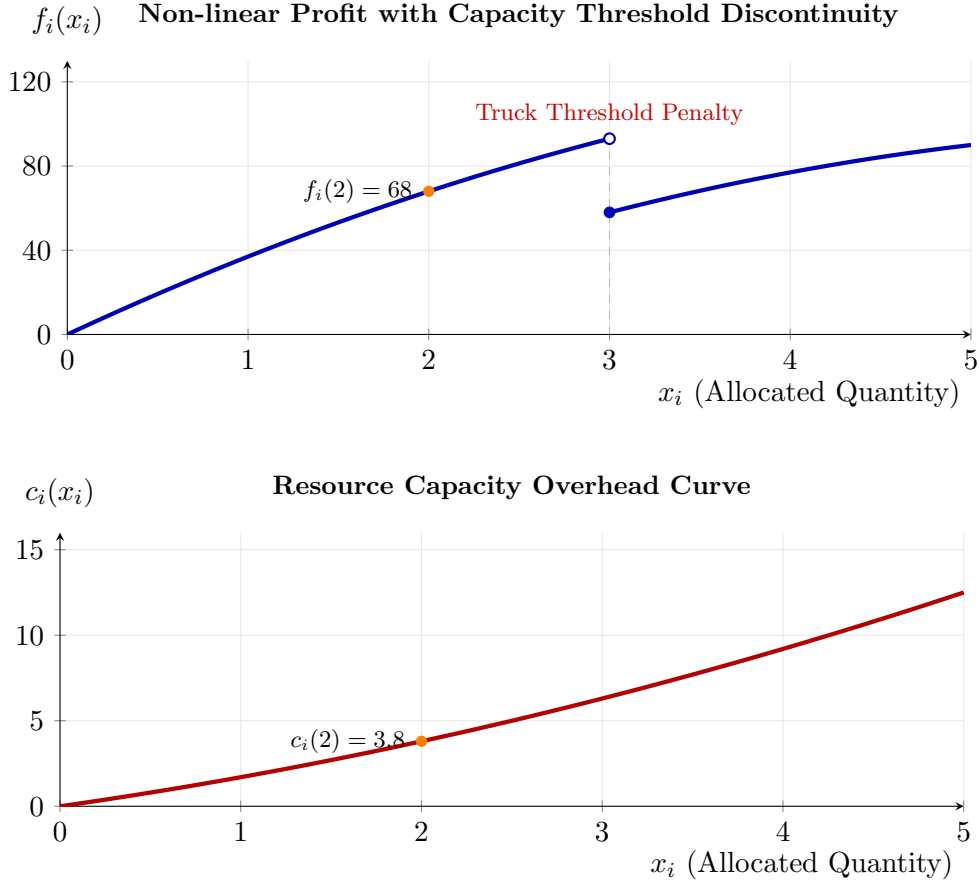

This problem can be solved in multiple ways: For example, the problem can be expressed as a Mixed-Integer Program and be solved directly with solvers. Although, when the benefit function $f_i(x_i)$ is non-linear this is not the case, see Figure \ref{fig:discontinuous_knapsack_profit} for an example.

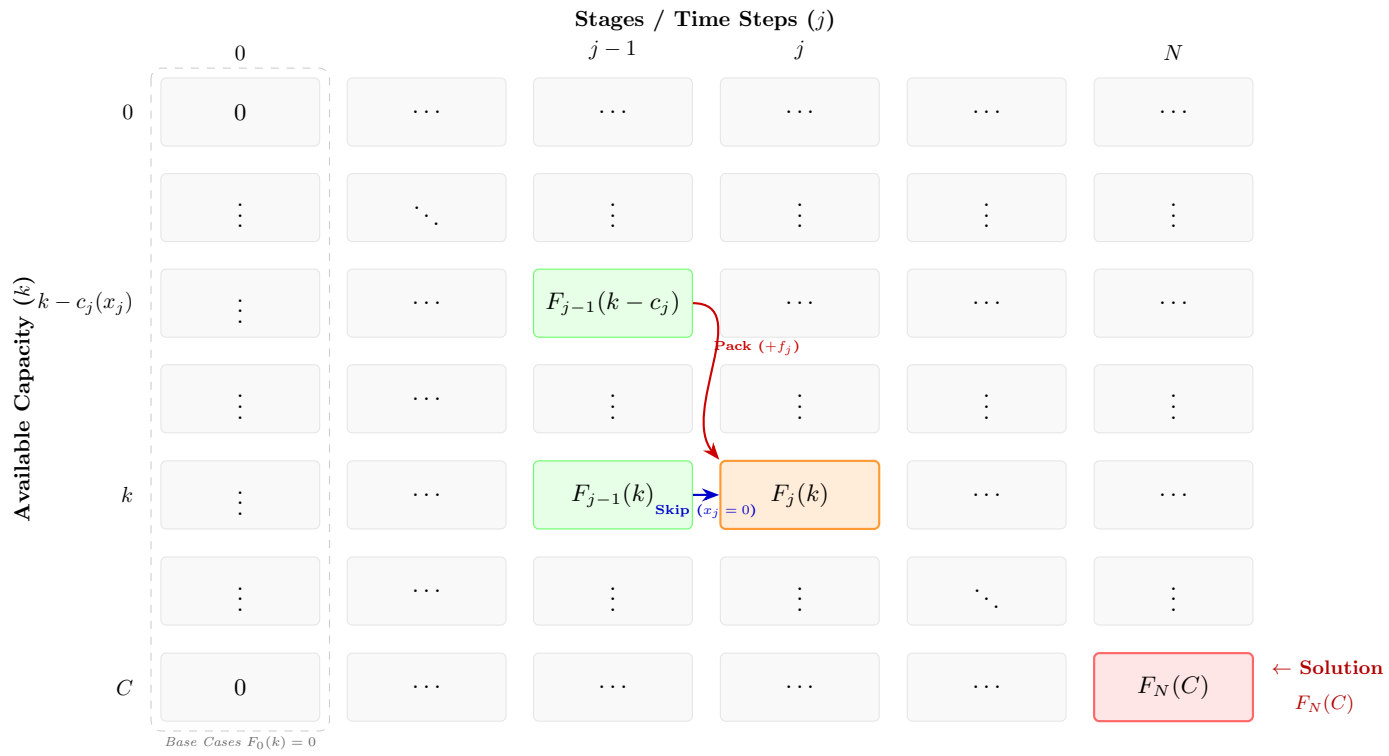
\begin{figure}[htbp]
    \centering
    \begin{minipage}{\linewidth}
        \centering
        \begin{tikzpicture}[
            >=Stealth,
            scale=0.82, transform shape 
        ]

            \matrix [
                matrix of math nodes,
                nodes in empty cells,
                row sep=0.35cm,
                column sep=0.35cm,
                nodes={draw=gray!20, fill=gray!5, minimum width=2.1cm, minimum height=0.9cm, font=\footnotesize, anchor=center, rounded corners=2pt}
            ] (M) {
                0      & \dots  & \dots  & \dots  & \dots  & \dots  \\
                \vdots & \ddots & \vdots & \vdots & \vdots & \vdots \\
                \vdots & \dots  & |[fill=green!10, draw=green!60, font=\bfseries\footnotesize]| F_{j-1}(k-c_j) & \dots  & \dots  & \dots  \\
                \vdots & \dots  & \vdots & \vdots & \vdots & \vdots \\
                \vdots & \dots  & |[fill=green!10, draw=green!60, font=\bfseries\footnotesize]| F_{j-1}(k)     & |[fill=orange!15, draw=orange!80, font=\bfseries\footnotesize, thick]| F_j(k) & \dots  & \dots  \\
                \vdots & \dots  & \vdots & \vdots & \ddots & \vdots \\
                0      & \dots  & \dots  & \dots  & \dots  & |[fill=red!10, draw=red!60, font=\bfseries\footnotesize, thick]| F_N(C) \\
            };

            \node[above=0.15cm of M-1-1, font=\small\bfseries] {$0$};
            \node[above=0.15cm of M-1-3, font=\small\bfseries] {$j-1$};
            \node[above=0.15cm of M-1-4, font=\small\bfseries] {$j$};
            \node[above=0.15cm of M-1-6, font=\small\bfseries] {$N$};
            \node[above=0.6cm of M-1-3.north east, font=\small\bfseries, xshift=0.2cm] {Stages / Time Steps ($j$)};

            \node[left=0.3cm of M-1-1, font=\small\bfseries] {$0$};
            \node[left=0.3cm of M-3-1, font=\small\bfseries] {$k-c_j(x_j)$};
            \node[left=0.3cm of M-5-1, font=\small\bfseries] {$k$};
            \node[left=0.3cm of M-7-1, font=\small\bfseries] {$C$};
            
            \node[left=2.2cm of M-4-1, font=\small\bfseries, rotate=90, anchor=center] {Available Capacity ($k$)};

            \draw[->, thick, blue!80!black] (M-5-3.east) -- 
                node[midway, below, font=\tiny\bfseries, text=blue!80!black] {Skip ($x_j=0$)} 
                (M-5-4.west);

            \draw[->, thick, red!80!black] (M-3-3.east) to[out=0, in=130] 
                node[midway, above right, font=\tiny\bfseries, text=red!80!black, xshift=-0.15cm] {Pack ($+f_j$)} 
                (M-5-4.north west);

            \draw[draw=gray!40, dashed, rounded corners] ([xshift=-0.15cm, yshift=0.15cm]M-1-1.north west) rectangle ([xshift=0.15cm, yshift=-0.15cm]M-7-1.south east);
            \node[below=0.1cm of M-7-1, font=\tiny\itshape, gray!80!black] {Base Cases $F_0(k)=0$};
            
            \node[right=0.15cm of M-7-6, font=\footnotesize\bfseries, text=red!70!black, align=left] {$\leftarrow$ Solution\\ \ \ \ $F_N(C)$};

        \end{tikzpicture}
    \end{minipage}
    \caption{Matrix representation of the forward-recursion Dynamic Programming state space for the Knapsack problem. Columns represent sequential decision horizons (stages/time), while rows represent the remaining sub-problem capacity. Computing the cell $F_j(k)$ evaluates the maximum option between state transitions in stage $j-1$.}
    \label{fig:knapsack_dp_matrix}
\end{figure}

For a general benefit function, we present an alternative procedure to achieve the same result. This procedure is called \emph{forward-recursion} or \emph{Dynamic Programming} depending on the textbook. In fact, this is similar to the production example or the dynamic inventory model presented in Chapter \ref{sec:inventory_management}). First define the optimal benefit for a knapsack of capacity $k\le C$ using the first $j\le N$ objects:
\begin{eqnarray}
  F_j(k):=\max_{\sum_{i=1}^{j}c_i(x_i)\le k,x_i\in \calx}\left(\sum_{i=1}^{j}f_i(x_i)\right),
\end{eqnarray}
The solution of the problem is by definition $F_N(C)$. Now, the trick is again to create a one-step recursion that relates these quantities to their predecessors (as we are going forwards instead of backwards). This recursion is:
\begin{eqnarray}
  F_j(k)=\max_{x_j\in\calx}\left\{f_j(x_j)+F_{j-1}(k-c_j(x_j))\right\},
\end{eqnarray}
The intuition is clear: the maximum benefit of a knapsack of capacity $k$ using the first $j$ objects is equal to the largest benefit of the object $j$, that is, $f_j(x_j)$ plus the optimal benefit of a knapsack using the previous $j-1$ objects with capacity $k-c_j(x_j)$. The only thing left is to define the boundary conditions to start the recursion. These can be defined as:
\begin{eqnarray}
 	F_0(k)=\begin{cases}
 		0, & \text{for } k\ge 0 \\
    -\infty, & \text{for } k<0.
 	\end{cases}
\end{eqnarray}
Then, the solution can be computed recursively by filling out the elements $F_j(k)$ in a forward fashion. The underlying assumption is that the capacity function $c_i$ is a non-decreasing function of $x_i$. And the number of operations is of order $O(CN|\calx|)$. See Figure \ref{fig:knapsack_dp_matrix} for a schematic of the DP recursion.

 \chapter{Vehicle Routing Problems} 
\label{sec:vehicle_routing}
\section{Introduction}
In previous chapters we have described matching the supply and demand in abstract sets, either in the time or space domain. In this chapter we combine both into a single optimization paradigm the optimization of both at the same time. The most common instance of this problem is optimizing transportation routes for vehicles. In a practical sense, this is a key problem to match supply and demand as the feasibility of the supply chain operation is not a given. For example, how many vehicles are needed to transport X amount of material within a Y timeframe is far from evident. Computationally, these belong to a family of problems that are hard to solve. Similar to the set covering problems studied in the previous chapter, these are combinatorial problems where time and space have to be allocated simultaneously.

While solving these problems to optimality is very hard, we will study some heuristic solutions that will allow us to get feasible solutions in computationally tractable time. These algorithms are a combination of the optimization algorithms that we have studied for warehouse allocation, minimizing the traveling time while fulfilling the demand.

\section{The Traveling Salesman Problem}

The cornerstone of transportation algorithms is known as the Traveling Salesman Problem (TSP), and it describes the following situation: imagine a salesman that has a list of locations that they have to visit (as a side-note, the forefather of Operations Research in the U.S., Phillip Morse, described himself primarily as a ``salesman'', so in a way, everyone involved with Operations Research or Operations Management is in fact a kind of salesman, even unknowingly). This salesman wants to visit the locations in an order that minimizes their cost in some way (time, effort, transportation costs, etc).

Intuitively, the reader can already guess how a formulation would work: the objective function is to minimize a cost that will be associated to the specific \emph{tour} that the salesman will take, subject to constraints on the validity of the tour. The challenge is then to express everything in a single formulation.

Similar to previous chapters, we start by defining locations as vertices in a graph $G=V,E$, with a collection of vertices $V=\{1,\dots,N\}$ and edges $E\subseteq \{(i,j):i,j\in V\}$. As a first primitive we consider the cost of traveling between nodes $i$ and $j$ as $c_{ij}$. For example, this cost can follow a similar structure to the warehouse location problem where it was dependent on a random variable $T_{i,j}$ and a constant. Traversing the graph, can be uderstood as a tour, a tour $\tau=(i_1,\dots,i_N)$ denotes a sequence of $N$ vertices that will be visited in that specific order: starting at $i_1\in V$, then going to $i_2\in V$, and so on until $i_N\in V$ is reached where the tour ends. Then, the objective in the TSP is to find a tour that minimizes the traveling cost.

A natural way to describe whether an edge belongs to the minimum cost tour in this graph is to define a decision variable $x_{i,j}=1$ if edge $(i,j)$ is part of the optimal tour and 0 otherwise. Then, the only challenge left is to express the constraints that define what a valid tour is.

The first property is that from every edge there is only one entering edge and one outgoing edge. Mathematically, this can be written as the constraints $\sum_{i\neq j}x_{i,j}=1$ for all $j\in V$ and $\sum_{j\neq i}x_{i,j}=1$ for all nodes $i\in V$. To better visualize this, let the solution be stored in a matrix $X=[x_{i,j}]$ where the element $i,j$ is $x_{i,j}$. For example, consider the points in Figure \ref{fig:tsp_irregular_example} representing the optimal tour for 4 arbitrary points in a plane, and their corresponding solution $X$.

\begin{figure}[htbp]
    \centering
    \begin{subfigure}[b]{0.45\textwidth}
        \centering
        \begin{tikzpicture}[scale=1.5, every node/.style={circle, draw, inner sep=2pt}]
            \node (1) at (0, 2) {1};     
            \node (2) at (0.5, 0) {2};   
            \node (3) at (2.5, 0.5) {3}; 
            \node (4) at (2, 2.3) {4};   

            \draw[thick, ->, >=stealth] (1) -- (4);
            \draw[thick, ->, >=stealth] (4) -- (3);
            \draw[thick, ->, >=stealth] (3) -- (2);
            \draw[thick, ->, >=stealth] (2) -- (1);
        \end{tikzpicture}
        \caption{Optimal tour $\tau = (1, 4, 3, 2)$ on an irregular polygon.}
    \end{subfigure}
    \hfill
    \begin{subfigure}[b]{0.45\textwidth}
        \centering
        \[
        X = \begin{pmatrix}
            0 & 0 & 0 & 1 \\
            1 & 0 & 0 & 0 \\
            0 & {1} & 0 & 0 \\
            0 & 0 & {1} & 0
        \end{pmatrix}
        \]
        \caption{Corresponding adjacency matrix $x_{i,j}$}
    \end{subfigure}
    \caption{Visualizing the relationship between a physical TSP tour on an irregular polygon and its mathematical representation in the decision matrix $X$.}
    \label{fig:tsp_irregular_example}
\end{figure}

Enforcing the constraints that there is only one incoming and outgoing edge from every tour is equivalent to making the matrix $X$ row and column sum to be equal to $\mathbf{1}$. That is, $\mathbf{1}X=\mathbf{1}$ and $X\mathbf{1}^\intercal=\mathbf{1}^\intercal$. With these, an initial formulation of the TSP minimizing the expected travel cost is (letting $c_{i,i}=\infty$ to avoid solutions where a point goes to itself):

\begin{eqnarray}
& \min& \ex\left\{\sum_{i=1}^N\sum_{j=1}^Nx_{i,j}c_{i,j}\right\}\\
&\text{s.t.}&X=[x_{i,j}], x_{i,j}\in\{0,1\},\nonumber\\
    &&\mathbf{1}X=\mathbf{1},\,\,\, X\mathbf{1}^\intercal=\mathbf{1}^\intercal.\nonumber
\end{eqnarray}

When solving this problem as is there are two possible outcomes: the first, is that the problem results in a valid tour that will be optimal. The second, is that these constraints might not be enough to enforce a tour that visits all points and the result might be a collection of subtours that do not visit all points, see Figure \ref{fig:tsp_subtours} for an example of this situation.

\begin{figure}[htbp]
    \centering
    \begin{subfigure}[b]{0.45\textwidth}
        \centering
        \begin{tikzpicture}[scale=1.2, every node/.style={circle, draw, inner sep=2pt}]
            \node (1) at (0.5, 3) {1};
            \node (3) at (0, 1.5) {3};
            \node (4) at (1.5, 2) {4};
            
            \node (2) at (3, 2.5) {2};
            \node (5) at (2.5, 0.5) {5};
            \node (6) at (4, 1) {6};
            
            \draw[thick, ->, >=stealth] (1) -- (3);
            \draw[thick, ->, >=stealth] (3) -- (4);
            \draw[thick, ->, >=stealth] (4) -- (1);
            
            \draw[thick, ->, >=stealth] (2) -- (5);
            \draw[thick, ->, >=stealth] (5) -- (6);
            \draw[thick, ->, >=stealth] (6) -- (2);
        \end{tikzpicture}
        \caption{Disjoint Subtours: $\{1,3,4\}$ and $\{2,5,6\}$}
    \end{subfigure}
    \hfill
    \begin{subfigure}[b]{0.45\textwidth}
        \centering
        \[
        X = \begin{pmatrix}
            0 & 0 & 1 & 0 & 0 & 0 \\
            0 & 0 & 0 & 0 & 1 & 0 \\
            0 & 0 & 0 & 1 & 0 & 0 \\
            1 & 0 & 0 & 0 & 0 & 0 \\
            0 & 0 & 0 & 0 & 0 & 1 \\
            0 & 1 & 0 & 0 & 0 & 0
        \end{pmatrix}
        \]
        \caption{Matrix $X$ satisfying assignment constraints}
    \end{subfigure}
    \caption{An example where the assignment constraints $\mathbf{1}X=\mathbf{1}$ and $X\mathbf{1}^\intercal=\mathbf{1}^\intercal$ are satisfied, yet the solution consists of two disconnected subtours.}
    \label{fig:tsp_subtours}
\end{figure}

Clearly, the solution satiesfies the constraints we originally imposed, but it ended up producing two subtours instead of 1 tour that visits all the nodes. Here, the most practical approach is to impose a constraint that after seeing this solution, imposes that the two subtours must communicate with each other. In this example, let $v_s$ be an indicator vector equal to 1 if node i belongs to subtour $s$. For example, for the subtour composed by the nodes $s=\{1,3,4\}$, $v_s=(1,0,1,1,0,0)$. Likewise, $1-v_s$ is an indicator of the nodes not belonging to $s$, that is $V\setminus s$. Enforcing that the subtour $s$ is connected outside, amounts to the statement: make that at least one of the outgoing connections in $s$ goes outside it. The vector $v_s X$ represents all the connections starting at $s$ and $X(\mathbf{1}-v_s)^\intercal$ means all the connections going to elements outside of $s$. Then, we can put both conditions together into a single expression and force that there is at least one such connection from $s$ to outside it, that is, $v_sX(\mathbf{1}-v_s)^\intercal\ge 1$. If we add this constraint to our original formulation we end up with:

\begin{eqnarray}
& \min& \ex\left\{\sum_{i=1}^N\sum_{j=1}^Nx_{i,j}c_{i,j}\right\}\\
&\text{s.t.}&X=[x_{i,j}], x_{i,j}\in\{0,1\},\nonumber\\
    &&\mathbf{1}X=\mathbf{1},\,\,\, X\mathbf{1}^\intercal=\mathbf{1}^\intercal,\nonumber\\
    &&v_sX(\mathbf{1}-v_s)^\intercal\ge 1.\nonumber
\end{eqnarray}

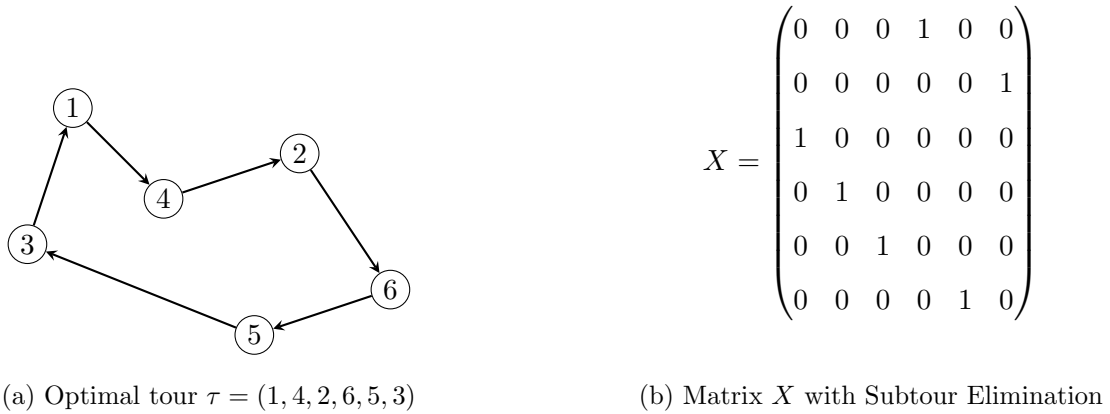
\begin{figure}[htbp]
    \centering
    \begin{subfigure}[b]{0.45\textwidth}
        \centering
        \begin{tikzpicture}[scale=1.2, every node/.style={circle, draw, inner sep=2pt}]
            \node (1) at (0.5, 3) {1};
            \node (3) at (0, 1.5) {3};
            \node (4) at (1.5, 2) {4};
            \node (2) at (3, 2.5) {2};
            \node (5) at (2.5, 0.5) {5};
            \node (6) at (4, 1) {6};
            
            \draw[thick, ->, >=stealth] (1) -- (4);
            \draw[thick, ->, >=stealth] (4) -- (2);
            \draw[thick, ->, >=stealth] (2) -- (6);
            \draw[thick, ->, >=stealth] (6) -- (5);
            \draw[thick, ->, >=stealth] (5) -- (3);
            \draw[thick, ->, >=stealth] (3) -- (1);
        \end{tikzpicture}
        \caption{Optimal tour $\tau = (1, 4, 2, 6, 5, 3)$}
    \end{subfigure}
    \hfill
    \begin{subfigure}[b]{0.45\textwidth}
        \centering
        \[
        X = \begin{pmatrix}
            0 & 0 & 0 & 1 & 0 & 0 \\
            0 & 0 & 0 & 0 & 0 & 1 \\
            1 & 0 & 0 & 0 & 0 & 0 \\
            0 & 1 & 0 & 0 & 0 & 0 \\
            0 & 0 & 1 & 0 & 0 & 0 \\
            0 & 0 & 0 & 0 & 1 & 0 
        \end{pmatrix}
        \]
        \caption{Matrix $X$ with Subtour Elimination}
    \end{subfigure}
    \caption{The optimal solution after imposing $v_s X(\mathbf{1}-v_s)^\intercal \ge 1$. The matrix now represents a single connected cycle visiting all six nodes.}
    \label{fig:tsp_optimal_full}
\end{figure}

After solving this formulation, we finally end up with a solution that not only yields a valid tour, but it is the optimal one. See Figure \ref{fig:tsp_optimal_full} for the same instance of the previous problem with subtours. The solution of adding all subtours as constraints is known as the Dantzig-Fulkerson-Johnson formulation of the TSP (see \cite{dantzig1954solution}). Their algorithm would amount to add constraints for all possible subtours $s\subseteq V$. Needless to say, this is an exponentially large set of constraints that in real-life is rarely implemented as is. Instead, modern TSP solvers have a feature called \emph{lazy constraints} where the solver sequentially adds the subtour elimination constraints as they emerge in the optimization process.

If there is not access to a solver that can add these constraints automatically, the algorithmic procedure is just to solve the problem sequentially and add constraints if the optimal solution is not a valid tour that visits all the vertices. In practice, this procedure converges relatively quickly. In Algorithm \ref{alg:tsp_compact} we summarize the procedure.

\begin{algorithm}
    \caption{Subtour Elimination Algorithm}
    \label{alg:tsp_compact}
    \begin{algorithmic}[1]
        \State \textbf{Input:} Costs $C$, nodes $V$. \textbf{Output:} Optimal tour $X^*$.
        \Procedure{TSP-Iterative}{$V, C$}
        \State $\mathcal{S},\mathcal{K} \gets \emptyset$
        \While{$|\mathcal{K}| \neq 1$}
            \State $X^* \gets \text{argmin} \left\{ \sum c_{i,j}x_{i,j} : \mathbf{1}X=\mathbf{1}, X\mathbf{1}^\intercal=\mathbf{1}^\intercal, \mathcal{S} \right\}$
            \State $\mathcal{K} \gets \text{DisjointCycles}(X^*)$
            \State $\mathcal{S} \gets \mathcal{S} \cup \{v_s X (\mathbf{1}-v_s)^\intercal \ge 1 \mid s \in \mathcal{K}\}$
        \EndWhile
        \State \textbf{return} $X^*$
        \EndProcedure
    \end{algorithmic}
\end{algorithm}

The difference between this procedure and what the solvers do is that the solver can incorporate the new constraints while keeping the information of the branch and bound tree (in a very simplified summary, MIP solvers solve a sequence of linear programs by sequentially adding constraints to a ``tree'' of candidate solutions until an LP with an integral solution is reached).

\subsection{Incorporating Randomness and Risk Awareness}
We can incorporate the techniques that we learned previously in the warehouse location optimization Chapter \ref{sec:networks} to incorporate randomness and create routes that are robust to variance in the cost. A classical example is the case when the cost $c_{i,j}$ is the random transportation time $T_{i,j}$ between node $i$ and node $j$. This is an important case because as we have seen, points that have in expected value a smaller traveling time could potentially have a high variance which could make some travels experience longer travel times. Which depending on the risk aversion of the decision-maker could be possible unnaceptably long. For example, consider a medical transportation problem where a patient needs to be picked up at their home and transported to a hospital. The goal is to bring the patient as fast as possible while hedging the risk of a potentially fatal longer travel. In general supply chains, time is always a tight constraint, where products have to be delivered under strict delivery windows.

The first of these models is the same technique that we studied to solve a mean-variance version of the problem. The variance of the sum $\sum_{i=1}^N\sum_{j=1}^Nx_{i,j}c_{i,j}$ is equal to the same expression for the variance of a sum of random variables that we studied previously in Chapter \ref{sec:networks} $\sum_{i,j\in V^2}\sum_{k,l\in V^2}x_{i,j}\cov(c_{i,j},c_{k,l})x_{k,l}$, or more succintly, $x\cov(c)x^\intercal$ with $x$ as the flattened vector of $X$ (and $c$ similarly defined). With these, the optimization problem can be stated as:

\begin{eqnarray}
& \min& \ex\left\{\sum_{i=1}^N\sum_{j=1}^Nx_{i,j}c_{i,j}\right\}+\lambda \var\left\{\sum_{i=1}^N\sum_{j=1}^Nx_{i,j}c_{i,j}\right\}\\
&\text{s.t.}&X=[x_{i,j}], x_{i,j}\in\{0,1\},\nonumber\\
    &&\mathbf{1}X=\mathbf{1},\,\,\, X\mathbf{1}^\intercal=\mathbf{1}^\intercal,\nonumber\\
    &&v_sX(\mathbf{1}-v_s)^\intercal\ge 1.\nonumber
\end{eqnarray}

Where the subtour constraints can be added on the fly as in algorithm \ref{alg:tsp_compact} if they are random.

\section{The Vehicle Routing Problem}

We have explored and defined the problem with only one salesman, representing a single resource like a vehicle, that must visit all $N$ nodes. In this subsection we extend the formulation to optimize $K$ vehicle routes. The simplest way to extend our original formulation is to add a special node, numbered $0$, representing a warehouse or depot where all vehicles depart to make their routes and return once they have finished their assigned routes. 

Given that each location still will be visited only once, we can try to retain most of the structure of the problem that we have already studied. The only node that is going to behave differently, is the depot node $0$, as all vehicles will depart and return to the same location. The number of vehicles departing and arriving to node $0$ is exactly $K$. That is, $\sum_{j\neq 0}x_{0,j}=K$ for the outgoing vehicles, and $\sum_{j\neq0}x_{j,0}=K$ for the entering vehicles. Note that the number of nodes is now $N+1$ instead of $N$. Similar to the previous section, this leads to a basic and incomplete formulation, missing the equivalent of the subtour elimination constraints:

\begin{eqnarray}
& \min& \ex\left\{\sum_{i=0}^N\sum_{j=0}^Nx_{i,j}c_{i,j}\right\}\\
&\text{s.t.}&X=[x_{i,j}], x_{i,j}\in\{0,1\},\nonumber\\
    &&(K,\mathbf{1})X=(K,\mathbf{1}),\nonumber\\
    &&X(K,\mathbf{1})^\intercal=(K,\mathbf{1})^\intercal.\nonumber
\end{eqnarray}

Similar to the previous section, solving this problem directly might result in one of two outcomes: The first is that the optimization solves the problem and there is a valid assignment of $K$ routes where all vehicles have non-overlapping routes starting and ending at node $0$. The other, is that there might be isolated islands of customers that are not served from a truck that departed from the warehouse. These \emph{ghost} islands that are not served from a truck, have again to be connected to a truck that visits the depot. For example, see Figure \ref{fig:vrp_without_subtour} for an example of this phenomenon. In this example there are two trucks that visit customes 1 to 4, but customers 5 and 6 are not visited.

\begin{figure}[htbp]
    \centering
    \begin{subfigure}[b]{0.45\textwidth}
        \centering
        \begin{tikzpicture}[scale=1.2, every node/.style={circle, draw, inner sep=2pt}]
            \node (0) at (0, 0) {\textbf{0}}; 
            
            \node (1) at (-1, 1) {1};
            \node (2) at (1, 1) {2};
            \draw[thick, ->, >=stealth] (0) -- (1);
            \draw[thick, ->, >=stealth] (1) -- (2);
            \draw[thick, ->, >=stealth] (2) -- (0);
            
            \node (3) at (-1, -1) {3};
            \node (4) at (1, -1) {4};
            \draw[thick, ->, >=stealth] (0) -- (3);
            \draw[thick, ->, >=stealth] (3) -- (4);
            \draw[thick, ->, >=stealth] (4) -- (0);
            
            \node (5) at (3, 0.5) {5};
            \node (6) at (3, -0.5) {6};
            \draw[thick, ->, >=stealth] (5) -- (6);
            \draw[thick, ->, >=stealth] (6) -- (5);
        \end{tikzpicture}
        \caption{Two valid routes and one isolated island $\{5,6\}$.}
    \end{subfigure}
    \hfill
    \begin{subfigure}[b]{0.45\textwidth}
        \centering
        \[
        X = \begin{pmatrix}
            0 & 1 & 0 & 1 & 0 & 0 & 0 \\
            0 & 0 & 1 & 0 & 0 & 0 & 0 \\
            1 & 0 & 0 & 0 & 0 & 0 & 0 \\
            0 & 0 & 0 & 0 & 1 & 0 & 0 \\
            1 & 0 & 0 & 0 & 0 & 0 & 0 \\
            0 & 0 & 0 & 0 & 0 & 0 & 1 \\
            0 & 0 & 0 & 0 & 0 & 1 & 0
        \end{pmatrix}
        \]
        \caption{Matrix $X$ with $x_{i,j} \in \{0,1\}$.}
    \end{subfigure}
    \caption{The row/column sum for the depot is 2, but customers 5 and 6 are never reached by a truck starting at 0.}
    \label{fig:vrp_without_subtour}
\end{figure}

Similar to the previous section, we can force the solver to eliminate the customer's only islands by adding subtour elimination constraints. In this case, the island of customers $s=\{5,6\}$ is represented by the the vector $v_s=(0,0,0,0,0,1,1)$. Then, adding the constraint $v_sX(\mathbf{1}-v_s)^\intercal\ge 1$ to the original problem will eliminate that specific isolated island. Then, after solving the problem with the new constraint might again have new islands, or yield a valid solution that will be optimal.

We can apply the same strategy as in the TSP, to solve a formulation adding a subtour elimination constraint to eliminate the islands. The algorithmic procedure to solve the problem would be identical to the TSP, with the only change of changing the DisjointCycles subroutine, to an IsolatedCycles one that instead detects the ``customers only" islands generated by the solution in the current step. The master problem is remarkably similar to the TSP one:

\begin{eqnarray}
& \min& \ex\left\{\sum_{i=0}^N\sum_{j=0}^Nx_{i,j}c_{i,j}\right\}\\
&\text{s.t.}&X=[x_{i,j}], x_{i,j}\in\{0,1\},\nonumber\\
    &&(K,\mathbf{1})X=(K,\mathbf{1}),\nonumber\\
    &&X(K,\mathbf{1})^\intercal=(K,\mathbf{1})^\intercal,\nonumber\\
        &&v_sX(\mathbf{1}-v_s)^\intercal\ge 1. \nonumber
\end{eqnarray}
As a final note, similar techniques for risk hedging can be applied to the same problem where the shape of the covariance is remarkably similar.

\subsection{On Capacities and Similar Issues}
So far, we have studied algorithms that minimize a cost measure $c_{i,j}$ that is independent of the routes and the vehicles. An avid reader might be rightfully asking questions about different operational issues that happen with vehicles such as capacities, time-window arrivals and other operational constraints. We present two schools of thought regarding these issues and present our biased assessment of the right modeling approach.

We call the first school of thought the school of the \emph{masochists}. This modeling school wants to include in a single formulation all the constraints of the problem and through indicators, counter variables and the big-M method ensure that the solution of the problem (however long it might take, thus the masochism part) returns a feasible and optimal solution. Except in the smallest instances, this approach tends to fail as the linear programming relaxation of the resulting problem is \emph{weak}. In the sense that the solver starts very far away from a set of cuts that can even yield a feasible solution. An example of masochism is to try to formulate the staff scheduling problem studied in Chapter \ref{sec:scheduling} with modeling every single worker (out of $H$) with an indicator variable if worker $i$ is active and hour $t$ and try to brute-force the formulation by writing a large number of constraints ensuring that the shift is no longer than 8 hours and so on. Masochists are often limited by their own lifespan in their endeavor of trying to force structure by enumeration and straitjackets. That being said, for small problem instances one can get away with a formulation that gets the job done in a reasonable time.

The second school of thought corresponds naturally to the \emph{sadists}. The sadists start with a very lean formulation that often yields unfeasible solutions, and through a process of hammering down constraints, eventually a feasible solution is reached. The exposition in this section corresponds to some degree of sadism in the formulation taste. The sadism comes from the fact that the number of interations needed for the process to converge to a feasible and optimal solution is uncertain. While there is some empirical evidence that when starting with a \emph{strong} formulation, the number of iterations to reach optimality, that is, of adding constraints (out of exponentially many) consists of an unknown number of iterations, there is a proof that this will always be the case. Moreover, the ultimate bottleneck is the inherent comobinatorial nature of the problem.

The latter paradigm leads to an interesting strategy to deal with certain common patterns in this type of problems. The first one is capacities. Normally, most vehicles are constrained by a physical limit on their capacity, for imagine the case where all vehicles have a capacity $C$. Going back to our usual exposition, imagine a vector of demands $\bm{D}=(D_1,\dots,D_N)$ that needs to be satisfied.  After solving the VRP the resulting tour for each vehicle $k=1,\dots,K$ is $\tau_k$, or the set of nodes $s_k$. Then, it could be the case that the sum of all the demands exceeds the capacity of some vehicle(s) $k$. That is, $\sum_{i\in s_k}\ex[D_i]\ge C$ for some vehicle(s) $k$, in which case the solution is unfeasible. This kind of unfeasibility can be again tackled by the same technique to eliminate the subtours that violated the feasibility of the TSP or VRP. The key observation is that in the subset of nodes $s_k$ in the tour $\tau_k$ there are not enough vehicles to service their cumulative required demand. The minimum number of vehicles that satisfies the demand is $m(s_k):=\lceil\sum_{i\in s_k}\ex [D_i]/C\rceil$. Then, it must be the case that for this specific group of nodes, the number of outgoing vehicles is at least $m(s_k)$. The logic is that this enforces that at least this amount of vehicles enters the specific subset of nodes. Then, the constraint that enforces this condition is similarly $v_{s_k}X(1-v_{s_k})^\intercal\ge m(s_k)$ for all vehicles $k$.

The exact same algorithm for the TSP and VRP can be applied by adding \emph{lazy constraints} with an additional subroutine called \emph{CapacityUnfeasibility} checking the subsets of customers that violate the capacity constraints and adding the constraints $v_{s_k}X(1-v_{s_k})^\intercal\ge m(s_k)$ accordingly, until feasibility is reached. On the capacities, they can be calculated as a number that hedges the risk that the demand exceeds certain level as in the Newsvendor model in Chapter \ref{sec:inventory_management}.

\subsection{Time-window constraints: Column generation}

A common business requirement in transportation problems is that deliveries have to be performed at specific hours for certain clients. For example, transporting perishables or food products often requires deliveries to occur during specific time slots. This is a challenging constraint to incorporate in our formulation as the temporal dimension is somewhat abstracted in the decision variables and there is not a direct mechanism built-in to keep track of the times the deliveries occur. 

One might think to use the same approach of adding constraints to the main problem when the deliveries do not occur at the specified time-window. This approach is not the strongest as the constraints just force to shuffle the order of the tours, which amounts to a slow and \emph{thin} cut in the geometry of the solution. In the capacity case, the cuts (of the feasible set) are stronger as they force a minimum number of vehicles within a certain subset, which in subsequent runs of the optimization problem greatly reduces the space of candidate solutions.

To make an interesting parallel, in the staff scheduling model in Chapter \ref{sec:scheduling} the decision variables were the number of resources (people) assigned to each of the possible schedules. This was possible because the possible schedules were a relatively small set due to their structure. The resulting formulation had a clean and direct interpretation. As we discussed in the previous subsection, our modeling approach for the routes is equivalent to create indicator variables for each of the workers and whether or not they are active during a certain hour. In hindsight, this seems like a silly modeling choice when the cleaner formulation is known.

If we wanted to formulate the TSP with time-window constraints using the same approach as the staff scheduling formulation we would set a problem to choose from a universe of routes that already satisfies the time-window constraints, and then we would let the formulation choose the one that minimizes the cost function. The formulation would look something like this:

\begin{eqnarray}
  &\min& \ex\left[\sum_{\tau\in\Omega}c_\tau y_\tau\right]\\
  &\text{s.t.}& \sum_{\tau\in\Omega}a_{i\tau}y_\tau= 1,\,\,\text{ for all }\,\, i\in V,\nonumber\\
  &&y_\tau\in\{0,1\}.\nonumber
\end{eqnarray}
$\Omega$ represents all the tours where locations are visited on time (each location $i\in V$ has a time window $[s_i,t_i]$ for arrival). $a_{i\tau}$ is the indicator that tour $\tau$ visits location $i$. The first constraint represents that each location must be visited exactly once. Since all tours $\tau$ in the set $\Omega$ are feasible, that is, they visit the locations in their corresponding arrival windows. This problem is ensured to find the optimal solution of the problem.

The issue is that in this case it is difficult to know beforehand the set of tours that are feasible (in the staff scheduling case, the columns of the matrix $A$ were evident from the context of the problem). The goal is then to generate \emph{columns} in this problem so that a feasible tour that is optimal is found. Moreover, there are exponentially many of these columns, which makes enumeration difficult.

Suppose that we are able to generate some routes that are feasible. The question is then, how can we define a procedure to create new routes to eventually find an optimal one. Here, we are going to use a technique that we used already in Chapter \ref{sec:scheduling} to solve a dynamic programming problem to generate new feasible (and better routes).

To build such a procedure, we need a way to peek into the structure of the problem to get a numerical estimate on how to modify a route and then apply a similar logic to the Knapsack iteration to build feasible routes that have lower costs. Duality offers an interesting insight to build such estimate: Imagine that we solve the linear problem relaxation of the formulation (with the constraints $0\le y_\tau\le 1$ instead of the binary ones). By duality, each of the constraints $\sum_{\tau\in\Omega}a_{i\tau}y_\tau= 1$ has a corresponding dual variable $\pi_i$, the so-called \emph{shadow price} of the constraint. In this case, the variable represents the increase in cost of adding location $i$ to a tour. When this increase is negative, it means that routes passing through that location reduce the cost. Imagine for a route $\tau$ with cost $c_\tau$ the corresponding dual prices in the relaxation are $\sum_{i\in V}a_{i\tau}\pi_i$. Then, if $c_\tau<\sum_{i\in V}a_{i\tau}\pi_i$ adding this route will result in a reduced cost, because in the current solution we are \emph{overpaying} to visit all the customers. Define the reduced cost $\bar c_\tau$ of a route as:

\begin{eqnarray}
\bar c_\tau := c_\tau-\sum_{i\in V}a_{i\tau}\pi_i,
\end{eqnarray}
Routes $\tau$ where $\bar c_\tau$ is negative are good condidates to include in the main optimization problem. When no such routes exist, we know that the optimal route is within the set of routes that we have considered. Then, we have found an optimal route.

Next, we discuss how to generate feasible routes. As discussed before, this is hard on face value, because just generating a route is a difficult combinatorial problem in itself. The structure of the problem is in fact similar to the \emph{knapsack} problem studied in the previous chapter. We define $F_j(k)$ as the minimum reduced cost of reaching node $j$ in a valid tour arriving at time $k$. As before, the key is to establish this recursion in a one-step way as:

\begin{eqnarray}
  F_j(k)=\min_{i\neq j\in V}\{(c_{i,j}-\pi_i)+F_i(k-t_{ij})\},
\end{eqnarray}
In principle, solving the Dynamic Program defined by this recursion should result in tours $c_\tau=F_0(k)<0$ with negative reduced cost that in principle are candidates to be added to the master program. The challenge is ensuring that at every step of the DP, the resulting tours are cycle free and feasible. The first step in ensuring feasibility is that the tours respect the delivery windows. This can be easily enforced as $F_j(k)=\infty$ if $k\notin [s_j,t_j]$, that is, $k$ is not in the delivery window.

We initialize the DP with the boundary conditions $F_0(0)=0$ and $F_k(j)=\infty$ for all other $j,k$. In other words, the tours start at the warehouse at time 0 and all the other entries will be updated as the DP recursion continues.

To populate the entries of the DP, we iterate on time from $k=0,\dots,T_{max}$ and through the nodes $j=1,\dots,N$. What can happen, is that for a given node $j$, we might have visited it before (at a time earlier than $k$) and also with lower reduced cost. Therefore it makes sense that this is not a better tour than one we have found already. This can be stated in the updated DP recursion:

\begin{eqnarray}
  F_j(k)=\begin{cases}
    \min_{i\neq j\in V}\{(c_{i,j}-\pi_i)+F_i(k-t_{ij})\}&\text{if }F_j(k)<\min_{t<k}F_j(t),k\in [s_j,t_j].\\
    \infty&\text{otherwise}
  \end{cases}
\end{eqnarray}
This updated recursion ensures that the resulting tours are prunned and non-improving tours are not considered. We summarize the full algorithm in \ref{alg:vrptw_column_generation}. In principle, once these columns are generated, we can combine the procedure with the capacity variation to obtain capacitated and time-window constraints.

\begin{algorithm}
    \caption{Column Generation for VRPTW}
    \label{alg:vrptw_column_generation}
    \begin{algorithmic}[1]
        \State \textbf{Input:} Costs $C$, nodes $V$, time windows $[s_i, t_i]$, depot $0$.
        \State \textbf{Output:} Optimal set of routes $y^*$.
        \Procedure{VRPTW-ColumnGen}{$V, C, [s, t]$}
            \State $\Omega' \gets \{\text{initial feasible routes}\}$ \Comment{e.g., one route per customer}
            \State $\bar{c}_{min} \gets -\infty$
            \While{$\bar{c}_{min} < 0$}
                \State Solve RMP Linear Relaxation:
                \State \quad $\min \sum_{\tau \in \Omega'} c_\tau y_\tau$
                \State \quad s.t. $\sum_{\tau \in \Omega'} a_{i\tau} y_\tau = 1, \forall i \in V$
                \State Extract dual variables $\pi_i$ from coverage constraints.
                \State $(\tau_{new}, \bar{c}_{min}) \gets \text{PricingSubproblem}(V, C, [s, t], \pi)$
                \If{$\bar{c}_{min} < 0$}
                    \State $\Omega' \gets \Omega' \cup \{\tau_{new}\}$
                \EndIf
            \EndWhile
            \State $y^* \gets \text{solve-MIP}(\Omega')$ \Comment{Solve with $y_\tau \in \{0,1\}$}
            \State \textbf{return} $y^*$
        \EndProcedure
        
        \Statex
        \Procedure{PricingSubproblem}{$V, C, [s, t], \pi$}
            \State Find route $\tau$ minimizing $\bar{c}_\tau = \sum_{(i,j) \in \tau} (c_{ij} - \pi_i)$
            \State s.t. $\text{Arrival}_i \in [s_i, t_i]$ for all $i \in \tau$
            \State \textbf{return} $\tau_{best}, \bar{c}_{best}$
        \EndProcedure
    \end{algorithmic}
\end{algorithm}

As a last note, if the time windows are not very tight. There might be some resulting tours with cycles. In this case, the simplest approach would be to eliminate them during the DP recursion, by checking that at every step a node is not revisited in the last $\ell$ steps.


\chapter{Systemic Demand Uncertainty and Chain Effects}
\section{Introduction}
Supply chains are complex networks composed of entities that order and selsl products: for example a reseller or a distributor orders from a factory and resells product to other customers, that in turn might be further reselling or be final customers. In our exposition we have emphasized on the random nature of the demand. In this chapter we study quantitatively how risk propagates in a supply chain network. 
 
Intuitively, when producing or selling a product, the optimal behavior is to produce/order the average plus a buffer to cover for the variance of the demand depending on its distribution according to the profit structure of the business in an optimal way. In a supply chain, the risk propagates and \emph{bubble-up} throughout the network, making the risk of a stockout higher. Lead times, that is, the time orders take to be fulfilled also adds a layer of propagation of risk that could affect the supply chain adversely. We formalize these insights into a powerful paradigm (based on Stochastic Networks theory, see \S 7.3 \cite{chen2001fundamentals}) that will allow us to study these effects in detail quantitatively.

\section{Basic Model}

To start, consider a single business buying/producing certain product and selling it to many different customers (potentially other businesses), all participants constitute a network with $N$ nodes. As we studied in Chapter \ref{sec:inventory_management}, when the demand is normally distributed (which in this case simplifies the analysis greatly), the optimal quantity for a business $i$ is to order a quantity:
\begin{eqnarray}
	 q_i^{NV}=\mu_i+z_i\sigma_i,
\end{eqnarray}
Where $\mu_i$ is the average demand, $\sigma_i$ is the demand standard deviation and $z_i=\Phi^{-1}(\rho_i)$ is the optimal factor in the Newsvendor model, where $\rho_i$ is the critical ratio given by the profit structure in the objective function (for example $\rho_i=(p_i-c_i)/p_i$ in the plain vanilla specification in Chapter \ref{sec:newsvendor}).

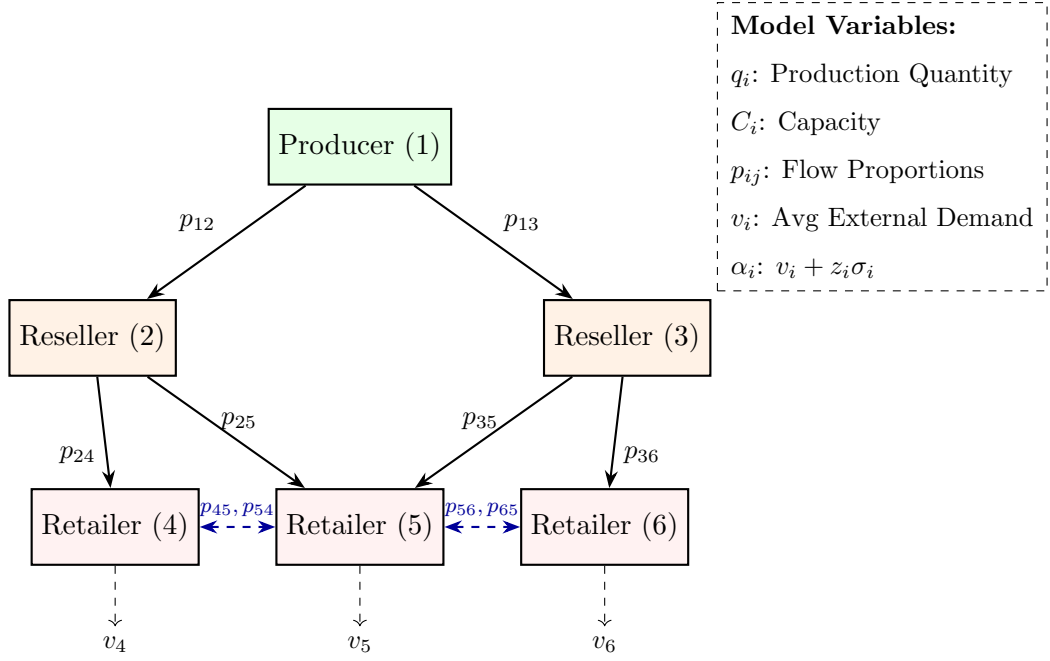
\begin{figure}[htb]
    \centering
\begin{tikzpicture}[
    node distance=2cm and 1.5cm,
    entity/.style={rectangle, draw, thick, minimum width=2.2cm, minimum height=1cm, fill=blue!5},
    arrow/.style={-Stealth, thick},
    transship/.style={Stealth-Stealth, dashed, thick, blue!60!black}
]

    
    \node[entity, fill=green!10] (P1) {Producer ($1$)};

    \node[entity, below left=1.5cm and 1.2cm of P1, fill=orange!10] (RS1) {Reseller ($2$)};
    \node[entity, below right=1.5cm and 1.2cm of P1, fill=orange!10] (RS2) {Reseller ($3$)};

    \node[entity, below=4cm of P1, fill=red!5] (RT2) {Retailer ($5$)};
    \node[entity, left=1cm of RT2, fill=red!5] (RT1) {Retailer ($4$)};
    \node[entity, right=1cm of RT2, fill=red!5] (RT3) {Retailer ($6$)};

    
    \draw[arrow] (P1) -- node[above left, font=\small] {$p_{12}$} (RS1);
    \draw[arrow] (P1) -- node[above right, font=\small] {$p_{13}$} (RS2);

    \draw[arrow] (RS1) -- node[left, font=\small, pos=0.7] {$p_{24}$} (RT1);
    \draw[arrow] (RS1) -- node[right, font=\small, pos=0.4] {$p_{25}$} (RT2);
    
    \draw[arrow] (RS2) -- node[left, font=\small, pos=0.4] {$p_{35}$} (RT2);
    \draw[arrow] (RS2) -- node[right, font=\small, pos=0.7] {$p_{36}$} (RT3);

    \draw[transship] (RT1.east) -- node[above, font=\scriptsize] {$p_{45}, p_{54}$} (RT2.west);
    
    \draw[transship] (RT2.east) -- node[above, font=\scriptsize] {$p_{56}, p_{65}$} (RT3.west);

    \draw[dashed, ->] (RT1.south) -- ++(0,-0.8) node[below, font=\small] {$v_4$};
    \draw[dashed, ->] (RT2.south) -- ++(0,-0.8) node[below, font=\small] {$v_5$};
    \draw[dashed, ->] (RT3.south) -- ++(0,-0.8) node[below, font=\small] {$v_6$};

    \node[right=of P1, xshift=2cm, align=left, font=\small, draw, dashed, inner sep=5pt] {
        \textbf{Model Variables:}\\
        $q_i$: Production Quantity\\
        $C_i$: Capacity\\
        $p_{ij}$: Flow Proportions\\
        $v_i$: Avg External Demand\\
        $\alpha_i$: $v_i + z_i\sigma_i$
    };

\end{tikzpicture}
\caption{Example of a supply chain network with a single producer (represented by a node 1), 2 resellers (nodes 2,3) and 3 retailers (nodes 4,5,6) that potentially sell to each other as well as servicing the external demands ($v_4,v_5,v_6$).}
\label{fig:supply_chain_network}
\end{figure}

\noindent Node $i$ has an average demand $\mu_i$ that needs to be satisfied in equilibrium as:
\begin{eqnarray}
\mu_i=\sum_{j}p_{ij} q_j+v_i,
\end{eqnarray}
Where $q_j$ is the effective quantity produced by node $j$ and $p_{ij}\in (0,1)$ is the proportion of the demand of $j$ produced by node $i$ and $v_i$ is the external average demand. Then, assuming independence between the nodes (here the variance is assumed exogenous from the external demand), the variance $\sigma^2_i$ is equal to:
\begin{eqnarray}
	\label{eq:net_variance}
\sigma^2_i=\sum_{j}p_{ij}^2\sigma^2_j+w^2_i,
\end{eqnarray}
Where $w_i^2$ is the \emph{external variance} of the demand (which can be 0 if node $i$ does not serve an external demand). Then, the effective quantity that node $i$ needs to produce is given by:
\begin{eqnarray}
q_i=\min\left(\sum_{j}p_{ij} q_j+\alpha_i,C_i\right),
\end{eqnarray}
Where $\alpha_i=v_i +z_i\sigma_i$. When the constraint holds with equality, this is the optimal quantity according to the Newsvendor model (as $q_i=\mu_i+z_i\sigma_i=\sum_{j}p_{ij} q_j+\alpha_i$) or the maximum physical capacity $C_i$ of node $i$. Since this is true for all nodes in the network, then, the maximum flow of product in the supply chain network is obtained by solving the following LP:
\begin{eqnarray}
	&\max&|q|\\
	&\text{s.t. }&(I-P)q^\intercal\le \alpha^\intercal,\nonumber\\
	&&\bm 0\le q\le C.\nonumber
\end{eqnarray}
Where $P$ is a sub-stochastic matrix with entries $p_{ij}$ and $\alpha$ is a vector with entry $i$ equal to $\alpha_i=v_i+z_i\sigma_i$. For most nodes $v_i$ and $w_i$ will be assumed to be 0 (think of producers/resellers that do not sell directly to consumers) and it will be positive for retailers (or nodes that face external demand directly). $C$ is a vector with entry $i$ equal to $C_i$, the capacity of that node. In Figure \ref{fig:supply_chain_network} we see an example of a supply chain network with a producer, resellers and retailers.

To solve the problem there are multiple strategies: the simplest one is to solve directly the LP obtaining a solution $q^*$ for the maximum flow of product. Since there are $N$ nodes, a basic feasible solution of the problem will have exactly $N$ binding constraints in the optimal solution (ignoring potential degeneracies). For example, for the constraints $q\le C$ the nodes will be divided into two sets: one set of nodes $U\subseteq\{1,\dots,N\}$ that will be producing at capacity $q_U=C_U$, these nodes are \emph{un-hedged} against the risk implied by the Newsvendor model. In other words, they are limited by their capacity and with a non-negligible chance might suffer a stock-out of product. Conversely, the other set of nodes $\calh=\{1,\cdots,N\}\setminus U$ is producing at the maximum possible given the external demand and the structure of the network.

Then, the other flow constraints should hold with equality for the other variables, that is, $(I_\calh-P_\calh)q^\intercal_\calh+(I_{\calh,U}-P_{\calh,U})C^\intercal_U=\alpha^\intercal_\calh$. Which simplifies to $q^\intercal_\calh=(I_\calh-P_\calh)^{-1}(\alpha^\intercal_\calh+P_{\calh,U}C^\intercal_U)$. Then, it must be that the optimal is given by:
\begin{eqnarray}
q^\intercal_U=C^\intercal_U,\,\,\,\, q^\intercal_\calh=(I_\calh-P_\calh)^{-1}(\alpha^\intercal_\calh+P_{\calh,U}C^\intercal_U).
\end{eqnarray}
The optimality of this solution can be verified using duality checking that the solution satisfies the complimentary slackness conditions. 

More importantly, the solution highlights the structure of the problem: namely, some nodes will be capped by their capacity and will not be producing enough to optimally hedge their risk. 

\subsection{The bullwhip effect and other sensitivities}
In the supply chain literature there is an empirical principle known as the \emph{bullwhip effect}. The bullwhip effect is the downstream increase in demand throughout the supply chain network caused by an increase in customer demand, which in turn ripples increasingly throughout the network.

In our model this and other effects can be quantified, we will see that the bullwhip effect is not only caused by increases in demand, but it is also exacerbated by the structure of the supply chain network and how it also propagates the variance of the external demands. The first step is to calculate the sensitivities, which is a fancy name for the derivatives of the solution with respect to some variables. The sensitivity with respect to $\alpha_i$ (capturing both the sensitivity to the demand and the volatility) is given by:
\begin{eqnarray}
\frac{\partial q_\calh}{\partial \alpha_i}&=&(I_\calh-P_\calh)^{-1}e^\intercal_i\,\,\,\text{   for   }i\in\calh,\\
\frac{\partial q_\calh}{\partial C_j}&=&(I_\calh-P_\calh)^{-1}P_{\calh,U}e^\intercal_j\,\,\,\text{   for   }j\in U.
\end{eqnarray} 
Where $e_i$ is a basis vector with entry $i$ equal to 1 and the rest equal to zero. In both sensitivities, we can see the role of the matrix $(I_\calh-P_\calh)^{-1}$. This matrix is commonly called the \emph{network multiplier} matrix. Intuitively, recall that for a geometric series we have the well known result:
\begin{eqnarray}
1+\delta+\delta^2+\delta^3+\cdots=(1-\delta)^{-1},
\end{eqnarray}
Which is true as long as $|\delta|<1$, which intuitively means that the terms $\delta^k$ for $k\rightarrow\infty$ not only decrease to 0, but their order of magnitude vanish as the partial sum goes to infinity (the mathematical notion of convergence of a series). Imagine the following situation to get some intuition: An oil producer, generages a quantity $q=D+5\%q$ of oil, that is, their external demand $D$ plus 5\% of their own output that they consume to re-refine the oil. Rearranging the terms yields $q=\frac{D}{1-5\%}=D(95\%)^{-1}$. The total quantity is the demand multiplied by a \emph{feedback} factor $(95\%)^{-1}\approx105.3\%$. This can be seen as a feedback loop in the production caused by output that reenters the production process. The same happens in a network case with matrices, where we have that:
\begin{eqnarray}
	I+P+P^2+P^3+\cdots = (I-P)^{-1},
\end{eqnarray}
It is also expected that for the result to hold, the higher powers of $P^k$ go to 0 as $k\rightarrow\infty$ and vanish from the partial sums. Luckily, this holds true for substochastic matrices in our base model\footnote{In general, what is required for the result to hold is that the largest eigenvalue of the matrix is less than 1 in absolute value.}. Moreover, the same intuition applies, each product of the matrix can be seen as a feedback loop of an increase of the required production coming from neighboring nodes, that in turn increase the production requirements for their neighbor nodes. This back and forth stabilizes in the form of multipliers (like in the above paragraph example of $(95\%)^{-1}$), and the entries of $(I_\calh-P_\calh)^{-1}$ represent exactly the multiplier of node in a row to a node in a column. That is why the sensitivity $(I_\calh-P_\calh)^{-1}e^\intercal_i$ is the row sum of of all columns entering into node $i\in \calh$. Likewise, the sensitivity $(I_\calh-P_\calh)^{-1}P_{\calh,U}e^\intercal_j$ represents the weighted (by $P_{\calh,U}$) multipliers of the hedged nodes into the unhedged ones.

By the chain rule, it is straightforward to calculate the effect of a change in the demand as $\frac{\partial q_\calh}{\partial v_i}=\frac{\partial q_\calh}{\partial \alpha_i}\frac{\partial \alpha_i}{\partial v_i}$ and $\frac{\partial \alpha_i}{\partial v_i}=1$ recalling $\alpha_i=v_i+\sigma_iz_i$. Therefore, $\frac{\partial q_\calh}{\partial v_i}=(I_\calh-P_\calh)^{-1}e^\intercal_i$ is the same sensitivity as before.

{\bf The bullwhip effect} is known as the empirical phenomenon when an increase of the demand ripples with increasing magnitude accross the supply chain. With our framework, we can study permanent increases in the demand, as well as increases in the volatility that as we will see, quantitatively ripple differently in the network. Suppose, the demand in one of the customer facing nodes increases, that is, $v_i$ increases for some node $i$. We are then interested in quantifying the impact of this increase in the demand.

To make the analysis easier, suppose the structure of the network is as follows: the graph is an acyclic directed graph and the matrix $P$ is substochastic, where nodes are ordered in a way that they provide their output to nodes numbered higher than them. For example, in the graph in Figure \ref{fig:supply_chain_network} the node 1 represents a producer, which supplies their output to the reseller nodes 2 and 3. In turn, node 2 sells their output to retailers 4 and 5, and node 3 sells to 5 and 6. Note that the retailer nodes, 4 sells to 5, and 5 sells to 4 which violates our assumption. So, imagine the same network, but without node 5 selling to 4 and withouth 6 selling to 5. Later we will see that this is not really a necessary condition, but it makes the next argument easier to explain.

Imagine there is a permanent shock to the demand in a node that increases it to a new level, for example, in the graph in Figure \ref{fig:supply_chain_network} without the cycles, imagine that node 6 demand increases to a new level $v'_6=v_6+\epsilon$ from the previous $v_6$. This increase affects first node $6$, which increases their production by an amount $\epsilon$. Next, the nodes that are affected are nodes 3 (a reseller) and 5 (another retailer). These in turn increase their demand by a magnitude $\epsilon p_{36}$ and $\epsilon p_{56}$ respectively. This continues happening so forth throughout all the network (and note that sometimes the effect will go forward again, as after node 3 increases their production, so will the output that node 5 does).

A way to analyze this effect numerically, is to look at the \emph{network multipliers}. Moreover, suppose that all nodes have infinite capacity $C_i\rightarrow \infty$, so that all nodes are producing at their hedged level, so that $q=(I-P)^{-1}\alpha^\intercal$. Since the network is free of cycles, the multiplier matrix $(I-P)^{-1}$. The multiplier effect of a shock in node $j$ into node $i$ is the element $[(I-P)^{-1}]_{i,j}$. In the expanded form this is:
\begin{eqnarray}
	[(I-P)^{-1}]_{i,j}=I_{i,j}+P_{i,j}+P^2_{i,j}+\cdots+P^{N-1}_{i,j},
\end{eqnarray}
First note that there are $N$ elements in the summation, as there no longer than $N$ steps between 2 nodes. Next, see that each term exactly reflects how an increase in the demand travels throughout the network. The first node affected is $I_{i,j}$ that will feel the full effect only if $i=j$ (for example, the node 6 in the previous paragraph). Next, the effect is the nodes that can be reached by one step $P_{i,j}$ (the nodes 3 and 6). Next, the ones that can be reached exactly in two steps $P^2_{i,j}$, in 3 steps $P^3_{i,j}$ and so on (in fact, the same analysis holds for matrices with cycles, instead of $N+1$ there is a large enough term $K$ from which further terms vanish). 

This is exactly what the bullwhip effect is, nodes that are directly affected see a lower order of magnitude effect, than deeper nodes (producers or reseller) as the shock is transmitted and amplified throughout all the network.

\subsection{Volatility Propagation}
 
In the literature the Bullwhip effect is also associated to the propagation of variance throughout the supply chain (affecting again the innermost nodes rather than the consumer facing ones). In our model this mechanism is captured by the structure of the variance. We can write the variance Equation in (\ref{eq:net_variance}) as:
\begin{eqnarray}
	(\sigma^2)^\intercal=(P\odot P)	(\sigma^2)^\intercal+(w^2)^\intercal,
\end{eqnarray}
Where $\sigma^2=(\sigma_1^2,\dots,\sigma_N^2)$, $w^2=(w_1^2,\dots,w_N^2)$ and $\odot$ is the element-wise product of two matrices of the same dimension, such that the $i,j$-th element of the matrix $(P\odot P)$ is $p_{ij}^2$. The equation becomes $(I-(P\odot P))(\sigma^2)^\intercal=(w^2)^\intercal$, or:
\begin{eqnarray}
	(\sigma^2)^\intercal=(I-(P\odot P))^{-1}(w^2)^\intercal.
\end{eqnarray}
This expression should be already familiar from the previous exposition. The matrix $(I-(P\odot P))^{-1}$ is again a multiplier, but this time of the variance, denote as $M_\sigma=(I-(P\odot P))^{-1}$. Similar to the previous case, this shows how the variance of the external demand $w^2$ propagates inside the network. Taking element-wise square root in both sides yields the equivalent expression $\sigma^\intercal=\sqrt{M_\sigma(w^2)^\intercal}$.

The change in the variance of nodes by a change in the external demand standard deviation $w_j$ in node $j$ is given by:
\begin{eqnarray}
	\frac{\partial \sigma^\intercal}{\partial w_j}=\text{diag}\left(1/\sqrt{M_\sigma(w^2)^\intercal}\right)m_jw_j,
\end{eqnarray}
Where $m_j$ is a column vector of $M_\sigma$. Then, we can finally express the sensitivity of the demand with respect to changes in the external variance. We have by the chain rule that $\frac{\partial q_\calh}{\partial w_i}=\frac{\partial q_\calh}{\partial \alpha_i}\frac{\partial \alpha_i}{\partial w_i}$ and in turn $\frac{\partial \alpha_i}{\partial w_i}=z_i\frac{\partial \sigma_i}{\partial w_i}=z_i\frac{m_{ii}}{\sqrt{M_\sigma(w^2)^\intercal}_i}w_i$ or:
\begin{eqnarray}
	\frac{\partial q_\calh^\intercal}{\partial w_i}=(I_\calh-P_\calh)^{-1}e^\intercal_i\left[\frac{m_{ii}}{\sqrt{M_\sigma(w^2)^\intercal}_i}z_iw_i\right].
\end{eqnarray}
The intuition is quite interesting: on one hand there is a network multiplier effect as before, but there is a further variance multiplier effect (in square brackets) that further exacerbates a shock to the standard deviation of the external demand. Note that the shock is also proportional to the current standard deviation $w_i$, accentuating shocks even further.

Outdated interpretations of the Bullwhip effect attribute it to the lack of information and coordination between nodes or noise coming from forecasting and lead times. Suggesting that the lack of better information systems is the main reason behind the bubbling-up of demand shocks.

In this chapter we showed mathematically that even with perfect information and coordination (all the nodes in the network can see the production quantities and uncertainty of adjacent nodes), {\bf the Bullwhip effect is an unavoidable and an inherent topological property of supply chain networks}. Even in the perfect information case, with instantaneous order fulfillment, hedging against risk implies the Bullwhip effect.

 \begin{exm}[{\bf Rare earth supply chain}]
	A supply chain network composed of a producer (node 1) and a seller/country (node 2). The external average demand $v_2=100$ tons per year, the external standard deviation of the demand $w_2=10$ tons. The seller buys all of its product from the producer, that is $p_{12}=1$, the routing matrix is:
	\begin{eqnarray*}
		P=\begin{pmatrix}
			0 & 1 \\
			0 & 0
		\end{pmatrix},
	\end{eqnarray*}
	With a safety level $z_1=z_2=1.645$ what are the optimal production quantities in this network? While maintaining the same average demand, how much an increase of 1 ton of the standard deviation of the demand increases the quantities for the producer?
	
	\noindent For the first question, we have the following quantities:
	\begin{eqnarray*}
		q_2=v_2+z_2w_2&=&100+1.645(10)=116.45\\
		q_1=p_{12}q_2+z_1\sqrt{p_{12}\sigma^2_2}&=&116.45+1.645(10)=132.9
	\end{eqnarray*}
	Even when the network only has 2 nodes, the producer has to produce way above node 1 to be optimally hedged against the demand fluctuation of the seller.
	
	\noindent For the second question, we make use of the network multiplier matrix:
	\begin{eqnarray*}
		(I-P)=\begin{pmatrix}
			1 & -1 \\
			0 & 1
		\end{pmatrix}\Rightarrow
		(I-P)^{-1}=\begin{pmatrix}
			1 & 1 \\
			0 & 1
		\end{pmatrix},
	\end{eqnarray*}
	Also note that $P=P\odot P$, then the network multiplier matrix $(I-P)^{-1}$ and the variance multiplier $M_\sigma$ are the same. For node 1 we have that the sensitivity is:
	
	\begin{equation}
		\frac{\partial q_1}{\partial w_2} = [(I - P)^{-1}]_{1,2} \times \left[ \frac{m_{22}}{\sigma_2} z_2 w_2 \right]=1\times\left[\frac{1}{10}1.645(10)\right]=1.645
	\end{equation}
	\noindent Then, an increase of the demand volatility by 1 ton increases production by 1.645 tons in the producer node.
\end{exm}

\section{Dynamic Model}
In practice, knowing only the new equilibrium quantities is important, but also the timing on how the network adapts to shocks and their magnitude. In this section we extend our model to a dynamic setting to answer questions related to the timing and magnitude of the \emph{Bullwhip effect} and other shocks. The standard way to incorporate a dynamic version, is to solve a so-called \emph{Skorokhod Problem} (see \cite{skorokhod1961stochastic}). 

Consider the evolution of the inventory if there are no capacity constraints, a so-called \emph{free inventory process}:
\begin{eqnarray}
  X_i(t)=Z_i(0)+(\text{Production}_i-\text{Internal Demand}_i-\text{External Demand}_i)t,
\end{eqnarray}
Where $Z_i(0)$ is an initial level of inventory (not to be confused by the safety margin $z_i$). Call the expression in parenthesis $\theta_i$ showing the net-rate that the inventory changes. In matrix-form, the free inventory process can be written as:
\begin{eqnarray}
X(t)=Z(0)+\theta t,
\end{eqnarray}
In this case, $\theta = C-\alpha-(\lambda\wedge C)$ as the maximum production capacity of a node is their capacity, that is, $\lambda\le C$, recalling $C=(C_1,\dots,C_N)$. $\theta_i$ denotes the rate the inventory increases if $\theta_i>0$ or decreases $\theta_j<0$, eventually reaching 0. Let $\lambda$ be the effective production rate of each node.

Let $Z(t)$ be the actual inventory process. As the inventory cannot be negative, that is, $Z(t)\ge 0$, the Skorokhod problem trick is to create a function $Y(t)$ that cancels out when the inventory is negative. As this must happen for all the nodes we have the relation:
\begin{eqnarray}
  Z(t)=X(t)+Y(t)-Y(t)P=X(t)+Y(t)(I-P),
\end{eqnarray}
In this case, $Y(t)$ appears first as it is cancelling out when $X(t)$ is negative, but it also must account when other nodes have a shortfall that affects node $i$, therefore the need for the term $-Y(t)P$, or $-\sum_jY_j(t)p_{ji}$. This keeps track of the propagation when other nodes have shortfalls, which is a mathematical way of tracking the contagion between nodes.

It must also be the case that when $Z(t)>0$, then $dY(t)=0$, because there is still enough inventory to satisfy incoming orders. Likewise, when $Z(t)=0$, then $dY(t)>0$, meaning that there is not enough production capacity to satisfy the demand. All the conditions can be summarized in the program:
\begin{eqnarray}
  &&Z(t)=X(t)+Y(t)(I-P),\\
  &\text{s.t.}&Y(0)=0,\,\,\,\dot{Y}(t)\ge0,\,\,\,Z(t)\ge0,\,\,\,Z(t)\dot{Y}(t)=0.\nonumber
\end{eqnarray}
Where $\dot{Y}$ denotes the time derivative of $Y$, that is, $\dot{Y}(t)=dY(t)/dt$.
In the long term, it must be the case that the equilibrium of the dynamic model is the same as the static one. To see why, first considering the case where all the production capacities are large enough so that any production level is possible. In this case, the evolution of the system is given by the equation:
\begin{eqnarray}
  \dot{Z}(t)=\dot{X}(t)+\dot{Y}(t)=\theta+\dot{Y}(t)(I-P),
\end{eqnarray}
If $Z(t)>0$ ends in a positive level, then $dY(t)=0$. Likewise, if $Z(t)$ reaches an equilibrium point, then $dZ(t)=0$. Leading to the equation:
\begin{eqnarray}
  0=\theta=C(I-P)-\alpha\,\,\, \underset{C=\lambda}{\Rightarrow}\,\,\, \lambda=\alpha(I-P)^{-1},
\end{eqnarray}
Where $\lambda$ represents the effective production rate of each node. In the limit, this is the same production quantity that we derived in the static case. That is $\lim_{t\rightarrow \infty}\lambda(t)=q^*$. In practice, we will see that because of their physical capacity, some nodes will be unable to produce their required quantity and they will be again partition into disjoint hedged $\calh$ and un-hedged $U$ nodes. In the case where nodes hit a capacity limit, this refines to the so-called \emph{flow equation}:
\begin{eqnarray}
  \label{eq:flow}
  \lambda = \alpha+(\lambda\wedge C)P,
\end{eqnarray}
Where $\wedge$ is the element-wise minimum operator. When $\lambda\wedge C=\lambda$ the previous equation holds.

Following our previous analysis, at any point $t$ the nodes are divided into two disjoint sets of hedged $\calh$ and un-hedged $U$ sets. For a given time $t$ that nodes are divided into hedge $\calh$ and un-hedged sets $U$. Then, the flow Equation \ref{eq:flow} refines to:
\begin{eqnarray}
  \lambda_\calh&=&\alpha_\calh+\lambda_\calh P_\calh+C_U P_{U,\calh},\\
  \lambda_U&=&\alpha_U+\lambda_\calh P_{\calh,U}+C_UP_U,
\end{eqnarray}
Solving yields the following flows at time $t$:
\begin{eqnarray}
  \lambda_\calh&=&(\alpha_\calh+C_UP_{U,\calh})(I-P_\calh)^{-1},\\
  \lambda_U&=&\alpha_U+(\alpha_\calh+C_UP_{U,\calh})(I-P_\calh)^{-1}P_{\calh,U}+C_UP_U,
\end{eqnarray}
Note that as time moves forward, the set of hedged nodes might decrease as a result of increased demand. The rates the inventory is accumulated are given by:
\begin{eqnarray}
  \theta = C-\alpha-(\lambda\wedge C)P=C-\lambda
\end{eqnarray}
\subsection{Sequential dynamics}
In principle, we have all the ingredients to describe the evolution over time of the supply chain network. The nodes start at $t=0$ with some inventory position $Z(0)\ge0$ with production rates $\lambda(t)$ dividing the nodes into disjoint sets $\calh(t)$ and $U(t)$. Suppose, the sets are initially known as well, that is $\calh(0)$ and $U(0)$. Then, we just need to describe the evolution of the system until the time $\tau$ when the sets change (possibly because a set of nodes had to increase production until a time that their production limit is reached). When that time $\tau$ is reached we can just restart the clock and simulate the system again.

The evolution of the system is linear between time 0 and the time where the sets of nodes change. For all nodes already at capacity $Z_U(0)=0$ and $Z_\calh(0)>0$. We also have the following dynamics for the hedged nodes:
\begin{eqnarray}
  Z_\calh(t)=Z_\calh(0)+\theta_\calh t+Y_U(t)(I-P_{U,\calh})
\end{eqnarray}
Where $Y_U(t)=(\lambda_U-C_U)t$ and $\theta_\calh=C_\calh-\lambda_\calh$. 

The system keeps evolving linearly according to these dynamics until the time $\tau$ where a node in the set $\calh$ cannot satisfy all their production obligations. This time is given by:
\begin{eqnarray}
  \tau=\min_{j\in\calh,\theta_j<0}\left\{\frac{Z_j(0)}{-\theta_j}\right\},
\end{eqnarray}
This follows from the fact that a node in $\calh$ with $\theta_j<0$ implies that $C_j<\lambda_j$, meaning that even at maximum capacity they will eventually reach their production limit. This will happen when $Z_j(0)+\theta_jt=0$ or at $t=-Z_j(0)/\theta_j$.


\chapter{Queueing Models}
\label{sec:queueing}
\section{Introduction}
In this chapter we continue studying operations in the time dimension. For example, in a production activity or a warehouse operation, the order that tasks are performed affects the final result. For example, in a production line, some stations might process jobs faster than others, and the throughput will be affected by how the production line is configured. In a warehouse operation, when many trucks arrive at the same time, the number of loading stations and workers will affect the effective throughput of the operation. Knowing how to assign and configure manufacturing systems to be faster and more resilient is a key strategic tool to develop supply chain operations.

Mathematically, we leverage our study of networks and probabilistic theory to describe common problems in queueing systems.

\section{A network view of a queueing system}
Queueing systems are defined as a probabilistic system where jobs, customers or tasks \emph{evolve} inside a well limited system: for example, a production line queue, cars passing through a bridge or a flower plantation.

The primitive unit in all these systems is the number of jobs (or customers, tasks, etc) in the system. At any time $t$ the number of jobs $N=0,1,2,\dots$ is a unique value, this feature makes us think that in time-homogeneous systems (the arrivals, processing and departures are time-independent), we can think of a probabilistic distribution of observing the number in system as a proportion of time:
\begin{eqnarray}
  \pi_i:=\lim_{t\rightarrow\infty}\int_{0}^t\frac{\ind\{X_s=i\}}{t}ds
\end{eqnarray}
Where $X_t$ is the stochastic process counting the number in system at time $t$.

The time-homogenous process by excellence is the Markov Chain. We previously studied it in Chapter \ref{sec:inventory_management}. In fact, this is the same definition of stationary distribution that we described in the study of inventory processes. Then, it is a reasonable goal to define a type of Markov Chains that can model manufacturing processes.
The main class of Markov Chains used to model manufacturing processes are known as \emph{Continuous Time Markov Chains} or CTMC. The simplest way to construct such chain requires two ingredients: First, a transition matrix $P$ between states $N$, with elements $p_{ij}$ denoting the probability of transitioning from state $i$ to state $j$. Second, an average time the chain stays in a state that is also time homogeneous.

For a random time $T$ to be time-homogeneous it is necessary that for any $s,t\ge 0$, we have $\pr(T\ge t)=P(T\ge t+s|T\ge s)$. Meaning that the tail CDF is the same regardless of time-shifts. Using Bayes's rule and rearranging yields the relation:
\begin{eqnarray}
  \pr(T\ge t)\pr(T\ge s)=\pr(T\ge t+s)\,\,\text{ or }\,\, R(t)R(s)=R(t+s),
\end{eqnarray}When defining $R(t):=P(T\ge t)$. Then, the goal is to find a distribution that satisfies this property. This implies that for $\Delta t=t/k$ we have $R(t)=R(\Delta t)^k$. In the discrete case, $P(T\le i)=1-R(i+1)$ and also $R(i+1)=R(1)^{i+1}$. Putting these two relations together yields $P(T\le i)=1-R(1)^{i+1}$, letting $1-p=R(1)$, we will call this a \emph{geometric distribution}.

\begin{figure}[htbp]
    \centering
    \begin{tikzpicture}
        \begin{groupplot}[
            group style={group size=2 by 1, horizontal sep=2.0cm},
            width=0.48\textwidth,
            height=6.5cm,
            grid=major,
            ymin=0, ymax=1.1,
            ylabel style={font=\small},
            xlabel style={font=\small},
            title style={font=\bfseries\small, yshift=2pt}
        ]

        \nextgroupplot[
            title={Continuous: Exponential Tail $R(t) = e^{-\mu t}$},
            xlabel={Time ($t$)},
            ylabel={Survival Probability $P(T > t)$},
            xmin=0, xmax=4.5,
            xtick={0, 1.5, 3.5},
            xticklabels={0, $s$, $s+t$}
        ]
        \addplot[blue, ultra thick, domain=0:4.5, samples=100] {exp(-0.65*x)};
        
        \draw[dashed, black!60, thick] (axis cs:1.5,0) -- (axis cs:1.5,0.377);
        \draw[dashed, black!60, thick] (axis cs:3.5,0) -- (axis cs:3.5,0.105);
        
        \node[anchor=south] at (axis cs:0.75, 0.05) {\footnotesize Interval $s$};
        \node[anchor=south] at (axis cs:2.5, 0.05) {\footnotesize Interval $t$};
        
        \node[draw, fill=white, rounded corners, align=left, font=\scriptsize] at (axis cs:2.8,0.75) {
            $\frac{R(s+t)}{R(s)} = R(t)$\\[4pt]
            Ratio is invariant\\
            to the shift $s$.
        };

        \nextgroupplot[
            title={Discrete: Geometric Tail $R(k) = (1-p)^k$},
            xlabel={Steps ($k$)},
            ylabel={Survival Probability $P(T > k)$},
            xmin=-0.5, xmax=6.5,
            xtick={0, 2, 5},
            xticklabels={0, $s$, $s+k$}
        ]
        \addplot[ybar, fill=orange!40, draw=orange!80!black, thick, domain=0:6, samples=7, ybar legend] {(1-0.35)^x};
        
        \draw[dashed, black!60, thick] (axis cs:2,0) -- (axis cs:2,0.422);
        \draw[dashed, black!60, thick] (axis cs:5,0) -- (axis cs:5,0.116);

        \node[anchor=south] at (axis cs:1.0, 0.02) {\footnotesize Steps $s$};
        \node[anchor=south] at (axis cs:3.5, 0.02) {\footnotesize Steps $k$};
        
        \node[draw, fill=white, rounded corners, align=left, font=\scriptsize] at (axis cs:4.2,0.75) {
            $R(s+k) = R(s)R(k)$\\[4pt]
            Each step drops by\\
            a constant factor $(1-p)$.
        };

        \end{groupplot}
    \end{tikzpicture}
    \caption{The visual architecture of memorylessness. The continuous Exponential curve (left) and discrete Geometric bars (right) represent the unique configurations where the relative remaining lifetime layout above a survival milestone $s$ maintains a constant, self-similar proportion directly identical to the original starting profile at point $0$.}
    \label{fig:memoryless_property}
\end{figure}
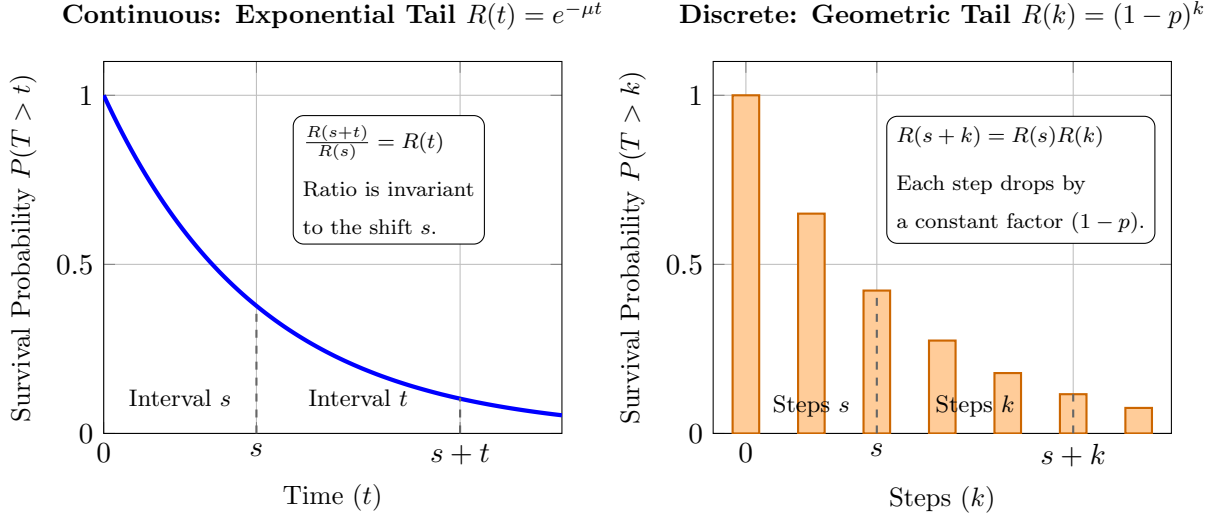

In the continuous case, something very similar can be done. For example, $R(t)=R(1)^t$. Changing the base of the exponential\footnote{$x=e^{\log x}$, then $x^t=e^{t\ln x}$.} we have $R(t)=R(1)^t=e^{t\ln R(1)}$. Letting $\mu=-\ln(R(1))>0$ we have a distribution with tail cdf $R(t)=e^{-\mu t}$, we call this distribution the \emph{exponential distribution}. See Figure \ref{fig:memoryless_property} for a visual representation of the survival probabilities.

Having sorted out the time-homogeneous aspect of the average time in a state (what we need to do is simply calibrate the parameters of the discrete geometric distribution $p$ or the exponential $\mu$ so that they match a target expected value). For example, in the exponential distribution case, we have that the average time in a state $i$, denoted as $T_i$ is equal to $\ex(T_i)=\int_0^\infty e^{-\mu_i t}dt=\left[-\frac{1}{\mu_i}e^{-\mu_i t}\right]_{0}^\infty=1/\mu_i$.

Now, we can think of a very small time interval $\Delta t$ and think of the transition probability $P_{ij}(t)$ of our Markov Chain where time passes continously. The first observation, is that from the memoryless property the probability of leaving a state is a constant times the time interval, for example, in the exponential case this is $\mu_i\Delta t$. Therefore, the probability of staying in the same state is $1-\mu_i\Delta_t$. If we extend this logic to all the states, we end up with this relation:
\begin{eqnarray}
  P_{ij}(t+\Delta t)\approx P_{ij}(t)(1-\mu_j \Delta t)+\sum_{k\neq j}P_{ik}(t)(p_{kj}\mu_k\Delta t),
\end{eqnarray}
The interpretation of this equation is that at time $t+\Delta t$ one of two outcomes occured: first, the chain simply stayed at state $j$. Or second, the chain transition to other state k, left that state (with probability $\mu_k \Delta t$) to state $j$ (with probability $p_{kj}$). Rearranging this term, we get:
\begin{eqnarray}
  \frac{P_{ij}(t+\Delta t)-P_{ij}(t)}{\Delta t}\approx-P_{ij}(t)\mu_j+\sum_{k\neq j}P_{ik}(t)(p_{kj}\mu_k),
\end{eqnarray}
Which is the definition of the derivative (or change) when $\Delta t$ is small, leading to the equation:
\begin{eqnarray}
  \frac{d P_{ij}(t)}{dt}=-P_{ij}(t)\mu_j+\sum_{k\neq j}P_{ik}(t)(p_{kj}\mu_k),
\end{eqnarray}
When grouping all the elements, we get the equation $\frac{dP(t)}{dt}=P(t)Q$, where $Q$ is a matrix with elements $q_{ii}=-\mu_i$ and $q_{ij}=p_{ij}\mu_i$. For readers familiar with calculus, the differential equation should be familiar: if we were to imagine $P(t)$ is a single variable equation (which it IS NOT, is a matrix), we could rearrange terms to $\frac{dP(t)}{P(t)}=Qdt$, integrate on both sides, getting $\ln P(t)-\ln P(0)=Qt$. Letting $P(0)=1$ (or actually the identity matrix), then $\ln(P(0))=0$. Taking exponential on both sides yields $P(t)=e^{Qt}$. Turns out that this intuition is also correct in the matrix case\footnote{Instead of integrating, use the fact that $e^{Qt}=I+Qt+(Qt)^2/2!+(Qt)^3/3!+\cdots$ Guess that $P(t)=e^{Qt}$, then $\frac{P(t)}{dt}=Q+Q^2t+Q^3t^2/2!+\cdots$. Factor the $Q$ to get exactly $\frac{P(t)}{dt}=QP(t)$. Lastly, $P(0)=I$, that is the system is in the initial state with probability 1.}.

What we get from this is remarkable: the fact that to describe all the transitions we just need a matrix of constants $Q$, and that the transition is only proportional to the time interval $\Delta t$. The fact that $P(t)=e^{Qt}$ also gives us for free the equilibrium distribution. Assuming a very small time interval $\Delta t$, we have $P(\Delta t)\approx I+Q\Delta t$ from the Taylor expansion in the previous footnote. We know that the equilibrium distribution follows:
\begin{eqnarray}
  \pi P(\Delta t)=\pi,
\end{eqnarray}
Then, we can just replace to get $\pi(I+Q\Delta t)=\pi$, which simplifies to the equilibrium distribution in CTMC:
\begin{eqnarray}
  \pi Q=\bm 0,
\end{eqnarray}
Solving this system gives us exactly the equilibrium distribution of a system that evolves according the the probabilities $p_{ij}$ and stays in state $i$ for an exponential amount of time $T_i$.

As we studied in previous chapters, getting a close expression for the equilibrium distribution is a powerful modeling tool to understand the behaviour of a stochastic system. We illustrate this with an example:

\begin{exm}
  A strait (or a man-made waterway) is a passage connecting two large bodies of water. Imagine a country wants to inspect and charge ships as they cross. They can inspect up to 5 ships at the same time, in 5 stations. Each inspection takes in average $\ex(\tau)=2$ hours per ship per station. The demand for the use of the strait depends on the fee ships are charged to pass. The random quantity of arriving ships $a$ depends on this fee as $\ex(a)=e^{-5p}$ ships per hour where $p$ is the fee a ship is charged per unit of time in the strait. Given $p$ what is the average revenue that the country receives per hour assuming $a(p)$ ships arrive per hour? Also, what is the average amount of time a ship spends in this system? If instead of charging per hour, the country charges per ship, what's the new revenue?

  The first step in solving this problem is defining the state space, that is, the set of states that we are going to analyze. From the statement of the problem there can be up to 5 ships in processing at any given time. Nonetheless, it is easy to imagine that a ship can arrive while 5 others are in processing, and potentially one more and so on. Also, the system can sometimes be empty. Then, it makes sense to let the state of the system be the integers $N=0,1,\dots$ representing the number of ships in the system. The strategy will be to estimate the average number of ships in the system (at any time), that is, $\ex(N)=0\pi_0+1\pi_1+2\pi_2\dots$. The system earns $p\ex(N)$ in the long term. To solve for $\pi$, we need to find it with $\pi Q=\bm 0$.
  To do this, we simply need to describe how the system evolves: from any state with $i>0$ ships in the system there are two posible outcomes, either a ship finish their service from a station or a new ship arrives. The arrival of ships is $e^{-5p}$ per hour, so $\lambda_i:=q_{i,i+1}=e^{-5p}$. Likewise, because of the memoryless property of the exponential distribution, the minimum of $k$ exponential random variables $min(\tau_1,\dots,\tau_k)$ is another exponential distribution\footnote{Note that $\pr(min(\tau_1,\dots,\tau_k)\ge t)=\pr(\tau_1\ge t,\dots,\tau_k\ge t)=\pr(\tau_1\ge t)\cdots\pr(\tau_k\ge t)=e^{-t/\ex{\tau_1}}\cdots e^{-t/\ex{\tau_k}}=e^{-t\left[1/\ex{\tau_1}+\cdots+1/\ex{\tau_k}\right]}$. That is, an exponential distribution with a rate equal to the sum of the rates of the other exponentials.} with rate equal to the sum of the rates of the $k$ exponentials. In our case, there are at most 5 stations or less active, then the transition occurs at a rate $\mu_i:=q_{i,i-1}=\min(i,5)\frac{1}{\ex(\tau)}=\min(i,5)/2$. Finally, $q_{i,i}=-(q_{i,i-1}+q_{i,i+1})$. With this entries of $Q$ solving the equilibrium $\pi Q=\bm0 $ gives the equilibrium distribution of the system, from which $p\ex(N)$ is the rate the system earns.

  This type of process is commonly called a \emph{Birth-Death Process}, as from any state $i$, in principle there will be a transition to either $i+1$ or $i-1$ with rates $\lambda_i$ and $\mu_i$ respectively. Computing the equilibrium distribution is also easy for these processes, as the equilibrium equation looks like:
  \begin{eqnarray}
    -\pi_0\lambda_0+\pi_1\mu_1&=0\nonumber\\
    \pi_0\lambda_0-\pi_1\mu_1-\pi_1\lambda_1+\pi_2\mu_2&=0\nonumber\\
                                                       &\vdots\nonumber
  \end{eqnarray}
  The first equation implies $\pi_0\lambda_0=\pi_1\mu_1$, which in turn cancels out the same term in the second equation leaving $\pi_1\lambda_1=\pi_2\mu_2$. The same continues occuring, thus having $\pi_k\lambda_k=\pi_{k+1}\mu_{k+1}$. Replacing $\pi_1=\pi_0\frac{\lambda_0}{\mu_1}$ and $\pi_2=\pi_1\frac{\lambda_1}{\mu_2}=\pi_0\frac{\lambda_0\lambda_1}{\mu_1\mu_2}$ from which we can deduce that $\pi_k=\pi_0\frac{\lambda_0\cdots\lambda_{k-1}}{\mu_1\cdots\mu_k}$. Lastly, since $\pi$ must be a valid distribution (adding to 1) then we have $1=\pi_0+\pi_1+\cdots$, then $1=\pi_0[1+\lambda_0/\mu_1+\cdots]$. Which yields $\pi_0=[1+\lambda_0/\mu_1+\cdots]^{-1}$. After calculating these, the average money that the country earns is $p\ex(N)=p[0\pi_0+1\pi_1+2\pi_2+\cdots]$.

  Estimating the average waiting time seems hard, as in principle it would require to average many possible paths for different number of ships in the system. But we can try to use what we know to make this estimation. For any ship arriving at random, this ship will experience a random amount of time in the system $W$. This ship would have to pay $p\ex(W)$ in average. In an hour, there are $\ex(a)=e^{-5p}$ ships arriving per hour, then, the amount they pay should be $\ex(q)p\ex(W)$. Note that the amount the ships pay is in fact the same revenue we found before. Then, it must be true that $p\ex(N)=\ex(q)p\ex{W}$. This is known as \emph{Little's law}:
  \begin{eqnarray}
    \ex(W)=\frac{\ex(N)}{\ex(a)},
  \end{eqnarray}
  This law is extremely useful, because we got the expected value of the waiting time $\ex(W)$ for free, just knowing the arrival rate $\ex(a)$ and the average number in the system $\ex(N)$. In general, with this relation, we need only two of the quantities to guess the third. In most textbooks this law is stated as $\ex(N)=\ex(a)\ex(W)$. In Figure \ref{fig:littles_law_sample_path_perfect} there is a graphical representation of the principle.

\begin{figure}[htbp]
    \centering
    \begin{tikzpicture}[
        every node/.style={font=\small},
        axis/.style={-={Stealth[scale=1.0]}, thick},
        curve/.style={blue, ultra thick},
        gridline/.style={help lines, color=gray!30, dashed}
    ]

        \begin{scope}[local bounding box=leftplot]
            \node[anchor=south, font=\bfseries] at (2, 3.9) {View 1: Time Snapshot};
            \node[anchor=south, font=\scriptsize, text width=5.2cm, align=center, color=gray] at (2, 3.0) 
                {\emph{Counting the proportion of time the system spends at each state}};

            \draw[gridline] (0,0) grid (4,3);
            \draw[axis] (0,0) -- (4.8,0);
            \node[anchor=north west] at (4.3, -0.1) {Time ($t$)}; 
            \draw[axis] (0,0) -- (0,3.3);
            \node[anchor=south, left=3pt] at (0,3.3) {Ships ($N$)}; 
            
            \node[anchor=north] at (0,0) {$0$};
            \node[anchor=north] at (3.8,0) {$T$};
            \node[anchor=east] at (0,1) {$1$};
            \node[anchor=east] at (0,2) {$2$};
            \node[anchor=east] at (0,3) {$3$};

            \draw[curve] (0,0) -- (0.5,0) -- (0.5,1) -- (1.2,1) -- (1.2,2) -- (1.8,2) -- (1.8,1) 
                         -- (2.2,1) -- (2.2,3) -- (2.9,3) -- (2.9,2) -- (3.3,2) -- (3.3,0) -- (3.8,0);
                         
            \fill[blue!15, opacity=0.6] (0.5,0) -- (0.5,1) -- (1.2,1) -- (1.2,2) -- (1.8,2) -- (1.8,1) 
                         -- (2.2,1) -- (2.2,3) -- (2.9,3) -- (2.9,2) -- (3.3,2) -- (3.3,0) -- cycle;

            \node[color=blue!80!black, font=\bfseries\footnotesize] at (1.8, 0.6) {Area = $\ex(N) \times T$};
        \end{scope}

        \draw[-={Stealth[scale=1.5]}, line width=1mm, gray!40] (5.2, 1.5) -- (6.6, 1.5) 
            node[midway, above=2pt, black, font=\bfseries\footnotesize] {Identical}
            node[midway, below=2pt, black, font=\bfseries\footnotesize] {Area};

        \begin{scope}[shift={(7.8,0)}]
            \node[anchor=south, font=\bfseries] at (2, 3.9) {View 2: Ledger Spreadsheet};
            \node[anchor=south, font=\scriptsize, text width=5.2cm, align=center, color=gray] at (2, 3.0) 
                {\emph{Summing up horizontal wait times ($W_i$) row by row for $C$ ships}};

            \draw[gridline] (0,0) grid (4,3);
            \draw[axis] (0,0) -- (4.8,0);
            \node[anchor=north west] at (4.1, -0.1) {Duration ($W$)}; 
            \draw[axis] (0,0) -- (0,3.3);
            \node[anchor=south, left=3pt] at (0,3.3) {Ledger Row}; 

            \node[anchor=east] at (0,0.5) {\scriptsize Ship 1};
            \node[anchor=east] at (0,1.2) {\scriptsize Ship 2};
            \node[anchor=east] at (0,1.9) {\scriptsize Ship 3};
            \node[anchor=east] at (0,2.6) {\scriptsize Ship $C$};

            \draw[draw=orange!80!black, fill=orange!20, ultra thick] (0, 0.3) rectangle (1.8, 0.7);
            \node[font=\tiny, font=\bfseries] at (0.9, 0.5) {$W_1$};

            \draw[draw=orange!80!black, fill=orange!20, ultra thick] (0, 1.0) rectangle (2.5, 1.4);
            \node[font=\tiny, font=\bfseries] at (1.25, 1.2) {$W_2$};

            \draw[draw=orange!80!black, fill=orange!20, ultra thick] (0, 1.7) rectangle (1.1, 2.1);
            \node[font=\tiny, font=\bfseries] at (0.55, 1.9) {$W_3$};

            \draw[draw=orange!80!black, fill=orange!20, ultra thick] (0, 2.4) rectangle (3.1, 2.8);
            \node[font=\tiny, font=\bfseries] at (1.55, 2.6) {$W_C$};
            
            \node[color=orange!80!black, font=\bfseries\footnotesize, anchor=west] at (1.3, 1.9) {Area = $\sum_{i=1}^C W_i = C \times \ex(W)$};
        \end{scope}

    \end{tikzpicture}
    \caption{Geometric derivation of Little's Law by equating total system revenue over an observation interval $T$. View 1 scales the time-weighted average number of ships in the system ($\ex(N) \times T$) by a fee $p$. View 2 scales the spreadsheet accounting record of individual random turnaround times experienced by $C$ arriving ships ($\sum W_i = C \times \ex(W)$) by the same fee $p$. Equating both operational perspectives balances the ledger, yielding $p\ex(N)T = pC\ex(W)$, which simplifies directly to $\ex(N) = \ex(a)\ex(W)$ given the average arrival rate $\ex(a) = C/T$.}
    \label{fig:littles_law_sample_path_perfect}
\end{figure}
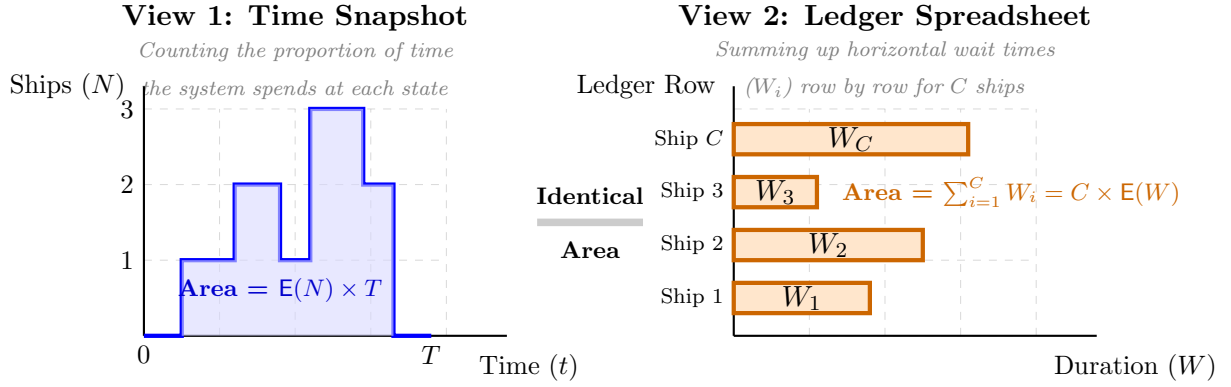

For the last situation, the country notes an interesting fact: as long as the system is flowing, that is, the rate of arrivals $\ex(a)$ is lower than the maximum rate of service $5/\ex(\tau)$. Meaning that the line of ships will not grow indefinitely, it must be the case that the rate of departures is the same as the rate of arrivals. The logic follows from the fact that it cannot be the case that the rate of departures is higher than the rate of arrivals as this would imply that somehow more ships are created during the inspection, which is impossible. On the other hand, the rate of departures cannot be less than that of departures as it would imply a collapse of the system longterm (remember, everything is time-homogenous so far).

If we think about why this is the case, recall that in the derivation of the equilibrium distribution, we had that for every state $\pi_k\lambda_k=\pi_{k+1}\mu_{k+1}$. Which in plain words means that in equilibrium between two states, the inflow and outflow of that partition is exactly the same. This property is often called \emph{detailed balance}\footnote{This property is described as $\pi_kq_{k,k+1}=\pi_{k+1}q_{k+1,k}$} and means that the chain behaves the same forward or backwards (think that there is a perfect symmetry of the flows forwards or backwards). Imagining an alternative process, where time is reversed, that is arrivals are departures and departures are arrivals, will distribute exactly the same. The key is that there is a strong symmetry that allows this.

The conclusion is extremely powerful, because if departures behave the same probabilistically than arrivals then it means that when you connect queues together, you can study them each individually as the output will behave the same as the input (in this case a so-called Poisson process, which only means that the time between arrivals follows an exponential distribution).

If for example the country wished to connect multiple queues in tandem (say $\ell$) or in a network where the next tandem is selected with a fixed probability (a so-called routing matrix), starting with Poisson arrivals at a rate $\ex(a)=e^{-5p}$. Then, the equilibrium distribution of the number in each queue is just the product of the equilibrium distributions of each queue as:
\begin{eqnarray}
  \pi(N_1,\cdots,N_{\ell})=\pi_1(N_1)\times\cdots\times\pi_{\ell}(N_\ell),
\end{eqnarray}
Likewise, this also implies that after the tandem of queues, the departure rate is identical in distribution to the arrivals. Therefore, the revenue of charging ships only once is $p\ex(a)=pe^{-5p}$. The price that maximizes the revenue is given by equating the derivative to 0, that is, $pe^{-5p}5=e^{-5p}$ or $p=1/5$.
\end{exm}

\section{Time-varying arrivals}
Time homogeneity (in both arrivals and processing times) is often an unrealistic assumption. Many manufacturing processes are subject to time-varying demand patterns. For example, ports have different demand patterns depending of the time of the year and time of day. Likewise, certain manufacturing processes see a spike of the demand at certain times of the year. Toy manufacturers see an increase towards the holidays. Certain rare events such as pandemics increase the demand for medical supplies and so on. All of these examples emphasize the need to extend the analysis to systems where the arrival changes over time.

The need to consider time-varying arrivals seem to contradict the time-homogenous structure of the (Continuous) Markov Chains that we have studied. To see why, suppose that the random arrivals $a$ change in average over time. From our equilibrium analysis, we get away by considering time by having that the interarrival times follow an exponential distribution, where the probability of seeing a new arrival in an interval $\Delta t$ is $\lambda\Delta t$ for an exponential random variable with rate $\mu$. If we wanted to change this rate to be time dependent as $\lambda(t)$ it would break the equilibrium conditions everywhere, as $\frac{d P(t)}{dt}$ is not a constant anymore.

Before throwing the time-homogenous time into the garbage, one can imagine a way to salvage the model, while still allowing some degree of change in the arrivals: The strategy is to consider dividing time into relatively \emph{long} intervals where the arrivals remain with a constant rate $r_0$ and then change to another rate $r_1$. For example, this can happen every 12 hours. You can extend the strait example to imagine that there is a higher demand rate $r_0$ during the daytime hours (say 7am-7pm) than during the other 12 night hours.

We can imagine in this mental experiment, the system is going to behave almost like the homogeneous case in the middle of each period. For example, at noon the system pretty much should look as if the arrival rate is exactly $r_1$ all the time. The problem is going to be at the times that the rate changes. But if we are willing to ignore this, then the approximation should still work with some caveats.
For example, Little's Law states that the average number in system $\ex(N)$ is equal to $\ex(W)\ex(a)$. In this case, we should expect that $\ex(N)=\ex(W)r_0$ during the daytime range and $\ex(N)=\ex(W)r_1$ during the night, assuming the service rate is also the same. If the processing is very fast compared to the arrival rate, meaning that the system is empty most of the time. Then, these approximations should hold relatively well. The problem is that if the processing is slow, if $r_0>r_1$, at the beginning of the night shift at 7pm there might be still jobs in process from the earlier part of the day.

This phenomenon is in fact very common in a lot of manufacturing and retail systems: processing is high-traffic, that is the system sees a number of arrivals that it can barely process. In a way, this makes economic sense, as a faster system means spending more money on processing resources (workers, servers, stations, etc), and unless there is a reward for higher throughput, the company will spend the minimum amount of resources that allows to meet the demand or target throughput. In other words, it makes economic sense to have a system that is barely stable. In a way it makes sense that some unprocessed work accumulates over less busy periods.

\begin{figure}[h]
    \centering
    \begin{tikzpicture}
        \begin{groupplot}[
            group style={group size=1 by 2, vertical sep=2.2cm}, 
            width=0.85\textwidth,
            height=4.5cm,
            grid=major,
            xmin=-4, xmax=4,
            xtick={-4, 0, 4},
            xticklabels={3pm, {7pm \\ ($t=0$)}, 11pm},
            xticklabel style={align=center, font=\small},
            ylabel style={font=\small},
            xlabel style={font=\small}
        ]

        \nextgroupplot[
            ylabel={Arrival Rate $\lambda(t)$},
            ymin=0, ymax=12,
            ytick={3, 9},
            yticklabels={$r_1$ (Night), $r_0$ (Day)},
            title={\textbf{Time-Varying Step Demand Profile}}
        ]
        \addplot[blue, ultra thick, domain=-4:0] {9};
        \addplot[blue, ultra thick, domain=0:4] {3};
        \draw[blue, thick, dashed] (axis cs:0,9) -- (axis cs:0,3);

        \nextgroupplot[
            ylabel={Avg. Jobs $\ex[N(t)]$},
            xlabel={Time of Day},
            ymin=0, ymax=45,
            ytick={10, 40}, 
            yticklabels={$\ex[N_{\text{night}}]$, $\ex[N_{\text{day}}]$},
            title={\textbf{System Congestion \& Backlog Hangover ($r_0 > r_1$)}}
        ]
        \addplot[red, ultra thick, domain=-4:0] {40};
        
        \addplot[red, ultra thick, domain=0:4, samples=100] {10 + 30*exp(-0.7*x)};
        
        \draw[gray, dashed, thick] (axis cs:-4,40) -- (axis cs:0,40);
        \draw[gray, dashed, thick] (axis cs:0,10) -- (axis cs:4,10);
        
        \node[anchor=north west, font=\scriptsize, color=black!70] at (axis cs:-3.8, 37.5) 
            {Day Equilibrium: $\ex[N] = \ex[W]r_0$};
            
        \node[anchor=north east, font=\scriptsize, color=black!70] at (axis cs:3.8, 7.5) 
            {Night Equilibrium: $\ex[N] = \ex[W]r_1$};
            
        \node[draw, fill=gray!5, rounded corners, font=\scriptsize, align=center, anchor=west] at (axis cs:2.2, 28) 
            {Backlog Residual\\ $N(0) > 0$ at $t=0$};
            
        \draw[-={Stealth[scale=0.7]}, thin, black!50] (axis cs:2.2, 28) -- (axis cs:0.03, 40.0);

        \end{groupplot}
    \end{tikzpicture}
    \caption{The mechanics of a non-homogeneous queueing system switching operational regimes. The upper axis illustrates a piece-wise constant arrival pattern splitting daytime demand ($r_0$) from nighttime demand ($r_1$). The lower axis visualizes the transient response of the expected inventory $\ex[N(t)]$. If processing speeds are tight, the system creates an inventory carryover at the 7pm boundary ($t=0$), meaning remnants of daytime congestion bleed into the night shift.}
    \label{fig:time_varying_arrivals}
\end{figure}
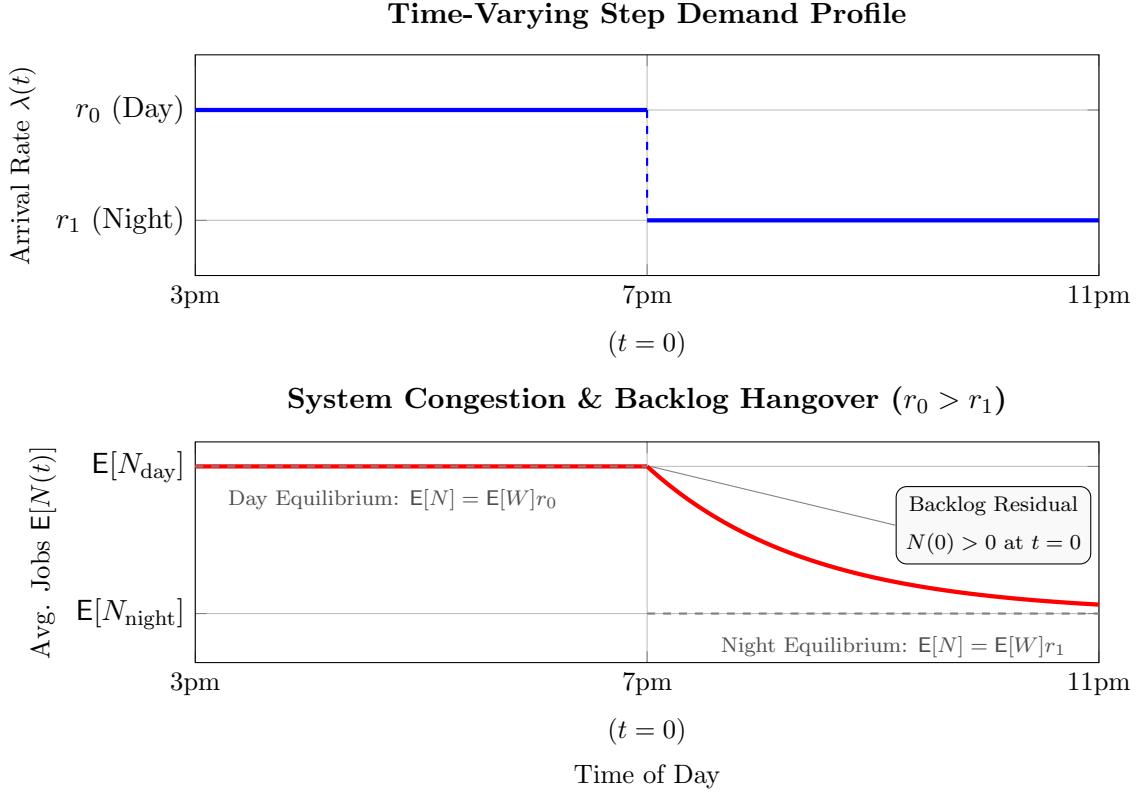

Going back to our mental experiment, at 7am the system should start almost empty. At 7pm, letting $t=0$, suppose that the system has $N(0)>0$ jobs that are still in processing when the arrival rates changes to the lower intensity $r_1$. Then, we would like to know how many jobs are in the system at 8pm, that is, $t=1$. We can do multiple analyses. First we start with the case where there are infinite servers (while unrealistic, there are some systems that behave like this). Suppose, each server servers at a rate $\mu$. See Figure \ref{fig:time_varying_arrivals} for a visual representation.

\subsection{Time-varying arrivals with infinite servers}

Since there are always available servers, we can treat the two populations separately: on one hand the previous jobs $N(0)$ plus the new jobs that arrive between 7pm and 8pm, that is, from $t=0$ to $t=1$. For the new jobs, because of the memoryless property, their time to service is the exponential rate. At $t=1$ we have:
\begin{eqnarray}
  \ex\left[\sum_{i=1}^{N(0)}\ind\{T_i\ge t=1\}\right]=\sum_{i=1}^{N(0)}\pr(T_i\ge 1)=N(0)e^{-1\mu},
\end{eqnarray}

For the new arrivals, between $t=0$ and $t=1$ we have that at each small interval $dt$ there is a new arrival with rate $r_1$ and each of these arrivals will still be in the system at 8 pm, that is, $t=1$ with probability equal to the exponential. For example, an arrival at $s\in(0,1)$ will still be in processing with probability $e^{-(1-s)\mu}$, which is exactly the definition of the reliability function $R(1-s)$. Since at each $dt$ there is an arrival with probability $r_1dt$, and the interval is so small that the probability that 2 arrivals occur is 0. We can do this calculation for each $dt$ as:
\begin{eqnarray}
  \int_0^1 e^{-(1-s)\mu}r_1ds=r_1\int_0^1 e^{-(1-s)\mu}ds=r_1\int_0^1 e^{-s\mu}ds=\frac{r_1}{\mu}(1-e^{-1\mu}),
\end{eqnarray}
Now that we know the expected numbers from each population (the new and old), we can put them together as:
\begin{eqnarray}
  \ex[N(1)|N(0)]=N(0)e^{-1/\ex(W)}+\ex(a(0,1))\ex(W)(1-e^{-1/\ex(W)})
\end{eqnarray}
Where $\ex(W)=\mu^{-1}$ and $\ex(a(0,1))=r_1$. This equation can be seen as a transient version of \emph{Little's Law} where when $\ex(W)\rightarrow 0$, meaning that the system is fast enough to process the jobs quickly and the system just converges to the expected Little's Law (as the term $e^{-1/\ex(W)}$ vanishes) and all that is left is $\ex[N(1)|N(0)]=\ex(a(0,1))\ex(W)$. Since this relation applies between any two time intervals, we can apply it to systems when the arrival rates change (and remain constant) in known intervals.

Then, for a situation when the rates change every period, say $r_1=\ex(a(0,1))$, $r_2=\ex(a(1,2)),\dots$, we can conclude the following relation:
\begin{eqnarray}
  \ex[N(t)|N(t-1)]= N(t-1)e^{-1/\ex(W)}+r_k\ex(W)(1-e^{-1/\ex(W)}),
\end{eqnarray}
If we take expectation on both sides, we end up with the relation:
\begin{eqnarray}
  \ex[N(t)]=\ex[N(t-1)]\alpha+\beta r_t,
\end{eqnarray}
With $\alpha=e^{-1/\ex(W)}$ and $\beta = \ex(W)(1-\alpha)$. Solving this equation recursively starting at $N(0)$, so that $\ex[N(1)]=\alpha N(0)+\beta r_1$. $\ex[N(2)]=\alpha(\alpha N(0)+\beta r_1)+\beta r_2$ and so on, yields $\ex[N(t)]=\alpha^t N(0)+\beta\sum_{i=1}^t\alpha^{i-t} r_i$, or the relation:
\begin{eqnarray}
  \ex[N(t)]=e^{-t/\ex(W)}N(0)+\ex(W)(1-e^{-1/\ex(W)})\sum_{i=1}^t r_ie^{-(i-t)/\ex(W)},
\end{eqnarray}
Intuitively this result can be seen as an average over the entire history of the process. The first term represents the oldest arrivals starting at time 0, and each term of the summation is the contribution of arrivals each time step. Similar as before, when $\ex(W)\rightarrow 0$, all the terms of the summation vanish, except the last one (as $e^{-(t-t)/\ex(W)}=1$), leaving $\ex[N(t)]=\ex(W)r_t$ which is Little's Law as expected.

This kind of recursion is useful when we have a system that has stable service rate, but time varying arrivals. These arrivals can be a seasonal demand that can be estimated with the methods in Chapter \ref{sec:demand} and we want to know how much resources are needed to sustain the operation in the future: imagine a factory with fixed and enough machinery which staffing levels are proportional to the number of jobs in the system $\ex[N(t)]$ for a distant $t$. This recursion would allow us to estimate the staffing load knowing the number of jobs today at $N(0)$.

\subsection{Time-varying arrivals with finite servers}

The infinite server case is a useful approximation in cases where the servers are not the bottleneck of the system, that is, relative to the amount of arrivals there is always enough processing capacity. This is the case in certain manufacturing and computing systems. When the processing capacity is limited by the number of servers, say $k$ of them each serving at a rate $\mu=1/\ex(W)$, the situation is different. Starting the $N(0)>0$ jobs in queue, when $N(0)\le k$ we are in a similar situation to the infinite server case, as the number in system will keep decreasing at an exponential rate towards equilibrium to Little's law. 

Suppose the arrival rate $\ex(a(0,1))$ is again $r_1$ and $N(0)>k$ meaning that the system starts with more jobs than available servers. At maximum capacity the effective service rate is $k\mu$ (all servers are active). In this case, the system will behave linearly, with arrivals growing at a rate of $r_1$ and exiting at a rate $k\mu$. In this case, the case behaves in a linear fashion as the number in system is $(r_1-k\mu)t$. If $r_1>k\mu$ the queue will keep growing as the server at maximum capacity is not able to process jobs fast enough to vacate the queue. Likewise, when $r_1\le k\mu$, the number in system will decrease linearly until either the arrival rate changes at $t=1$ to $r_2$ or the number in the system is equal to $k$. More explicitly, for a time $t\in(0,1)$ the number in system is in average $N(t)=N(0)+(r_1-k\mu)t$.

If for example, $(r_1-k\mu)t$ is decreasing fast enough it could happen that there is a $t<1$ that the system now has less than $k$ jobs, that is, a $N(t)<k$. In this case, as jobs get processed the system becomes exactly like the infinite server case as jobs that arrive will find a server and will behave statistically the same way as the previous section. This time will be $t^*=\frac{k-N(0)}{r_1-k\mu}$ if this is less than 1, then the change in the dynamic will occur by then and the system will start behaving as the infinite server regime. That is, $\ex [N(1)]=\ex[N(t^*)]e^{-(1-t^*)/\ex(W)}+ r_1\ex(W)(1-e^{-(1-t^*)/\ex(W)})$. This can be seen as restarting the system at time $t^*$.

As the rates change the system can alternate between the linear regime $\ex[N(t)]>k$ or the exponential regime $\ex[N(t)]\le k$. Even though it is called exponential regime, it does not mean that it is faster than the linear one. In fact, it is slower as less than $k$ servers are in operation while the arrival rate remains the same. For a system starting at time $t-1$ with $N(t-1)$ jobs there are the following cases:
\begin{itemize}
  \item $N(t-1)>k$ a linear system: The system evolves as $\ex[N(t)]=\ex[N(t-1)]+(r_t-k\mu)t$ either until the next rate change or until a time $\tau>t-1$ when the system enters the exponential regime.
  \item $N(t-1)\le k$ an exponential system: The system evolves as $\ex[N(t)]=\ex[N(t-1)]\alpha(\Delta t)+(1-\alpha(\Delta t))\ex(W) r_t$ until the change of arrival rate or until a time $\tau>t-1$ when the system enters the linear regime.
\end{itemize}

We illustrate this with an example:
\begin{exm}[World Cup Tickets]
  A cloud provider service hosts an application to sell World Cup tickets. They have staggered access to the website starting at 9am to 10 thousand premium customers, the system starts with $N(0)=5$ thousand and the arrival rate for the next hour is the other 5 thousand. The general population will have access starting at 10am, arriving at a rate of 20 thousand per hour until noon when the sale closes. Each person takes in average 10 minutes to complete their purchase. The cloud computing company wants to know how much computing power to allocate. They can choose between two configurations $k=3$ or $k=5$ thousand capacity servers. Plot the average number of customers in the system.

  The first aspect to determine is in what regime the system starts. In the case of the configuration with $k=5$ thousand servers the system starts in the exponential system. As time progresses, the average number in system $\ex[N(t)]$ converges to Little's law: The average time in system is $\ex[W]=1/\mu=1/6$ and the arrival rate is $\ex(a(0,1))=5$ thousand. Then, the system decays until convergence at $\ex(W)\ex(a(0,1))=5/6$ thousand (approximately 833 customers). This will continue until time $t=1$, corresponding to 10am when the arrival rate $\ex(a(1,3))$ changes to $20$ thousand per hour. In this case, since we start at the exponential regime at time 1 with 833 average customers. The new equilibrium implied by Little's law is $\ex(a(1,3))\ex(W)=20/6$ thousand (approximately 3,333 customers). Since this is below $k=5000$, the system should continue in the exponential regime until closing at time 3 (noon).

  The other configuration is slighly more interesting. First, the system starts in the linear regime as $N(0)=5$ thousand is greater than the available $k=3$ thousand servers. The slope of the system is $\ex(a(0,1))-k\mu=5-3(6)=-13$ thousand customers per hour. Starting at $N(0)=5$ the time that the system reach the exponential regime is $5-13t=3$ or $t=2/13\approx 0.15$ or around 9:10am. At this time the system enters the exponential regime and converges again to the equilibrium point of 833 customers. Likewise, at 10am the system starts in the exponential regime towards 3,333 customers. But before that, it must reach 3000 (the number of servers) at a time $e^{-6t}(5/6)+(1-e^{-6t})(20/6)=3$ or just $2/6=e^{-6t}(15/6)$ which simplifies to $t=\ln(15/2)/6\approx 0.33$ or around 10:20am. At this time, the system enters the linear regime growing a rate of $\ex(a(1,3))-k\mu=20-3(6)=2$ thousand customers per hour. This linear growth will continue until the sales closing time at noon, the system will never recover or reach equilibrium. See Figure \ref{fig:worldcup}.

\begin{figure}[h]
    \centering
    \begin{tikzpicture}
        \begin{groupplot}[
            group style={group size=2 by 1, horizontal sep=1.5cm},
            xlabel={Hours (from 9am)},
            grid=major,
            width=0.48\textwidth,
            height=6cm,
            xmin=0, xmax=3,
            ymin=0, ymax=6000,
        ]

        \nextgroupplot[title={\textbf{$k=5000$}}, ylabel={$\ex[N(t)]$}]
        \addplot[blue, thick, domain=0:1] {5000*exp(-6*x) + (5000/6)*(1 - exp(-6*x))};
        \addplot[red, thick, domain=1:3] {
            (5000*exp(-6) + (5000/6)*(1-exp(-6)))*exp(-6*(x-1)) + (20000/6)*(1-exp(-6*(x-1)))
        };
        \node[anchor=west] at (axis cs:0,5300) {\small $N(0)=5000$};

        \nextgroupplot[title={\textbf{$k=3000$}}]
        
        \addplot[orange, thick, domain=0:0.154, forget plot] {5000 - 13000*x};
        
        \addplot[blue, thick, domain=0.154:1] {3000*exp(-6*(x-0.154)) + (5000/6)*(1 - exp(-6*(x-0.154)))};
        
        \addplot[purple, thick, domain=1:1.336] {
            834*exp(-6*(x-1)) + (20000/6)*(1 - exp(-6*(x-1)))
        };
        
        \addplot[red, thick, domain=1.336:3] {3000 + 2000*(x - 1.336)};
        
        \addplot[dashed, black] coordinates {(0,3000) (3,3000)};
        \node[anchor=south east] at (axis cs:3,3000) {\tiny Capacity $k$};

        \end{groupplot}
    \end{tikzpicture}
    \caption{Transient behavior of the average number of customers in the system $\ex[N(t)]$ over a 3-hour window for two capacity configurations. The left plot ($k=5000$) shows a system that remains entirely in the exponential regime, successfully absorbing both the initial workload and the 10am demand surge well below capacity limits. The right plot ($k=3000$) exhibits alternating regimes: an initial linear drain down to the capacity boundary, a stable exponential phase chasing an 833-customer equilibrium, followed by a post-10am surge that crosses back over the capacity line into an unsustainable, linear growth phase.}
    \label{fig:worldcup}
\end{figure}
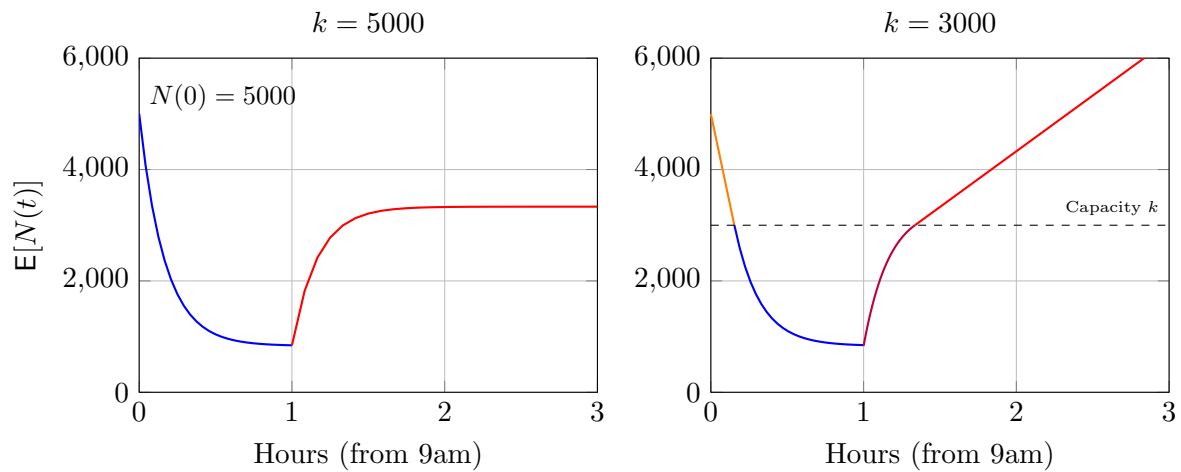

\end{exm}

\chapter{Optimization of Queueing Systems}
\section{Introduction}
We have built a decent understanding of queueing systems by analyzing the variants derived from the memoryless property with changes in the arrival rates and the number of servers. In real life, every queueing system is tied to some resource that needs to be optimized: the staff scheduling of a retail store, the server capacity of an online application, the number of active stations serving a manufacturing process. As mentioned in the previous Chapter, the matching of resources needs to make economic sense in the context of the business goal: satisfiying some target demand within a prescribed time, sustaining a level of service for all incoming customers, etc.

\section{Linear Allocation}
Going back to the staff scheduling models studied in \ref{sec:scheduling} imagine that you want to schedule a facility that serves customers, each taking a random time $\tau$ in service. This setting is typically a call center, a retail store, a restaurant, a factory or any system where when jobs accumulate past certain level they wait in queue. The resource is an integer number representing the number of servers available to serve customers: such as employees. The rate the system serves is $\mu$.

It would be great if there is a system where the service $\ex(W)$, say the waiting time a customer spends in the system: waiting in queue plus their service time behaves linearly. The simplest model is a single server, operating at a rate $\mu$. In this case, we can apply the same logic: a customer arriving sees in average $\ex(N)$ customers, each served in average at a time $1/\mu$ per customer, plus their own service time $1/\mu$. Then, a random customer waits in this system as:
\begin{eqnarray}
  \ex(W)=\ex(N)(1/\mu)+1/\mu,
\end{eqnarray}
But we also know that by Little's law that $\ex(N)=\ex(a)\ex(W)$, replacing this in the equation yields $\mu\ex(W)=\ex(a)\ex(W)+1$ or:
\begin{eqnarray}
  \ex(W)=\frac{1}{\mu-\ex(a)},
\end{eqnarray}
The implicit assumption is that the system is stable, that is, that the arrival rate $\ex(a)$ is less than the processing rate of the server $\mu$. This expression is great because it is almost linear on the server term $\mu$. Then, we could approximate the rate of the server $\mu$ to be proportional to the number of resources, say $x\tilde\mu$, choosing $\tilde \mu$ in such way that it matches the desired processing rate. Where $x$ is a vector choosing the number of resources assigned to a given collection of schedules such as in Chapter \ref{sec:scheduling}, this leads to the following scheduling formulation:

 \begin{eqnarray}
	&\min&f(x)\\
	&\text{s.t. }&\mu Ax^\intercal\ge \ex(\bm A),\nonumber\\
  	&&1\le(\tilde\mu Ax^\intercal -\ex(\bm A))\bar W,\nonumber\\
	&&x\in\calx.\nonumber
\end{eqnarray}
The first constraint just represents the stability condition of the system, ensuring that the demand $\ex(\bm A)$, which is a vector with elements $\ex(a(t-1,t))$ is below the processing capacity of the system. The second constraint represents limiting the average waiting time to some maximum level $\bar W$, such that, at all hours the service is in average below this level. As before, $f(x)$ can be choosen to be the payroll. Or even a more elaborate measure, we illustrate this with an example:
\begin{exm}[Car Wash Scheduling]
  An owner of a car wash wants to minimize their operating costs which are proportional to their payroll, expressed as a vector $c x^\intercal$ and the water and electricity costs of their car wash, which are $\gamma$ dollars per unit of time when each station is on (there are in total $K$ stations). While servicing all the demand, the owner wants to do a scheduling at minimum cost. How to model this problem?

  First, satisfying all the demand is the constraint $\mu Ax^\intercal\ge \ex(\bm A)$, which is the same stability condition of the queueing system. The cost part is more elaborated: the payroll cost is $c x^\intercal$ as in the Scheduling chapter, but the variable cost is in principle a complex calculation as it would require in principle estimating the average behaviour of the system for different staffing configurations $Ax^\intercal$. Then, we want to know how to express $f(x)=cx^\intercal+\text{Variable Cost}$. To do this, we can use what we learned in the Chapter \ref{sec:queueing} via Little's law, we know that the average number in system is given by:
  \begin{eqnarray}
    \ex(N(t))=\ex(a(t-1,t))\ex(W)=\ex(a(t-1,t))\frac{1}{\tilde \mu (Ax^\intercal)_t -\ex(a(t-1,t))},
  \end{eqnarray}
  Then, the cost over the periods the car-wash is open $t=1,2,\dots,T$ is approximately $\sum_{t=1}^T\gamma\ex(N(t))$. Then, the next step is to linearize the variable cost part of the objective $\sum_{t=1}^T\gamma\ex(N(t))$ on $x$ so that we can use an optimizer to solve the problem. A way to write this problem as a MIP is to estimate $\ex(N(t))$ for the range of possible values it can take, for example, let $\ex_k(N(t))=\frac{\ex(a(t-1,t))}{\tilde\mu k-\ex(a(t-1,t))}$ and then creating a variable $z_{k,t}$ selecting the corresponding value to 1 if $(Ax^\intercal)_t=k$. Then, the problem can be written as:

 \begin{eqnarray}
  &\min&cx^\intercal+\gamma\sum_{t=1}^T\ex_k(N(t))z_{k,t}\\
	&\text{s.t. }&\mu Ax^\intercal\ge \ex(\bm A),\nonumber\\
  && (Ax^\intercal)_t=\sum_{k=1}^K k z_{k,t},\,\,\text{for}\,\,t=1,\dots,T,\nonumber\\
  &&\sum_{k=1}^Kz_{k,t}=1\,\,\text{for}\,\,t=1,\dots,T,\nonumber\\
  &&x\in\calx,z_{i,k}\in\{0,1\}.\nonumber
\end{eqnarray}

\end{exm}

Clearly, this type of formulation has its pros and cons: on one hand it is a tractable formulation that in a single pass allows to solve a seemingly complex problem without multiple simulation runs in a difficult combinatorial space. On the other hand, it is still an approximation on many levels.

The first element, is that $\tilde \mu$ should be chosen carefully as one mega-server working always at a rate $k\mu$ is not the same as $k$ servers at a rate $\mu$. To see why, when there is less than $k$ jobs in the system, say $i<k$ the superserver operates at a rate $k\mu$ while the multiple servers operate at a rate $\min(i,k)\mu$. Then, it makes sense to deteriorate the rate $\tilde\mu$ to be slower to compensate for the speedup. A simple way is to define a factor $\tilde\mu_k=\eta_k\mu$ where $\eta_k<1$, thus slowing down the server. One such value for $\eta$ could be the probability that all servers are active $\eta_k=P(i\ge k)=\sum_{i=k}^\infty \pi_i$ in the original $M/M/k$ queue, this slows down the server to reflect that when $i<k$ some servers are idle in the original model, thus making the approximation better reflect the original behavior of $k$ servers into a single one with normalized rate $\mu\eta_k$.

\subsection{Mean-variance optimization}

Another common feature throughout out exposition is modeling the risk-aware version of the problem. In this case, we can exploit the same approximation to consider the variance of the number in system, that is $\var_k(N(t))$ with multiplier $\lambda$ (see Chapter \ref{sec:networks} to see different ways to choose the parameter). This leads to the equivalent risk aware formulation:

 \begin{eqnarray}
  &\min&cx^\intercal+\gamma\sum_{t=1}^T\ex_k(N(t))z_{k,t}+\lambda\sum_{t=1}^T\var_k(N(t))z_{k,t}\\
	&\text{s.t. }&\mu Ax^\intercal\ge \ex(\bm A),\nonumber\\
  && (Ax^\intercal)_t=\sum_{k=1}^K k z_{k,t},\,\,\text{for}\,\,t=1,\dots,T,\nonumber\\
  &&\sum_{k=1}^Kz_{k,t}=1\,\,\text{for}\,\,t=1,\dots,T,\nonumber\\
  &&x\in\calx,z_{i,k}\in\{0,1\}.\nonumber
\end{eqnarray}
Closed expressions for the Variance of queueing systems are widely known, for example, for the $M/M/1$ approximation, we have $\var_k(N(t))=\frac{\ex(a(t-1,t))k\tilde\mu}{[k\tilde\mu - \ex(a(t-1,t))]^2}$.

\section{Accuracy and non-exponential service times}
The approach presented is compelling in the sense that it allows to model service times at determistic epochs $t,t+1,...$ that are useful for planning horizons where the accuracy is tolerable (this is true for many manufacturing systems where noise averages out) and the resources scheduled have a difficult combinatorial nature (for example, scheduling people).

There are other systems that have a mathematically more amenable description of their resource space but require more accuracy. One key example of such system are computer services that deal with thousands of client request simultaneously. Large Language Models (LLMs) providers have to solve a server allocation problem to service their demand within a reasonable time in system for hundreds of thousands of concurrent requests. Here, our understanding of queueing systems shines to solve these kind of problems. Imagine the following problem:

\begin{eqnarray}
  &\min& \ex \left(\int_{0}^T[c(k_t)+\gamma(N_t)]dt\right)\\
  &s.t. & \pr(W>s)\le \alpha.\nonumber
\end{eqnarray}
Where $k_t$ is the number of active servers at time $t$, $N_t$ is the number in system at time $t$ and $c$ and $c(k_t)$ is a cost function dependent on the number of active servers, while $\gamma(N_t)$ is a variable cost function depending on the number of queries in the system. $W$ represents the random amount of time in system the customers face. The problem can be stated as finding the number of servers that in average keeps the system at minimum cost while with a low probability $\alpha$ makes that the time in system does not exceed a threshold $s$. This could easily represent the problem that a cloud provider has to solve to service an application or service (such as providing inference for LLMs). In principle we could apply the same strategy that we proposed in the exposition of the previous section, but in this case we can do better.

Recalling Chapter \ref{sec:scheduling} we can apply a similar dynamic programming principle, were we can solve the problem using the same iterative procedure. Start with a set of actions, the number of servers, that is, $k=0,1,\dots,K$ and discretize the time horizon at intervals of size $\Delta t$. The key as always is to set up a recursive relation for the problem, in this case:

\begin{eqnarray}
  V_t(N)=\min_{k}[c(k,N)\Delta t+V_{t+1}(f(N,k,r_t))],
\end{eqnarray}
Where $r_t:=\ex(a(t,t+\Delta t))$ and $f(N,k,r_t)$ represents the transition function of a system that currently has $N$ customers, $k$ active servers and a expected arrival rate $r_t$. The first and easiest part is defining the loss function $c(k,N)$ which is just $c(k)+\gamma(N)$.

For the dynamics of the system $f(N,k,r_t)$, we can exploit what we know about the transient behaviour of a queueing system. For a small $\Delta t$ we have that in expected value:

\begin{eqnarray}
  N_{next}=f(N,k,r_t) = \begin{cases}
N + (r_t - k\mu)\Delta t & \text{if } N \geq k \quad \text{(Linear/Saturated)} \\
N e^{-\mu \Delta t} + \frac{r_t}{\mu}(1 - e^{-\mu \Delta t}) & \text{if } N < k \quad \text{(Exponential/Recovery)}
\end{cases}
\end{eqnarray}

Similar to the previous chapters, a way to incorporate the service constraint $\pr(W>s)\le\alpha$ is to prune the paths where this occurs. The trick is to express this probability in terms of the primitives of the model.

There are many ways to achieve this, perhaps the simplest is the argument that we used to derive the waiting time of the single server queue. In this case, when the server is saturated $N>k$, the waiting time a new customer faces is the waiting time of people in queue $(N-k)$, each entering service at a rate $k\mu$. And after that, the customer enters their own service at a rate $\mu$. That is, the time the $N$-th customer stays in the system with $k$ active servers is in average:
\begin{eqnarray}
  \ex(W(N))=\frac{(N-k)^+}{k\mu}+\frac{1}{\mu},
\end{eqnarray}
Likewise, its variance is the variance of the exponential times, which are independent from each other, that is:
\begin{eqnarray}
  \var(W(N))=\frac{(N-k)^+}{(k\mu)^2}+\frac{1}{\mu^2},
\end{eqnarray}
Then, we can approximate the random waiting time at each step as:
\begin{eqnarray}
  W(N)\approx\frac{\max(N,k)}{k\mu}+Z\sqrt{\frac{\max(N,k)}{(k\mu)^2}}, 
\end{eqnarray}
Therefore, the probability $\pr(W(N)>s)\le \alpha$ or $\pr(W(N)\le s)\ge 1-\alpha$, is given by the expression:
\begin{eqnarray}
  sk\mu\ge \max(N,k)+\sqrt{\max(N,k)}\Phi^{-1}(1-\alpha),
\end{eqnarray}
Then, combinations of $N,k$ that violate this condition are not admissible in the DP iteration.
\subsection{Non-exponential service times}
Throughout our exposition we have considered exponential service times. The exponential service times have allowed us to use the memoryless property and analyze the system at each snapshot of time. Considering non-exponential service times tends to make control problems intractable as there is a need to keep track of either the elapsed or remaining service times as they are no longer memory-less. The reason is that the reliability function in general does not follow time-homogeneity, that is, $R(s)R(t)\neq R(s+t)$.

In recent years, a new framework called QPLEX (see \cite{dieker2025qplex}) allows to estimate transition dynamics without doing full simulations or having to keep track of an explosive state-space. The key insight is that given a specific number of customers in the system we sample a single customer in service remaining time (say $N_t=N$, that in the previous subsection we were taking as the expected value). Conditional on $N$ the remaining time in service can be seen as an i.i.d. random variable. Call $L$ as the random remaining service time of a job picked at random. Estimating the distribution of $L$ for a given time $t$ is what QPLEX does through a series of clever tricks. Call this distribution $\nu_t$. With this distribution as input, conditioning on a starting number of jobs $N_t$ we can express most of the primitives of the problem. For example, the number of expected departures is $N_t\nu_t(L=\Delta t|N_t)$. In fact, QPLEX estimates the distributions of the number in system at the next step $p_t$ and the remaining service durations $\nu_t$, call these pair of distributions $p_t,\nu_t$.

QPLEX is able to approximate these sequential distributions through clever conditionings. Using the distributions $p_t,\nu_t$ we can redo our dynamic programming procedure as:
\begin{eqnarray}
  V_t(N)=\min_{k}[\ex_{p_t}(c(k,N))\Delta t+\ex_{\nu_t}(V_{t+1}|N,k,r_t)].
\end{eqnarray}
\printbibliography
\end{document}